\documentclass[11pt,a4paper]{article}

\usepackage{amsmath,amssymb}
\newcommand{\e}{\varepsilon}
\newcommand{\D}{\Delta}
\newcommand{\La}{\Lambda}
\newcommand{\n}{\nabla}
\newcommand{\N}{\frac{N}{2}}
\newcommand{\q}{q_{l}}
\newcommand{\ui}{u_{l}}
\newcommand{\g}{\int_{\mathbb{R}^{N}}}
\newcommand{\p}{\partial}

\newcommand{\R}{\mathbb{R}}
\newcommand{\h}{\hookrightarrow}
\newcommand{\w}{\widetilde{T}}
\newcommand{\de}{\delta}

\newtheorem{definition}{Definition}
\newtheorem{theorem}{Theorem}
\newtheorem{proposition}{Proposition}

\newtheorem{remarka}{Remark}

\newtheorem{lemme}{Lemma}

\title{Existence of solutions for compressible fluid models of Korteweg type }
\author{Boris Haspot \thanks{Universit\'e Paris XII - Val de Marne 61, avenue
du G\'en\'eral de Gaulle 94 010 CRETEIL Cedex T\'el\'ephone : (33-1)
45 17 16 51 T\'el\'ecopie : (33-1) 45 17 16 49 e-mail :
haspot@univ-paris12.fr}}
\date{}

\begin{document}

\maketitle
\subsubsection*{Abstract}
This work is devoted to the
study of the initial boundary value problem for a general non
isothermal model of capillary fluids derived by J.E Dunn and
J.Serrin (1985) in \cite{fDS,fR}, which can be used as a phase transition model.\\
We distinguish two cases, when the physical coefficients depend only
on the density, and the general case. In the first case we can work
in critical scaling spaces, and we prove global existence of
solution and uniqueness for data close to a stable equilibrium. For
general data, existence and uniqueness is stated on a short time
interval.\\
In the general case with physical coefficients depending on density
and on temperature, additional regularity is required to control the
temperature in $L^{\infty}$ norm. We prove global existence of
solution close to a stable equilibrium and local in time existence
of solution with more general data. Uniqueness is also obtained.
\section{Introduction}
\subsection{Derivation of the Korteweg model}
We are concerned with compressible fluids endowed with internal
capillarity. The model we consider  originates from the XIXth
century work by van der Waals and Korteweg \cite{fK} and was
actually derived in its modern form in the 1980s using the second
gradient theory, see for instance \cite{fJL,fTN}.
\\
Korteweg-type models are based on an extended version of
nonequilibrium thermodynamics, which assumes that the energy of the
fluid not only depends on standard variables but on the gradient of
the density. Let us consider a fluid of density $\rho\geq0$,
velocity field $u\in\R^{N}$ ($N\geq2$), entropy density $e$, and
temperature $\theta=(\frac{\p e}{\p s})_{\rho}$. We note $w=\n\rho$,
and we suppose that the intern specific energy, $e$ depends on the
density $\rho$, on the entropy specific $s$, and on $w$. In terms of
the free energy, this principle takes the form of a generalized
Gibbs relation :
$$de=\widetilde{T}ds+\frac{p}{\rho^{2}}d\rho+\frac{1}{\rho}\phi^{*}\cdot dw$$
where $\w$ is the temperature, $p$ the pressure,
$\phi$ a vector column of $\R^{N}$ and $\phi^{*}$ the adjoint vector.\\
In the same way we can write a differential equation for the intern
energy per unit volume, $E=\rho e$,
$$dE=\w dS+gd\rho+\phi^{*}\cdot dw$$
where $S=\rho s$ is the entropy per unit volume and
$g=e-s\w+\frac{p}{\rho}$ is the chemical
potential.\\
In terms of the free energy, the Gibbs principle gives us:
$$dF=-Sd\w+gd\rho+\phi^{*}\cdot dw.$$
In the present chapter, we shall make the hypothesis that:
$$\phi=\kappa w.$$
The nonnegative coefficient $\kappa$ is called the capillarity and
may depend on both $\rho$ and $\w$.
\\
All the thermodynamic quantities are sum of their classic version
(it
means independent of $w$) and of one term in $|w|^{2}$.\\
In this case the free energy F decomposes into a standard part
$F_{0}$ and an additional term due to gradients of density:
$$F=F_{0}+\frac{1}{2}\kappa|w|^{2}.$$
We denote $v=\frac{1}{\rho}$ the specific volume and $k=v\kappa$.
Similar decompositions hold for $S$, $p$ and $g$:
$$
\begin{aligned}
&p=p_{0}-\frac{1}{2}K_{p}|w|^{2}\;\;\;\;\mbox{where:}\;\;
K_{p}=k^{'}_{v}\;\mbox{and}\;p_{0}=-(f_{0})^{'}_{v}\\
&g=g_{0}+\frac{1}{2}K_{g}|w|^{2}\;\;\;\;\mbox{where:}\;\;
K_{e}=k-\w k^{'}_{\w}\;\mbox{and}\;e_{0}=f_{0}-\w(f_{0})^{'}_{\w}.\\
\end{aligned}
$$
The model deriving from a Cahn-Hilliard like free energy (see the
pioneering work by J.E.Dunn and J.Serrin in \cite{fDS} and also in
\cite{fA,fC,fGP}), the conservation of mass, momentum and energy read:
$$
\begin{cases}
\begin{aligned}
&\p_{t}\rho+{\rm div}(\rho u)=0,\\
& \p_{t}(\rho u)+{\rm div}(\rho
u\otimes u+pI)={\rm div}(K+D)+\rho f,\\
&\p_{t}(\rho(e+\textstyle{\frac{1}{2}}u^{2}))+{\rm div}( u(\rho
e+\textstyle{\frac{1}{2}}\rho|u|^{2}+p))={\rm
div}((D+K)\cdot u-Q+W)+\rho f\cdot u,\\
\end{aligned}
\end{cases}
$$
with:
$$
\begin{aligned}
&D=(\lambda{\rm div}u)I+\mu(du+\n u),\;\mbox{is the diffusion tensor}\\[2mm]
&K=(\rho\,{\rm div}\phi)I-\phi w^{*},\;\mbox{is the Korteweg tensor}\\[2mm]
&Q=-\eta\n \w,\;\mbox{is the heat flux.}\\
\end{aligned}
$$
The term
$$W=(\p_{t}\rho+u^{*}\cdot\n\rho)\phi=-(\rho{\rm
div}u)\phi$$
is the intersticial work which is needed in order to ensure the entropy balance and
was first introduced by Dunn and Serrin in \cite{fDS}.\\
The coefficients $(\lambda,\mu)$ represent the
viscosity of the fluid and may depend on both the density $\rho$ and
the temperature $\w$. The thermal coefficient $\eta$ is a given non
negative function of
the temperature $\w$ and of the density $\rho$.\\
Differentiating formally the equation of conservation of the mass,
we obtain a law of conservation for $w$:
$$\p_{t}w+{\rm div}(uw^{*}+\rho du)=0\;.$$
One may obtain an equation for $e$ by using the  mass and
momentum conservation laws and the
relations:
$${\rm div}((-pI+K+D)u)=({\rm div}(-pI+K+D))\cdot u-p{\rm div}(u)+(K+D):\n u\;.$$
Multiplying the momentum equation by $u$ yields:
$$\big({\rm div}(-pI+K+D)\big)\cdot u=\big(\p_{t}(\rho u)+{\rm div}(\rho
uu^{*})\big)\cdot u=\p_{t}(\frac{\rho|u|^{2}}{2})+{\rm
div}(\frac{\rho|u|^{2}}{2}u)\;.$$ We obtain then:
$$\rho(\p_{t}e+u^{*}\cdot\n e)+p\,{\rm div}u=(K+D):\n u+{\rm div}(W-Q)\;.$$
In substituting $K$, we have (with the summation convention over
repeated indices):
$$K:\n u=\rho\,{\rm div}\phi\,{\rm div}u-\phi_{i}w_{j}\p_{i}u_{j}\;,$$
while:
$$-{\rm div}W={\rm div}((\rho{\rm div}u)\phi)=\rho({\rm div}\phi)({\rm
div}u)+(w^{*}\cdot\phi){\rm
div}u+\phi_{i}\,\rho\,\p^{2}_{j,i}u_{j}.$$ In using
$w_{j}=\p_{j}\rho$, we obtain:
$$
\begin{aligned}
K:\n u-{\rm div}W&=-(w^{*}\cdot\phi){\rm
div}u-\phi_{i}\p_{j}(\rho\p_{i}u_{j})\\
&=-(w^{*}\cdot\phi){\rm div}u-({\rm div}(\rho du))\cdot\phi.\\
\end{aligned}
$$
Finally, the equation for $e$ rewrites:
$$\rho(\p_{t}e+u^{*}\cdot\n e)+(p+w^{*}\cdot\phi){\rm div}u=D:\n u-({\rm
div}(\rho du))\cdot\phi-{\rm div}Q\;.$$ From now on, we shall
denote: $d_{t}=\p_{t}+u^{*}\cdot\n$.
\subsection{The case of a generalized
Van der Waals law }
From now on, we assume that there exist two functions
$\Pi_{0}$ and $ \Pi_{1}$ such that:
$$p_{0}=\w \Pi^{'}_{1}(v)+\Pi^{'}_{0}(v),$$
$$e_{0}=-\Pi_{0}(v)+\varphi(\w)-\w\varphi^{'}(\w).$$
We now suppose that the coefficients $\lambda$, $\mu$ depend on the
density and on the temperature, and in all the sequel the
capillarity $\kappa$  doesn't depend on the
temperature.\\
Moreover we suppose that the intern specific energy is an increasing
function of $\w$:
$$\Psi^{'}(\w)>0\;\mbox{with}\;\Psi(\w)=\varphi(\w)-\w\varphi^{'}(\w).\leqno{(A)}$$
We then set  $\theta=\Psi(\w)$ and we search to obtain an equation
on $\theta$.
In what follows, we assume that $\kappa$ depends only on the
specific volume.
\subsubsection*{Obtaining an equation for $\theta$ :} As:
$$e=-\Pi_{0}(v)+\theta+\frac{1}{2}\kappa|w|^{2},$$
we thus have:
$$d_{t}e=-\Pi_{0}^{'}(v)d_{t}v+d_{t}\theta+\frac{1}{2}\kappa^{'}_{v}|w|^{2}d_{t}v+\kappa w^{*}\cdot d_{t}w\;.$$
By a direct calculus we find:
$$d_{t}v=v{\rm div}u\;\;\;\text{and}\;\;\;w^{*}\cdot d_{t}w=-|w|^{2}{\rm div}u-{\rm
div}(\rho\,du)\cdot w\;.$$ Then we have:
$$d_{t}\theta=d_{t}e+v(p-\w\Pi^{'}_{1}(v)){\rm
div}u+\kappa|w|^{2}{\rm div}u+\kappa{\rm div}(\rho du)\cdot w.$$ And
in using
the third equation of the system, we get an equation on $\theta$:\\
$$d_{t}\theta+v{\rm div} Q+v\w\Pi^{'}_{1}(v){\rm div}
u=vD:\n u+{\rm div}(\rho du)\cdot(\kappa
w-v\phi)+(\kappa|w|^{2}-vw^{*}\cdot\phi){\rm div} u\;.$$ But as we
have $\phi=\kappa w$ and $k=v\kappa$ we conclude that:
$$d_{t}\theta-v{\rm
div}(\chi\n\theta)+v\Psi^{-1}(\theta)\Pi^{'}_{1}(v){\rm div} u=vD:\n
u$$ with: $\chi(\rho,\theta)=\eta
(\rho,\w)(\Psi^{-1})^{'}(\theta)\;.$
\subsubsection*{Obtaining a system for $\rho,\,u,\,\theta$:}
We obtain then for the momentum equation:
$$d_{t}u-\frac{{\rm
div}D}{\rho}+\frac{\n p_{0}}{\rho}=\frac{\rm div
K}{\rho}+\frac{1}{2}\frac{\n(K_{p}|w|^{2})}{\rho}$$
where $K_{p}=\kappa-\rho \kappa^{'}_{\rho}$.\\
And by a calculus we check that:
$${\rm
div}K+\frac{1}{2}\n(K_{p}|w|^{2})=\rho\n(\kappa\D\rho)+\frac{\rho}{2}\n(\kappa_{\rho}^{'}|\n\rho|^{2}).$$
Indeed we have:
$$
\begin{aligned}
I=&\n(\rho{\rm div}(\kappa\n\rho))-{\rm
div}(\kappa ww^{*})+\frac{1}{2}\n(K_{p}|w|^{2})\\[3mm]
=&[\rho\n(\kappa\D\rho)+\kappa\n\rho\D\rho+\rho\n(\kappa^{'}_{\rho}|\rho|^{2})+\kappa^{'}_{\rho}|\n\rho|^{2}\n\rho]-
[\kappa{\rm
div}(ww^{*})+\kappa^{'}_{\rho}\,w\cdot\n\rho\, w],\\[3mm]
&+[\frac{\kappa}{2}\n|w|^{2}+\frac{\kappa^{'}_{\rho}}{2}\n|w|^{2}\n\rho-\frac{1}{2}\n(\rho
\kappa^{'}_{\rho}|w|^{2})],\\[3mm]
=&[\rho\n(\kappa\D\rho)+\kappa\n\rho\D\rho+\frac{\rho}{2}\n(\kappa^{'}_{\rho}|\rho|^{2})]-[\kappa{\rm
div}(ww^{*})]+[\frac{\kappa}{2}\n|w|^{2}],\\[3mm]
=&[\rho\n(\kappa\D\rho)+\kappa\n\rho\D\rho+\frac{\rho}{2}\n(\kappa^{'}_{\rho}|\rho|^{2})]-[\kappa
w{\rm
div} w+\frac{\kappa}{2}\n|w|^{2}]+[\frac{\kappa}{2}\n|w|^{2}],\\[3mm]
=&\rho\n(\kappa\D\rho)+\frac{\rho}{2}\n(\kappa_{\rho}^{'}|\n\rho|^{2})\;.\\
\end{aligned}
$$
Finally we have obtained the following system:
$$
\begin{cases}
\begin{aligned}
&\p_{t}\rho+{\rm div}(\rho u)=0,\\[3mm]
&\p_{t}u+u\cdot\n u-\frac{{\rm
div}D}{\rho}-\nabla(\kappa\Delta\rho)+\frac{\n(P_{0}(\rho)+\Psi^{-1}(\theta)P_{1}(\rho))}{\rho}\\
&\hspace{7cm}=\n(\frac{\kappa^{'}_{\rho}}{2}|\n\rho|^{2}),\\[3mm]
&\p_{t}\theta+u\cdot\n\theta-\frac{{\rm
div}(\chi\n\theta)}{\rho}+\Psi^{-1}(\theta)\frac{P_{1}(\rho)}{\rho}{\rm
div}(u)=\frac{D:\n u}{\rho},\\
\end{aligned}
\end{cases}
\leqno{(NHV)}
$$
where: $P_{0}=\Pi^{'}_{0}$, $P_{1}=\Pi^{'}_{1}$ and
$\w=\Psi^{-1}(\theta)$.\\
We supplement $(NHV)$ with initial conditions:
$$\rho_{t=0}=\rho_{0}\geq 0\,\,\, u_{t=0}=u_{0},\,\,\mbox{and}\,\,\,\theta_{t=0}=\theta_{0}.$$
\subsection{Classical a priori-estimates}
Before getting into the heart of mathematical results, let us
derive the physical energy bounds of the $(NHV)$ system when
$\kappa$ is a constant and where the pressure just depends on the
density to simplify. Let $\bar{\rho}>0$ be a constant reference
density, and let $\Pi $ be defined
by:
$$\Pi(s)=s\biggl(\int^{s}_{\bar{\rho}}\frac{P_{0}(z)}{z^{2}}dz-\frac{P_{0}(\bar{\rho})}{\bar{\rho}}\biggl)$$
so that $P_{0}(s)=s\Pi^{'}(s)-\Pi(s)\, ,\,\Pi^{'}(\bar{\rho})=0$
and:
$$\p_{t}\Pi(\rho)+{\rm div}(u\Pi(\rho))+P_{0}(\rho){\rm div}(u)=0\;\;\;\mbox{in}\;\;
{\cal D}^{'}((0,T)\times\R^{N}).$$
Notice that $\Pi$ is convex as far as $P$ is non decreasing (since
$P_{0}^{'}(s)=s\Pi^{''}(s)$), which is the case for $\gamma$-type
pressure laws or for Van der Waals law above the critical temperature.\\
Multiplying the equation of momentum conservation in the system
$(NHV)$ by $\rho u$ and integrating by parts over $\R^{N}$, we
obtain the following
estimate:
$$
\begin{aligned}
&\int_{\R^{N}}\big(\frac{1}{2}\rho
|u|^{2}+\rho\theta+(\Pi(\rho)-\Pi(\bar{\rho}))+\frac{\kappa}{2}|\nabla\rho|^{2}\big)(t)dx
+2\int_{0}^{t}\int_{\R^{N}}\big(2\mu
D(u):D(u)\\
&\hspace{2cm}+(\lambda+\mu)|{\rm div} u|^{2}\big)dx
\leq\int_{\R^{N}}\big(\frac{|m_{0}|^{2}}{2\rho}+\rho_{0}\theta_{0}+(\Pi(\rho_{0})-\Pi(\bar{\rho}))
+\frac{\kappa}{2}|\nabla\rho_{0}|^{2}\big)dx.\\
\end{aligned}
$$
It follows that assuming that the initial total energy is finite:\\
$$\epsilon_{0}=\int_{\R^{N}}\big(\frac{|m_{0}|^{2}}{2\rho}+\rho_{0}\theta_{0}+(\Pi(\rho_{0})-\Pi(\bar{\rho}))
+\frac{\kappa}{2}|\nabla\rho_{0}|^{2}\big)dx<+\infty\,,$$ then we
have the a priori bounds:
$$\Pi(\rho)-\Pi(\bar{\rho}),\;\rho |u|^{2},\;\;\mbox{and}\;\;\rho\,\theta\in
L^{\infty}(0,\infty,L^{1}(\R^{N}))$$
$$\n\rho\in L^{\infty}(0,\infty,L^{2}(\R^{N}))^{N},\;\;\mbox{and}\;\;\n u\in
L^{2}(0,\infty,\R^{N})^{N^{2}}.$$
\section{Mathematical results}
We wish to prove existence and uniqueness results for $(NHV)$ in
functions spaces very close to energy spaces. In the non isothermal
non capillary case and $P(\rho)=a\rho^{\gamma}$, with $a>0$ and
$\gamma>1$, P-L. Lions in \cite{fL2} proved the global existence of
variational solutions $(\rho,u,\theta)$ to $(NHV)$ with $\kappa=0$
for $\gamma> \frac{N}{2}$ if $N\geq 4$, $\gamma\geq \frac{3N}{N+2}$
if $N=2,3$ and initial data $(\rho_{0},m_{0})$ such that:
$$\Pi(\rho_{0})-\Pi(\bar{\rho}),\;\;\frac{|m_{0}|^{2}}{\rho_{0}}\in
L^{1}(\R^{N}),\;\;\mbox{and}\;\;\rho_{0}\theta_{0}\in
L^{1}(\R^{N}).$$ These solutions are weak solutions in the classical
sense for the equation of mass conservation and for
the equation of the momentum.\\
On the other hand, the weak solution satisfies only an inequality
for the thermal energy equation.
\\
Notice that the main difficulty for proving Lions' theorem consists
in exhibiting strong compactness properties of the density $\rho$ in
$L^{p}_{loc}$ spaces required to pass to the limit in the pressure
term $P(\rho)=a\rho^{\gamma}$.\\
Let us mention that Feireisl in \cite{fF}  generalized the result to
$\gamma>\frac{N}{2}$ in establishing that we can obtain renormalized
solution without imposing that $\rho\in L^{2}_{loc}$, for this he
introduces the concept of oscillation defect measure evaluating the
lost of compactness.\\
We can finally cite a very interesting result from Bresch-Desjardins
in \cite{fBD},\cite{fBD1} where they show the existence of global weak
solution for $(NHV)$ with $\kappa=0$ in choosing specific type of
viscosity where $\mu$ and $\lambda$ are linked. It allows them to
get good estimate on the density in using energy inequality and to
can treat by compactness all the delicate terms. This result is very
new because the energy equation is verified really in distribution
sense. In \cite{fMV}, Mellet and Vasseur improve the results of
Bresch,Desjardins in generalize to some coefficient $\mu$ and
$\lambda$ admitting the vacuum in the case of Navier-Stokes
isothermal, they use essentially a gain of integrability on the
velocity.
\\
\\
In the case $\kappa>0$, we remark then that the density belongs to
$L^{\infty}(0,\infty,{\dot {H}}^{1}(\R^{N}))$. Hence, in contrast to
the non capillary case one can easily pass to the limit in the
pressure term. However let us emphasize at this point that the above
a priori bounds do not provide any $L^{\infty}$ control on the
density from below or from above. Indeed, even in dimension $N=2$,
$H^{1}$ functions are not necessarily locally bounded. Thus, vacuum
patches are likely
 to form in the fluid in spite of the presence of capillary forces,
which are expected to smooth out the density. Danchin and Desjardins
show in \cite{fDD} that the isothermal
 model has weak solutions if there exists $c_{1}$ and $M_{1}$ such that:
$$c_{1}\leq|\rho|\leq M_{1}\;\;\mbox{and}\;\;|\rho-1|<<1.$$
The vacuum is one of the main difficulties to get weak solutions,
and the problem remains open.
\\
In the isothermal capillary case with specific type of viscosity and
capillarity $\mu(\rho)=\mu\rho$ and $\lambda(\rho)=0$, Bresch,
Desjardins and Lin in \cite{fBDL} obtain the global existence of weak
solutions without smallness assumption on the data. We can precise
 the space of test functions used depends on the solution itself
which are on the form $\rho\phi$ with $\phi\in
C^{\infty}_{0}(\R^{N})$. The specificity of the viscosity allows to
get a gain of one derivative on the density: $\rho\in L^{2}(H^{2})$.
\\
\\
Existence of strong solution with $\kappa$, $\mu$ and $\lambda$
constant is known since the work by Hattori an Li in \cite{fH1},
\cite{fH2} in the whole space $\R^{N}$. In \cite{fDD}, Danchin and
Desjardins study the well-posedness of the problem for the
isothermal case with constant coefficients in critical Besov spaces.
\\
Here we want to investigate the well-posedness of the full non
isothermal problem in critical spaces, that is, in spaces which are
invariant by the scaling of Korteweg's system. Recall that such an
approach is now classical for incompressible Navier-Stokes equation
and yields local well-posedness (or global well-posedness for small
data) in spaces with minimal regularity.
\\
Let us explain precisely the scaling of Korteweg's system. We can
easily check that, if $(\rho,u,\theta)$ solves $(NHV)$, so does
$(\rho_{\lambda},u_{\lambda},\theta_{\lambda})$, where:
\\
$$\rho_{\lambda}(t,x)=\rho(\lambda^{2}t,\lambda x)\,,\,u_{\lambda}(t,x)=\lambda u(\lambda^{2}t,\lambda x)\,\,\,\,\mbox{and}\,\,\,\, \theta_{\lambda}(t,x)=\lambda^{2} \theta(\lambda^{2}t,\lambda x)$$
\\
provided the pressure laws $P_{0},P_{1}$ have been changed into
$\lambda^{2}P_{0},\lambda^{2}P_{1}.$
\begin{definition}
We say that a functional space is critical with respect to the
scaling of the equation if the associated norm is invariant under
the transformation:
$$(\rho,u,\theta)\longrightarrow(\rho_{\lambda},u_{\lambda},\theta_{\lambda})$$
(up to a constant independent of $\lambda$).
\end{definition}
This suggests us to choose initial data
$(\rho_{0},u_{0},\theta_{0})$ in spaces whose norm is invariant for
all $\lambda>0$ by
$(\rho_{0},u_{0},\theta_{0})\longrightarrow(\rho_{0}(\lambda\cdot),\lambda
u_{0}(\lambda\cdot),\lambda^{2}\theta_{0}(\lambda\cdot)).$
\\
\\
A natural candidate is the homogeneous Sobolev space
$\dot{H}^{N/2}\times (\dot{H}^{N/2-1})^{N}\times\dot{H}^{N/2-2}$,
but since $\dot{H}^{N/2}$ is not included in $L^{\infty}$, we cannot
expect to get $L^{\infty}$ control on the density when
$\rho_{0}\in\dot{H}^{N/2}$. The same problem occurs in the equation
for the temperature when dealing with the non linear term  involving
$\Psi^{-1}(\theta)$.
\\
This is the reason why, instead of the classical homogeneous Sobolev
space $\dot{H}^{s}(\R^{d})$, we will consider homogeneous Besov
spaces with the same derivative index
$B^{s}=\dot{B}^{s}_{2,1}(\R^{N})$ (for the corresponding definition
we refer to section 4).
\\
One of the nice property of $B^{s}$ spaces for critical exponent $s$
is that $B^{N/2}$ is an algebra embedded in $L^{\infty}$. This
allows to control the density from below and from above, without
requiring more regularity on derivatives of $\rho$. For similar
reasons, we shall take $\theta_{0}$ in $B^{\N}$ in the general case
where appear non-linear terms in function of the temperature.
\\
Since a global in time approach does not seem to be accessible for
general data, we will mainly consider the global well-posedness
problem for initial data close enough to stable equilibria (Section
5). This motivates the following definition:\\
\begin{definition}
\label{fbar} Let $\bar{\rho}>0$, $\bar{\theta}>0$. We will note in
the sequel:
$$q=\frac{\rho-\bar{\rho}}{\bar{\rho}}\;\;\mbox{and}\;\;{\cal T}=\theta-\bar{\theta}.$$
\end{definition}
One can now state the main results of the paper.\\
The first three theorems concern the global existence and uniqueness
of solution to the Korteweg's system with {\it small} initial data.
In particulary the first two results concern  Korteweg's system with
coefficients depending only on the density and where the intern
specific energy is a linear function of the temperature.
\begin{theorem} \label{ftheo1}Let $N\geq3$. Assume that the function $\Psi$ defined in $(A)$ satisfies
$\Psi(\w)=A\w$ with $A>0$ and that all the physical coefficients are
smooth functions depending only on the density.
Let $\bar{\rho}>0$ be such that:
$$\kappa(\bar{\rho})>0,\;\;\;\mu(\bar{\rho})>0,\;\;\;\lambda(\bar{\rho})+2\mu(\bar{\rho})>0,\;\;\eta(\bar{\rho})>0\;\;
\mbox{and}\;\;\p_{\rho}P_{0}(\bar{\rho})>0.$$
Moreover suppose that:
$$q_{0}\in\widetilde{B}^{\N-1,\N},\;u_{0}\in B^{\N-1},\;{\cal T}_{0}\in\widetilde{B}^{\N-1,\N-2}.$$
There exists  an $\e_{0}$ depending only on the  physical
coefficients
(that we will precise later) such that if:
$$\|q_{0}\|_{\widetilde{B}^{\frac{N}{2}-1,\N}}+\|u_{0}\|_{\widetilde{B}^{\N-1}}+
\|{\cal{T}}_{0}\|_{\widetilde{B}^{\N-1,\N-2}}\leq\e_{0}$$
then (NHV) has a unique global solution $(\rho,u,{\cal T})$ in
$E^{N/2}$ where  $E^{s}$ is defined by:
$$
\begin{aligned}
E^{s}=&[C_{b}(\R_{+},\widetilde{B}^{s-1,s})\cap
L^{1}(\R_{+},\widetilde{B}^{s+1,s+2})]\times[C_{b}(\R_{+},B^{s-1})^{N}\cap
L^{1}(\R_{+},B^{s+1})^{N}]\\
&\times[C_{b}(\R_{+},\widetilde{B}^{s-1,s-2})\cap
L^{1}(\R_{+},\widetilde{B}^{s+1,s})].\\
\end{aligned}
$$
\end{theorem}
\begin{remarka}Above, $\widetilde{B}^{s,t}$ stands for a Besov space with
regularity $B^{s}$ in low frequencies and $B^{t}$ in high frequencies (see definition
\ref{fBesov}).
\end{remarka}
The case $N=2$ requires more regular initial data because of
technical problems involving some nonlinear terms in the temperature
equation.
\begin{theorem} \label{ftheo2}Let $N=2$. Under the assumption of the theorem
\ref{ftheo1} for $\Psi$ and the physical coefficients, let $\e^{'}>0$
and suppose that:
$$q_{0}\in\widetilde{B}^{0,1+\e^{'}},\;u_{0}\in \widetilde{B}^{0,\e^{'}},\;{\cal T}_{0}\in\widetilde{B}^{0,-1+\e^{'}}.$$
There exists  an $\e_{0}$ depending only on the physical
coefficients
such that if:
$$\|q_{0}\|_{\widetilde{B}^{0,1+\e^{'}}}+\|u_{0}\|_{\widetilde{B}^{0,\e^{'}}}
+\|{\cal{T}}_{0}\|_{\widetilde{B}^{0,-1+\e^{'}}}\leq\e_{0}$$
then (NHV) has a unique global solution $(\rho,u,{\cal T})$ in
the space:
$$
\begin{aligned}
E^{'}=&[C_{b}(\R_{+},\widetilde{B}^{0,1+\e^{'}})\cap
L^{1}(\R_{+},\widetilde{B}^{2,3+\e^{'}})]\times[C_{b}(\R_{+},\widetilde{B}^{0,\e^{'}})^{2}\cap
L^{1}(\R_{+},\widetilde{B}^{2,2+\e^{'}})^{2}]\\
&\times[C_{b}(\R_{+},\widetilde{B}^{0,-1+\e^{'}})\cap
L^{1}(\R_{+},\widetilde{B}^{2,1+\e^{'}})].\\
\end{aligned}
$$
\end{theorem}
In the following theorem we are interested by showing the global
existence of solution for  Korteweg's system with general conditions
and small initial data. In order to control the non linear terms in
temperature more regularity is required. That's why we want control
the temperature in norm $L^{\infty}$.
\begin{theorem}
\label{ftheo3}Let $N\geq2$. Assume that $\Psi$ be a regular function
depending on $\theta$. Assume that all the coefficients are smooth
functions of $\rho$ and $\theta$ except $\kappa$ which depends only
on the density. Take $(\bar{\rho},\bar{T})$ such that:
$$\kappa(\bar{\rho})>0,\;\;\;\mu(\bar{\rho},\bar{T})>0,\;\;\;\lambda(\bar{\rho},\bar{T})+2\mu(\bar{\rho},\bar{T})>0,
\;\;\eta(\bar{\rho},\bar{T})>0\;\;\mbox{and}\;\;\p_{\rho}P_{0}(\bar{\rho},\bar{T})>0.$$
Moreover suppose that:
$$q_{0}\in\widetilde{B}^{\N-1,\N+1},\;u_{0}\in \widetilde{B}^{\N-1,\N},\;{\cal{T}}_{0}\in\widetilde{B}^{\N-1,\N}.$$
There exists  an $\e_{1}$ depending only on the physical
coefficients
such that if:
$$\|q_{0}\|_{\widetilde{B}^{\N-1,\N+1}}+\|u_{0}\|_{\widetilde{B}^{\N-1,\N}}+\|{\cal{T}}_{0}\|_{\widetilde{B}^{\N-1,\N}}
\leq\e$$ then (NHV) has a  unique global solution $(\rho,u,{\cal
T})$ in:
$$
\begin{aligned}
F^{\N}=&[C_{b}(\R_{+},\widetilde{B}^{\N-1,\N+1})\cap
L^{1}(\R_{+},\widetilde{B}^{\N+1,\N+3})]
\times[C_{b}(\R_{+},\widetilde{B}^{\N-1,\N})^{N}\\
&\cap L^{1}(\R_{+},\widetilde{B}^{\N+1,\N+2})^{N}]
[C_{b}(\R_{+},\widetilde{B}^{\N-1,\N})\cap
L^{1}(\R_{+},\widetilde{B}^{\N+1,\N+2})].\\
\end{aligned}
$$
\end{theorem}
In the previous theorem we can observe for the case $N=2$ that the
initial data are very close from the energy space of Bresch,
Desjardins and Lin in \cite{fBDL}.\\
In the following three theorems we are interested by the existence
and uniqueness of solution in finite time for {\it large} data. We
distinguish always the differents cases $N\geq3$ and $N=2$ if the
coefficients depend only on $\rho$, and the case where the
coefficients depend also on $\tilde{T}$.
\begin{theorem}
\label{ftheo4}
\smallskip
Let $N\geq3$,and  $\Psi$ and the physical coefficients be as in
theorem \ref{ftheo1}. We suppose that $(q_{0},u_{0},{\cal{T}}_{0})\in
B^{\N}\times (B^{\N-1})^{N}\times B^{\N-2}$ and that $\rho_{0}\geq
c$ for some $c>0$.
\\[1,5mm]
Then there exists a time $T$ such that system (NHV) has a unique
solution in $F_{T}$
$$
\begin{aligned}
F_{T}=&[\widetilde{C}_{T}(B^{\N})\cap
L^{1}_{T}(B^{\N+2})]\times[\widetilde{C}_{T}(B^{\N-1})^{N}\cap
L^{1}_{T}(B^{\N+1})^{N}]\\
&\times[\widetilde{C}_{T}(B^{\N-2})\cap L^{1}_{T}(B^{\N})]\;.\\
\end{aligned}
$$
\end{theorem}
For the same reasons as previously in the case $N=2$ we can not reach the critical level of regularity.
\begin{theorem}
\label{ftheo5}
\smallskip
Let $N=2$ and $\e^{'}>0$. Under the assumptions of theorem
\ref{ftheo1} for $\Psi$ and the physical coefficients we suppose that
$(q_{0},u_{0},{\cal{T}}_{0})\in \widetilde{B}^{1,1+\e^{'}}\times
(\widetilde{B}^{0,\e^{'}})^{2}\times \widetilde{B}^{-1,-1+\e^{'}}$
and $\rho_{0}\geq c$ for some $c>0$.
\\
Then there exists a time $T$ such that the system has a unique
solution in $F_{T}(2)$ with:
$$
\begin{aligned}
F_{T}(2)=&[\widetilde{C}_{T}(\widetilde{B}^{1,1+\e^{'}})\cap
L^{1}_{T}(\widetilde{B}^{3,3+\e^{'}})]\times[\widetilde{C}_{T}(\widetilde{B}^{0,\e^{'}})^{2}\cap
L^{1}_{T}(\widetilde{B}^{2,2+\e^{'}})^{2}]\\
&\times[\widetilde{C}_{T}(\widetilde{B}^{-1,-1+\e^{'}})\cap L^{1}_{T}(\widetilde{B}^{1,1+\e^{'}})]\;.\\
\end{aligned}
$$
\end{theorem}
In the last theorem we see the general system without conditions,
and like previously we need more regular initial data.
\begin{theorem}
\label{ftheo6}Under the hypotheses of theorem \ref{ftheo3} we suppose
that:
$$(q_{0},u_{0},{\cal{T}}_{0})\in
\widetilde{B}^{\N,\N+1}\times (B^{\N})^{N}\times B^{\N}\;\;\mbox{
and}\;\; \rho_{0}\geq c\;\;\mbox{for some}\;\;c>0.$$
Then there exists a time $T$ such that the system has a unique
solution in :
$$
\begin{aligned}
F^{'}_{T}=&[\widetilde{C}_{T}(\widetilde{B}^{\N,\N+1})\cap
L^{1}_{T}(\widetilde{B}^{\N+2,\N+3})]\times[\widetilde{C}_{T}(B^{\N})^{N}\cap
L^{1}_{T}(B^{\N+2})^{N}]\\
&\times[\widetilde{C}_{T}(B^{\N})\cap L^{1}_{T}(B^{\N+2})]\;.\\
\end{aligned}
$$
\end{theorem}
This chapter is structured in the following way, first of all we recall
in the section \ref{fsection3} some definitions on Besov spaces and
some useful theorem concerning Besov spaces. Next we will
concentrate in the section \ref{fsection4} on the global existence
and uniqueness of solution for our system $(NHV)$ with small initial
data. In subsection \ref{fsubsection41} we will give some
necessary conditions to get the stability of the linear part
associated to the system $(NHV)$. In subsection
\ref{fsubsection42} we will study the case where the
specific intern energy is linear and where the physical coefficients
are independent of the temperature. In our proof we will distinguish
the case $N\geq3$ and the case $N=2$ for some technical reasons. In
the section \ref{fsection5} we will examine the local existence and
uniqueness of solution with general initial data. For the same reasons as section \ref{fsection4} we will
distinguish the cases in function of the behavior of the
coefficients and of the intern specific energy.
\section{Littlewood-Paley theory and Besov spaces\label{fsection3}}
\subsection{Littlewood-Paley decomposition}
Littlewood-Paley decomposition  corresponds to a dyadic
decomposition  of the space in Fourier variables.\\
We can use for instance any $\varphi\in C^{\infty}(\R^{N})$,
supported in
${\cal{C}}=\{\xi\in\R^{N}/\frac{3}{4}\leq|\xi|\leq\frac{8}{3}\}$
such that:
$$\sum_{l\in\mathbb{Z}}\varphi(2^{-l}\xi)=1\,\,\,\,\mbox{if}\,\,\,\,\xi\ne 0.$$
Denoting $h={\cal{F}}^{-1}\varphi$, we then define the dyadic
blocks by:
$$\D_{l}u=2^{lN}\int_{\R^{N}}h(2^{l}y)u(x-y)dy\,\,\,\,\mbox{and}\,\,\,S_{l}u=\sum_{k\leq
l-1}\D_{k}u\,.$$ Formally, one can write that:
$$u=\sum_{k\in\mathbb{Z}}\D_{k}u\,.$$
This decomposition is called homogeneous Littlewood-Paley
decomposition. Let us observe that the above formal equality does
not hold in ${\cal{S}}^{'}(\R^{N})$ for two reasons:
\begin{enumerate}
\item The right hand-side does not necessarily converge in
${\cal{S}}^{'}(\R^{N})$.
\item Even if it does, the equality is not
always true in ${\cal{S}}^{'}(\R^{N})$ (consider the case of the
polynomials).
\end{enumerate}
However, this equality holds true modulo
polynomials hence homogeneous Besov spaces will be defined modulo
the polynomials, according to \cite{fBG}.
\subsection{Homogeneous Besov spaces and first properties}
\begin{definition}
\label{fBesov} For
$s\in\R,\,\,\mbox{and}\,\,u\in{\cal{S}}^{'}(\R^{N})$
we set:\\
$$\|u\|_{B^{s}}=(\sum_{l\in\mathbb{Z}}(2^{ls}\|\D_{l}u\|_{L^{2}}).$$
\end{definition}
A difficulty due to the choice of homogeneous spaces arises at this
point. Indeed, $ \|.\|_{B^{s}}$ cannot be a norm on
$\{u\in{\cal{S}}^{'}(\R^{N}),\|u\|_{B^{s}}<+\infty\}$ because
$\|u\|_{B^{s}}=0$ means that $u$ is a polynomial. This enforces us
to adopt the following definition for homogeneous Besov spaces, see
\cite{fBG}.
\begin{definition}
Let $s\in\R$.\\
Denote $m=[s-\N]$ if $s-\N\notin\mathbb{Z}$ and
$m=s-\N-1$ otherwise.
\begin{itemize}
\item If $m<0$, then we define $B^{s}$ as:
$$B^{s}=\biggl\{u\in{\cal{S}}^{'}(\R^{N})\;\big/\|u\|_{B^{s}}<\infty\,\,\,\mbox{and}\,\,\,u=
\sum_{l\in\mathbb{Z}}\D_{l}u\,\,\mbox{in}\,\,{\cal{S}}^{'}(\R^{N})\biggl\}\,.$$
\item If $m\geq 0$, we denote by ${\cal{P}}_{m}[\R^{N}]$ the set of
polynomials of degree less than or equal to $m$ and we set:
$$B^{s}=\biggl\{u\in{\cal{S}}^{'}(\R^{N})/{\cal{P}}_{m}[\R^{N}]\;\big/\|u\|_{B^{s}}<\infty\,\,\,and
\,\,\,u=\sum_{l\in\mathbb{Z}}\D_{l}u\,\,in\,\,{\cal{S}}^{'}(\R^{N})/{\cal{P}}_{m}[\R^{N}]\biggl\}\,.$$
\end{itemize}
\end{definition}
\begin{proposition}
\label{fderivation} The following properties hold:
\begin{enumerate}
\item Density: If  $|s|\leq \frac{N}{2}$ , then
$C^{\infty}_{0}$ is dense in $B^{s}$.
\item Derivatives: There exists a universal constant $C$
such that:
$$C^{-1}\|u\|_{B^{s}}\leq\|\n u\|_{B^{s-1}}\leq
C\|u\|_{B^{s}}.$$
\item Algebraic
properties: For $s>0$, $B^{s}\cap L^{\infty}$ is an algebra.
\item Interpolation:
 $(B^{s_{1}},B^{s_{2}})_{\theta,1}=B^{\theta
s_{2}+(1-\theta)s_{1}}$.
\end{enumerate}
\end{proposition}
\subsection{Hybrid Besov spaces and Chemin-Lerner spaces}
Hybrid  Besov spaces are functional spaces where regularity
assumptions are different in low frequency and high frequency, see
\cite{fDG}.
\\
They may be defined as follows:
\begin{definition}
Let $s,\,t\in\R$.We set:
$$\|u\|_{\widetilde{B}^{s,t}}=\sum_{q\leq 0}2^{qs}\|\D_{q}u\|_{L^{2}}+\sum_{q>
0}2^{qt}\|\D_{q}u\|_{L^{2}}\,.$$
Let $m=-[\N+1-s]$, we then define:
\begin{itemize}
\item $\widetilde{B}^{s,t}=\{u\in{\cal{S}}^{'}(\R^{N})\;\big/\|u\|_{\widetilde{B}^{s,t}}<+\infty\}$,
if $m<0$
\item $\widetilde{B}^{s,t}=\{u\in{\cal{S}}^{'}(\R^{N})/{\cal{P}}_{m}[\R^{N}]\;\big/\|u\|_{\widetilde{B}^{s,t}}<+\infty\}$
if $m\geq 0$.
\end{itemize}
\end{definition}
Let us now give some properties of these hybrid spaces and some results
on how they behave with respect to the product. The following
results come directly from the paradifferential calculus.
\begin{proposition}
\label{fimportant}
We recall some inclusion:
\begin{itemize}
\item We have $\widetilde{B}^{s,s}=B^{s}$.
\item If $s\leq t$ then $\widetilde{B}^{s,t}=B^{s}\cap B^{t}$,
if $s>t$ then $\widetilde{B}^{s,t}=B^{s}+B^{t}$.
\item If $s_{1}\leq s_{2}$ and $t_{1}\geq t_{2}$ then
$\widetilde{B}^{s_{1},t_{1}}\hookrightarrow
\widetilde{B}^{s_{2},t_{2}}$.
\end{itemize}
\end{proposition}
\begin{proposition}
For all $s,\,t>0$, we have:
$$\|uv\|_{\widetilde{B}^{s,t}}\leq
C(\|u\|_{L^{\infty}}\|v\|_{\widetilde{B}^{s,t}}+\|v\|_{L^{\infty}}\|u\|_{\widetilde{B}^{s,t}})\,.$$
For all $s_{1},\,s_{2},\,t_{1},\,t_{2}\leq \N$ such that
$\min(s_{1}+s_{2},\,t_{1}+t_{2})>0$ we have:
$$\|uv\|_{\widetilde{B}^{s_{1}+t_{1}-\N,s_{2}+t_{2}-\N}}\leq
C\|u\|_{\widetilde{B}^{s_{1},t_{1}}}\|v\|_{\widetilde{B}^{s_{2},t_{2}}}\,.$$
\end{proposition}
For a proof of this proposition see \cite{fDG}. We are now going to
define the spaces of Chemin-Lerner in which we will work, which are
a refinement of the spaces:
$$L_{T}^{\rho}(B^{s}):=L^{\rho}(0,T,B^{s}).$$
$\hspace{15cm}$
\begin{definition}
Let $\rho\in[1,+\infty[$, $T\in[1,+\infty]$ and $s\in\R$. We then
denote:
$$\|u\|_{\widetilde{L}^{\rho}_{T}(\widetilde{B}^{s_{1},s_{2}})}=
\sum_{l\leq0}2^{ls_{1}}\big(\int_{0}^{T}\|\D_{l}u(t)\|^{\rho}_{L^{2}}dt\big)^{1/\rho}+
\sum_{l>0}2^{ls_{2}}\big(\int_{0}^{T}\|\D_{l}u(t)\|^{\rho}_{L^{2}}dt\big)^{1/\rho}\,.$$
And we have in the case $\rho=\infty$:
$$\|u\|_{\widetilde{L}^{\infty}_{T}(\widetilde{B}^{s_{1},s_{2}})}=
\sum_{l\leq0}2^{ls_{1}}\|\D_{l}u\|_{L^{\infty}(L^{2})}+
\sum_{l>0}2^{ls_{2}}\|\D_{l}u\|_{L^{\infty}(L^{2})}\,.$$
\end{definition}
We note that thanks to Minkowsky inequality we have:
$$\|u\|_{L^{\rho}_{T}(\widetilde{B}^{s_{1},s_{2}})}\leq\|u\|_{\widetilde{L}^{\rho}_{T}
(\widetilde{B}^{s_{1},s_{2}})}\;\;\;\mbox{and}\;\;\;\|u\|_{L^{1}_{T}(\widetilde{B}^{s_{1},s_{2}})}=
\|u\|_{\widetilde{L}^{1}_{T}
(\widetilde{B}^{s_{1},s_{2}})}.$$
From now on, we will denote:
$$\|u\|^{-}_{\widetilde{L}^{\rho}_{T}(B^{s_{1}})}=\sum_{l\leq0}2^{ls_{1}}\big(\int_{0}^{T}\|\D_{l}u(t)\|^{\rho}
_{L^{p}}dt\big)^{1/\rho}$$
$$\|u\|^{+}_{\widetilde{L}^{\rho}_{T}(B^{s_{2}})}=
\sum_{l>0}2^{ls_{2}}\big(\int_{0}^{T}\|\D_{l}u(t)\|^{\rho}_{L^{p}}dt\big)^{1/\rho}.$$
Hence:
$$\|u\|_{\widetilde{L}^{\rho}_{T}(B^{s_{1}})}=\|u\|^{-}_{\widetilde{L}^{\rho}_{T}(B^{s_{1}})}
+\|u\|^{+}_{\widetilde{L}^{\rho}_{T}(\widetilde{B}^{s_{1},s_{2}})}$$
We then define the space:
$$\widetilde{L}^{\rho}_{T}(\widetilde{B}^{s_{1},s_{2}})=\{u\in
L^{\rho}_{T}(\widetilde{B}^{s_{1},s_{2}})/\|u\|_{\widetilde{L}^{\rho}_{T}(\widetilde{B}^{s_{1},s_{2}})}<\infty\}\,.$$
We
denote moreover by $\widetilde{C}_{T}(\widetilde{B}^{s_{1},s_{2}})$
the set of those functions of
$\widetilde{L}^{\infty}_{T}(\widetilde{B}^{s_{1},s_{2}})$ which are
continuous from $[0,T]$ to $\widetilde{B}^{s_{1},s_{2}}$. In the
sequel we are going to give some properties of this spaces
concerning the interpolation and their relationship with the heat
equation.
\begin{proposition}
\label{finterp} Let $s,\,t,\,s_{1},\,s_{2}\in\R$,
$\rho,\,\rho_{1},\,\rho_{2}\in[1,+\infty]$. We have:
\begin{enumerate}
\item Interpolation:
$$
\begin{aligned}
\|u\|_{\widetilde{L}^{\rho}_{T}(\widetilde{B}^{s,t})}\leq\|u\|_{\widetilde{L}^{\rho_{1}}_{T}
(\widetilde{B}^{s_{1},t_{1}})}^{\theta}
\|u\|_{\widetilde{L}^{\rho_{2}}_{T}(\widetilde{B}^{s_{2},t_{2}})}^{1-\theta}&\;\;\;\mbox{with}\;\;\theta\in[0,1]\;\;
\mbox{and}\;\;
\frac{1}{\rho}=\frac{\theta}{\rho_{1}}+\frac{1-\theta}{\rho_{2}},\\[2mm]
& s=\theta s_{1}+(1-\theta)s_{2},\; t=\theta
t_{1}+(1-\theta)t_{2}.\\
\end{aligned}
$$
\item Embedding:
$$\widetilde{L}^{\rho}_{T}(\widetilde{B}^{s,t})\hookrightarrow
L^{\rho}_{T}(C_{0})\;\;\mbox{and}\;\;
\widetilde{C}_{T}(B^{\frac{N}{2}})\hookrightarrow
C([0,T]\times\R^{N}).$$
\end{enumerate}
\end{proposition}
The $\widetilde{L}^{\rho}_{T}(B^{s}_{p})$ spaces suit particulary
well to the
study of smoothing properties of the heat equation. In \cite{fCT}, J-Y. Chemin proved the following proposition:\\
\begin{proposition}
\label{fchaleur} Let $p\in[1,+\infty]$ and
$1\leq\rho_{2}\leq\rho_{1}\leq+\infty$. Let u be a solution of:
$$
\begin{cases}
\begin{aligned}
&\p_{t}u-\mu\D u=f\\
&u_{t=0}=u_{0}\,.\\
\end{aligned}
\end{cases}
$$
Then there exists $C>0$ depending only on $N,\mu,\rho_{1}$ and
$\rho_{2}$ such that:
$$\|u\|_{\widetilde{L}^{\rho_{1}}_{T}(B^{s+2/\rho_{1}})}\leq
 C\|u_{0}\|_{B^{s}}+C\|f\|_{\widetilde{L}^{\rho_{2}}_{T}(B^{s-2+2/\rho_{2}})}\,.$$
\end{proposition}
To finish with, we explain how the product of functions behaves in
the spaces of Chemin-Lerner. We have the following properties:
\begin{proposition}
\label{fproduit} $\hspace{15cm}$
Let $s>0,\;t>0$,
$1/\rho_{2}+1/\rho_{3}=1/\rho_{1}+1/\rho_{4}=1/\rho\leq1$,
$u\in\widetilde{L}^{\rho_{3}}_{T}(\widetilde{B}^{s,t})\cap
\widetilde{L}^{\rho_{1}}_{T}(L^{\infty})$ and
$v\in\widetilde{L}^{\rho_{4}}_{T}(\widetilde{B}^{s,t})\cap
\widetilde{L}^{\rho_{2}}_{T}(L^{\infty})$.
\\[1,5mm]
Then $uv\in\widetilde{L}^{\rho}_{T}(\widetilde{B}^{s,t})$ and we
have:
$$\|uv\|_{\widetilde{L}^{\rho}_{T}(\widetilde{B}^{s,t})}\leq
C\|u\|_{\widetilde{L}^{\rho_{1}}_{T}(L^{\infty})}\|v\|_{\widetilde{L}^{\rho_{4}}_{T}(\widetilde{B}^{s,t})}
+\|v\|_{\widetilde{L}^{\rho_{2}}_{T}(L^{\infty})}\|u\|_{\widetilde{L}^{\rho_{3}}_{T}(\widetilde{B}^{s,t})}.$$
If $s_{1},s_{2},t_{1},t_{2}\leq \N$ , $s_{1}+s_{2}>0$,
$t_{1}+t_{2}>0$,
$\frac{1}{\rho_{1}}+\frac{1}{\rho_{2}}=\frac{1}{\rho}\leq1$,
$u\in\widetilde{L}^{\rho_{1}}_{T}(\widetilde{B}^{s_{1},t_{1}})$ and
$v\in\widetilde{L}^{\rho_{2}}_{T}(\widetilde{B}^{s_{2},t_{2}})$ then
$uv\in\widetilde{L}^{\rho}_{T}(\widetilde{B}^{s_{1}+s_{2}-\N,t_{1}+t_{2}-\N})$
and:
\\
$$\|uv\|_{\widetilde{L}^{\rho}_{T}(\widetilde{B}^{s_{1}+s_{2}-\frac{N}{2},t_{1}+t_{2}-\frac{N}{2}}_{2}),}\leq
C\|u\|_{\widetilde{L}^{\rho_{1}}_{T}(\widetilde{B}^{s_{1},t_{1}})}\|v\|_{\widetilde{L}^{\rho_{2}}_{T}(
\widetilde{B}^{s_{2},t_{2}})}\,.$$
\end{proposition}
For a proof of this proposition see \cite{fDG}. Finally we need an
estimate on the composition of functions in the spaces
$\widetilde{L}^{\rho}_{T}(\widetilde{B}^{s}_{p})$ (see the proof in
the appendix).
\begin{proposition}
\label{fcomposition} Let $s>0$ , $p\in[1,+\infty]$ and
$u_{1},u_{2}\cdots, u_{d}\in\widetilde{L}^{\rho}_{T}(B^{s}_{p})\cap
L^{\infty}_{T}(L^{\infty})$.
\\
\\
(i) Let $F\in W^{[s]+2,\infty}_{loc}(\R^{N})$ such that $F(0)=0$.
Then
$F(u_{1},u_{2},\cdots,u_{d})\in\widetilde{L}^{\rho}_{T}(B^{s}_{p})$.
\\
More precisely, there exists a constant $C$ depending only on $s, p,
N\;\mbox{and}\;F$ such that:
\\
$$
\begin{aligned}
\|F(u_{1},u_{2},\cdots
,u_{d})\|_{\widetilde{L}^{\rho}_{T}(B^{s}_{p})}\leq&
\,\,C(\|u_{1}\|_{L^{\infty}_{T}(L^{\infty})},\|u_{2}\|_{L^{\infty}_{T}(L^{\infty})},\cdots,
\|u_{d}\|_{L^{\infty}_{T}(L^{\infty})})\\[2mm]
&\hspace{2,3cm}(\|u_{1}\|_{\widetilde{L}^{\rho}_{T}(B^{s}_{p})}+\cdots+\|u_{d}\|_{\widetilde{L}^{\rho}_{T}(B^{s}_{p})}).\\
\end{aligned}
$$
\\
\\
(ii) Let
$u\in\widetilde{L}^{\rho}_{T}(\widetilde{B}^{s_{1},s_{2}})$,
$s_{1},\,s_{2}>0$ then we have $F(u)
\in\widetilde{L}^{\rho}_{T}(\widetilde{B}^{s_{1},s_{2}})$ and
$$\|F(u)\|_{\widetilde{L}^{\rho}_{T}(\widetilde{B}^{s_{1},s_{2}})}\leq C(\|u\|_{L^{\infty}_{T}(L^{\infty})})
\|u\|_{\widetilde{L}^{\rho}_{T}(\widetilde{B}^{s_{1},s_{2}})}.$$
\\
(iii) If $v,\,u\in\widetilde{L}^{\rho}_{T}(B^{s}_{p})\cap
L^{\infty}_{T}(L^{\infty})$ and $G\in
W^{[s]+3,\infty}_{loc}(\R^{N})$ then $G(u)-G(v)$ belongs to
$\widetilde{L}^{\rho}_{T}(B^{s}_{p})$ and there exists a constant C
depending only of $s, p ,N\;\mbox{and}\;G$ such that:
\\
$$
\begin{aligned}
\|G(u)-G(v)\|_{\widetilde{L}^{\rho}_{T}(B^{s}_{p})}\leq&
\,\,C(\|u\|_{L^{\infty}_{T}(L^{\infty})},\|v\|_{L^{\infty}_{T}(L^{\infty})})
\big(\|v-u\|_{\widetilde{L}^{\rho}_{T}(B^{s}_{p})}
(1+\|u\|_{L^{\infty}_{T}(L^{\infty})}\\[2mm]
&+\|v\|_{L^{\infty}_{T}(L^{\infty})})+\|v-u\|_{L^{\infty}_{T}(L^{\infty})}(\|u\|_{\widetilde{L}^{\rho}_{T}(B^{s}_{p})}
+\|v\|_{\widetilde{L}^{\rho}_{T}(B^{s}_{p})})\big).\\
\end{aligned}
$$
\\
\\
(iv) If $v,\,u\in\widetilde{L}^{\rho}_{T}(B^{s_{1},s_{2}}_{p})\cap
L^{\infty}_{T}(L^{\infty})$ and $G\in
W^{[s]+3,\infty}_{loc}(\R^{N})$ then $G(u)-G(v)$ belongs to
$\widetilde{L}^{\rho}_{T}(B^{s_{1},s_{2}}_{p})$ and it exists a
constant C depending only of $s, p ,N\;\mbox{and}\;G$ such that:
\\
$$
\begin{aligned}
\|G(u)-G(v)\|_{\widetilde{L}^{\rho}_{T}(\widetilde{B}^{s_{1},s_{2}}_{p})}\leq&
\,\,C(\|u\|_{L^{\infty}_{T}(L^{\infty})},\|v\|_{L^{\infty}_{T}(L^{\infty})})
\big(\|v-u\|_{\widetilde{L}^{\rho}_{T}(\widetilde{B}^{s_{1},s_{2}}_{p})}
(1+\|u\|_{L^{\infty}_{T}(L^{\infty})}\\[2mm]
&+\|v\|_{L^{\infty}_{T}(L^{\infty})})
+\|v-u\|_{L^{\infty}_{T}(L^{\infty})}(\|u\|_{\widetilde{L}^{\rho}_{T}(\widetilde{B}^{s_{1},s_{2}}_{p})}+
\|v\|_{\widetilde{L}^{\rho}_{T}(\widetilde{B}^{s_{1},s_{2}}_{p})})\big).\\
\end{aligned}
$$
\end{proposition}
The proof is an adaptation of a theorem by J.Y. Chemin and H.
Bahouri in \cite{fBC}, see the proof in the Appendix.
\section{Existence of solutions for small initial data\label{fsection4}}
\subsection{Study of the linear part\label{fsubsection41}}
This section is devoted to the study of the linearization of system
$(NHV)$ in  order to get conditions for the existence of solution.
We recall the system $(NHV)$ in the case where $\kappa$ depends only
on the density $\rho$:
$$\begin{cases}
\begin{aligned}
&\p_{t}\rho+{\rm div}(\rho u)=0,\\[2mm]
&\p_{t}(\rho\,u)+{\rm div}(\rho u\otimes u)-{\rm div}D-\rho\nabla(\kappa\Delta\rho)+(\nabla(P_{0}(\rho)+TP_{1}(\rho))\\
&\hspace{9cm}=\rho\nabla(\frac{\kappa^{'}_{\rho}}{2}|\nabla\rho|^{2}),\\
&\p_{t}\theta-\frac{{\rm div}(\chi\nabla\theta)}{\rho}+T\frac{P_{1}(\rho)}{\rho}{\rm div}(u)=\frac{D:\nabla u}{\rho}.\\[2mm]
\end{aligned}
\end{cases}
\leqno{(NHV)}$$
Moreover we have:
$$
\begin{aligned}
&{\rm div}(D)=(\lambda+\mu)\nabla {\rm div} u+\mu\D
u+\n(\lambda){\rm div}
u+(du+\n u)\n\mu,\\[3mm]
&\hspace{1,2cm}=(2\lambda+\mu)\nabla {\rm div} u+\mu\D
u+\p_{1}\lambda(\rho,\theta)\n\rho{\rm div}
u+\p_{2}\lambda(\rho,\theta)\n\theta{\rm
div} u\\[2mm]
&\hspace{3,5cm}+(du+\n u)\p_{1}\mu(\rho,\theta)\n\rho+(du+\n
u)\p_{2}\mu(\rho,\theta)\n\theta.\\
\end{aligned}
$$
 We transform the system to study it  in the neighborhood of
$(\bar{\rho},0,\bar{\theta})$. Using the notation of definition
\ref{fbar}, we obtain the following system where $F,G, H$ contain the
non linear part:
$$
\begin{cases}
\begin{aligned}
&\p_{t}q+{\rm div}u=F\,,\\
&\p_{t}u-\frac{\bar{\mu}}{\bar{\rho}}\Delta
u-\frac{(\bar{\lambda}+\bar{\mu})}{\bar{\rho}}\nabla {\rm
div}\,u-\bar{\rho}\bar{\kappa}\nabla\Delta
q+(P_{0}^{'}(\bar{\rho})+\bar{T}P_{1}^{'}(\bar{\rho}))\nabla
q\\
&\hspace{8cm}+\frac{P_{1}(\bar{\rho})}{\bar{\rho}\psi^{'}(\bar{T})}\nabla{\cal{T}}=G,\\
&\p_{t}{\cal{T}}-\frac{\bar{\chi}}{\bar{\rho}}\Delta{\cal{T}}+\frac{\bar{T}P_{1}(\bar{\rho})}{\bar{\rho}}{\rm
div}u§=§
H.\\
\end{aligned}
\end{cases}
\leqno{(M)}
$$
This induces us to study the following linear system:
$$
\begin{cases}
\begin{aligned}
&\p_{t}q+{\rm div} u=F\\
&\p_{t}u-\widetilde{\mu}\D u-(\widetilde{\mu}+\widetilde{\lambda})-\e\n\D q-\beta\n q-\gamma\n{\cal{T}}=G\\
&\p_{t}{\cal{T}}-\alpha\Delta{\cal{T}}+\delta{\rm div}u=H\\
&\p_{t}u-\widetilde{\mu}\D u=PG\\
\end{aligned}
\end{cases}
\leqno{(M^{'})}
$$
where $\nu,\,\e,\,\alpha,\,\beta,\,\gamma,\,\delta$ and
$\widetilde{\mu}$ are given real parameters. Note that system $(M)$
with right hand side considered as source terms enters in the class
of models $(M^{'})$, it is only a matter of setting:
$$\widetilde{\mu}=\frac{\bar{\mu}}{\bar{\rho}},\;\widetilde{\lambda}=\frac{\bar{\lambda}}{\bar{\rho}},\;\e=\bar{\rho}\,\bar{\kappa},\;\beta=P_{0}^{'}(\bar{\rho})+\bar{T}P_{1}^{'}(\bar{\rho}),\;\gamma=\frac{P_{1}(\bar{\rho})}{\bar{\rho}\psi^{'}(\bar{T})},\;
\alpha=\frac{\bar{\chi}}{\bar{\rho}},\;\delta=\frac{\bar{T}P_{1}(\bar{\rho})}{\bar{\rho}}\;
.$$
We transform
the system in setting:
$$d=\Lambda^{-1}{\rm div} \,u\;\;\mbox{and}\;\;\Omega=\Lambda^{-1}{\rm curl}\, u$$
where we set: $\Lambda^{s}h={\cal{F}}^{-1}(|\xi|^{s}\hat{h})$ (the $\rm curl$ is defined in the appendix).\\
We finally obtain the following system in projecting on divergence
free vector fields and on potential vector fields:
$$
\begin{cases}
\begin{aligned}
&\p_{t}q+\Lambda d=F,\\
& \p_{t}d-\nu\Delta d-\e \Lambda^{3}q-\beta \Lambda
q-\gamma\Lambda{\cal{T}}=\Lambda^{-1}{\rm div}\,
G,\\[2mm]
&\p_{t}{\cal{T}}-\alpha\Delta{\cal{T}}+\delta\Lambda d=H,\\[2mm]
&\p_{t}\Omega-\widetilde{\mu}\D \Omega=\Lambda^{-1}{\rm curl}\, G,\\[2mm]
&u=-\Lambda^{-1}\n d-\Lambda^{-1}{\rm div}\,\Omega.\\[2mm]
\end{aligned}
\end{cases}
\leqno{(M^{'}_{1})}
$$
The last equation is just a heat equation. Hence we are going to
focus  on the first three equations. However the last equation
gives us an idea of which spaces we can work with.\\
The first three equation can be read as follows:
$$
\p_{t}\left(\begin{array}{c}
\hat{q}(t,\xi)\\
\hat{d}(t,\xi)\\
\hat{{\cal{T}}}(t,\xi)\\
\end{array}
\right)+A(\xi)\left(\begin{array}{c}
\hat{q}(t,\xi)\\
\hat{d}(t,\xi)\\
\hat{{\cal{T}}}(t,\xi)\\
\end{array}
\right)=\left(\begin{array}{c}
\hat{F}(t,\xi)\\
\Lambda^{-1}{\rm div}\,\hat{G}(t,\xi)\\
\hat{H}(t,\xi)\\
\end{array}
\right) \leqno{(M^{'}_{2})}
$$
where we have:
$$A(\xi)=\left(\begin{array}{ccc}
0&|\xi|& 0\\
-\e|\xi|^{3}-\beta|\xi|& \nu|\xi|^{2}& -\gamma|\xi|\\
0 &\delta|\xi| &\alpha|\xi|^{2}\\
\end{array}
\right)\,.
$$
The eigenvalues of the matrix $-A(\xi)$ are of the form
$|\xi|^{2}\lambda_{\xi}$ with $\lambda_{\xi}$ being the roots of the
following polynomial:
$$P_{\xi}(X)=X^{3}+(\nu+\alpha)X^{2}+\biggl(\e+\nu\alpha+\frac{\gamma\delta+\beta}{|\xi|^{2}}\biggl)X+
\biggl(\alpha\e+\frac{\alpha\beta}{|\xi|^{2}}\biggl).$$ For very
large $\xi$, the roots tend to those of the following polynomial (by
virtue of continuity of the roots in function of the coefficients):
$$X^{3}+(\nu+\alpha)X^{2}+(\e+\nu\alpha)X+\alpha\e.$$
The roots are  $-\alpha$ and
$-\frac{\nu}{2}(1\pm\sqrt{1-\frac{4\e}{\nu^{2}}})$. \\
The system $(M^{'}_{1})$ is well-posed if  and only if for $|\xi|$
tending to $+\infty$ the
real part of the eigenvalues associated to $A(\xi)$ stay non positive. Hence, we must have:\\
$$\e,\;\nu,\;\alpha\geq0.$$
\\
Let us now state a necessary and sufficient condition for the global
stability of $(M^{'})$.
\begin{proposition}
\label{fanalyse} The linear system  $(M^{'})$ is globally stable if
and only if the following conditions are verified:
$$\nu,\,\e,\,\alpha\geq
0,\;\;\alpha\beta\geq0,\;\;\gamma\delta(\nu+\alpha)+\nu\beta\geq0,\;\;\gamma\delta+\beta\geq0.\leqno{(*)}$$
If all the inequalities are strict, the solutions tend to 0 in the
sense of distributions and the three eigenvalues
$\lambda_{1}(\xi),\,\lambda_{+}(\xi),\,\lambda_{-}(\xi)$ have the
following asymptotic behavior  when $\xi$ tends to 0:\\
$$\lambda_{1}(\xi)\sim-\big(\frac{\alpha\beta}{\beta+\gamma\delta}\big)|\xi|^{2},\;
\lambda_{\pm}(\xi)\sim-\big(\frac{\gamma\delta(\nu+\alpha)+\nu\beta}{2(\gamma\delta+\beta)}\big)|\xi|^{2}\pm
i|\xi|\sqrt{\gamma\delta+\beta}.$$ \label{fconditions}
\end{proposition}
{\bf Proof} :\\
\\
We already know  that the system is well-posed if and only if
$\nu,\alpha\geq0$. We want that all the eigenvalues have a
negative real part for all  $\xi$.\\
We have to distinguish two cases: either all the eigenvalues are
real or there are two
 complex conjugated eigenvalues.
\subsubsection*{First case:}
The eigenvalues are real. A necessary condition for negativity of
the eigenvalues is that $P(X)\geq0$ for $X\geq0$. We must have in
particular:
$$P_{\xi}(0)=\alpha\e+\frac{\alpha\beta}{|\xi|^{2}}\geq0\;\;\;\forall\xi\ne0.$$
This imply that $\alpha\beta\geq0$ and $\alpha\e\geq 0$. Hence, given that $\alpha\geq0$, we must have $\beta\geq0$ and $\e\geq0$.
For $\xi$ tending to 0, we have:
$$P_{\xi}(\lambda)\sim\frac{\lambda(\gamma\delta+\beta)+\alpha\beta}{|\xi|^{2}}.$$
Making  $\lambda$ tend to infinity, we must have
$P_{\xi}(\lambda)\geq0$ and so $\gamma\delta+\beta\geq0.$ \\
The converse is trivial.
\subsubsection*{Second case:}
$P_{\xi}$ has two complex roots $z_{\pm}=a\pm ib$ and one real root
$\lambda$, we have:
$$P_{\xi}(X)=(X-\lambda)(X^{2}-2aX+|z_{\pm}|^{2}).$$
A necessary condition to have the real parts negative is in the same
way that $P_{\xi}(X)\geq0$ for all $X\geq0$.\\
If $\gamma\delta+\beta>0$, we are in the case where $\xi$ tends to
$0$ (and we see that $P_{\xi}$ is increasing).\\
We can observe the terms of degree 2 and we get:
$\lambda+2a=-\alpha-\nu$ then $\lambda$ and $\alpha$ are non
positive if and only if $P_{\xi}(-\alpha-\nu)\leq0$ (for this it
suffices to rewrite $P_{\xi}$ like
$P_{\xi}(X)=(X-\lambda)(X^{2}-2aX+|z_{\pm}|^{2})$).
Calculate:
$$P_{\xi}(-\alpha-\nu)=-\nu\e-\nu^{2}\alpha-\frac{\nu\beta+\nu\gamma\delta+\alpha\gamma\delta}{|\xi|^{2}}.$$
With the hypothesis that we have made, we deduce that
$P_{\xi}(-\alpha-\nu)\leq0$ for $\xi$ tending to $0$ if and only if
$\nu\beta+\nu\gamma\delta+\alpha\gamma\delta\geq0$.
\subsubsection*{Behavior of the eigenvalues in low frequencies: }
\label{fbasse} Let us now study the asymptotic behavior of the
eigenvalues when $\xi$
tends to $0$ and  all the inequalities in $(A)$ are strict.\\
We remark straight away that the condition $\gamma\delta+\beta>0$
ensures the strict monotonicity of the function: $\lambda\rightarrow
P_{\xi}(\lambda)$ for $\xi$ small. Then there's only one real
eigenvalue $\lambda_{1}(\xi)$ and two complex eigenvalues
$\lambda_{\pm}(\xi)=a(\xi)\pm ib(\xi)$.\\
\\
Let $\e^{-}<-\frac{\alpha\beta}{\gamma\delta+\beta}<\e^{+}<0$. When
$\xi$ tends to $0$, we have:
$$P_{\xi}(\lambda)\sim|\xi|^{-2}(\lambda(\gamma\delta+\beta)+\alpha\beta)).$$
Then $P_{\xi}(\e^{-})<0$ and $P_{\xi}(\e^{+})>0$ and $P_{\xi}$ has a
unique real root included between $\e^{-}$ and $\e^{+}$. These
considerations give the asymptotic value of $\lambda_{1}(\xi)$.\\
Finally, we have:
$$\lambda_{1}(\xi)+2a(\xi)=-\alpha-\nu\;\mbox{and}\;-(a(\xi)^{2}+b(\xi)^{2})\lambda(\xi)=\alpha\xi+\frac{\alpha\beta}{|\xi|^{2}}\sim
\frac{\alpha\beta}{|\xi|^{2}},$$
whence the result.
\hfill{$\Box$}
\\
\\
We summarize this results in the following remark.
\begin{remarka}
According to the analysis made in proposition \ref{fanalyse}, we expect the
system $(M)$ to be locally well-posed close to the equilibrium
$(\bar{\rho},0,\bar{T})$ if and only if we have:
$$\mu(\bar{\rho},\bar{\theta})\geq0,\;\lambda(\bar{\rho},\bar{\theta})+2\mu(\bar{\rho},\bar{\theta})\geq0,\;
\kappa(\bar{\rho})\geq0,\;\;\;\mbox{and}\;\;\chi(\bar{\rho},\bar{T})\geq0.\leqno{(C)}$$
By the calculus we have:
$$\beta=\p_{\rho}p_{0}(\bar{\rho},\bar{T}),\;\gamma=\frac{\p_{T}p_{0}(\bar{\rho},\bar{T})}{\bar{\rho}\,\p_{T}
e_{0}(\bar{\rho},\bar{T})},\;
\delta=\frac{\bar{T}\,\p_{T}p_{0}(\bar{\rho},\bar{T})}{\bar{\rho}}.$$
We remark that $\gamma\delta\geq0$ if
$\p_{T}e_{0}(\bar{\rho},\bar{T})\geq0$. In the case where $\eta$
verifies $\eta(\bar{\rho},\bar{T})>0$, the supplementary condition
giving the global stability reduces to:
$$\p_{\rho}p_{0}(\bar{\rho},\bar{T})\geq0.\leqno{(D)}$$
\end{remarka}
Now that we know the stability conditions on the coefficients of the
system $(M^{'})$, we aim at proving estimates in
the space $E^{\N}$.
\\
We add a condition in this following proposition compared with the
proposition \ref{fconditions} which is: $\gamma\delta>0$, but it's
not so important because in the system $(NHV)$ we are interested in,
we have effectively
$\gamma\delta=\frac{1}{\bar{T}\Psi^{'}(\bar{T})}>0$.
\begin{proposition}:
\label{flinear1} Under the conditions of proposition \ref{fconditions}
with strict inequalities and with the condition $\gamma\delta>0$,
let $(q,d,{\cal{T}})$ be a solution of the system $(M^{'})$ on
$[0,T)$ with initial data $ (q_{0},u_{0},{\cal{T}}_{0}) $ such that:
$$ q_{0}\in \widetilde{B}^{s-1,s},d_{0}\in B^{s-1},{\cal{T}}_{0}\in\widetilde{B}^{s-1,s-2}\;\;\mbox{for some $s\in\R$}\;.$$
Moreover  we suppose that for some $1\leq r_{1}\leq+\infty$, we have:\\
$$F\in\widetilde{L}_{T}^{r_{1}}(\widetilde{B}^{s-3+\frac{2}{r_{1}},s-2+\frac{2}{r_{1}}}),\;\;
G\in\widetilde{L}_{T}^{r_{1}}(B^{s-3+\frac{2}{r_{1}}}),\;\;
H\in\widetilde{L}_{T}^{r_{1}}(\widetilde{B}^{s-3+\frac{2}{r_{1}},s-4+\frac{2}{r_{1}}}).$$
We then have  the following estimate for all $r_{1}\leq r\leq
+\infty$:
$$
\begin{aligned}
&\|q\|_{\widetilde{L}_{T}^{r}(\widetilde{B}^{s-1+\frac{2}{r},s+\frac{2}{r}})}+\|
{\cal{T}}\|_{\widetilde{L}_{T}^{r}(\widetilde{B}^{s-1+\frac{2}{r},s-2+\frac{2}{r}})}+\|
u\|_{\widetilde{L}_{T}^{r}(B^{s-1+\frac{2}{r}})}\lesssim\| q_{0}\|_{\widetilde{B}^{s-1,s}}+\| u_{0}\|_{B^{s-1}}\\[2mm]
&+\| {\cal T}_{0}\|_{\widetilde{B}^{s-1,s-2}} +\|
F\|_{L_{T}^{r_{1}}(\widetilde{B}^{s-3+\frac{2}{r_{1}},s-2+\frac{2}{r_{1}}})}+\|
G\|_{ L_{T}^{r_{1}}(B^{s-3+\frac{2}{r_{1}}})}+\|
H\|_{L_{T}^{r_{1}}(\widetilde{B}^{s-3+\frac{2}{r_{1}},s-4+\frac{2}{r_{1}}})}\;.\\
\end{aligned}
$$
\end{proposition}
{\bf Proof}:\\
\\
We are going to separate the case of the low, medium and high
frequencies, particulary the low and high frequencies which have a
different behavior, and  depend on the indice of Besov space.
\subsubsection*{1) Case of low frequencies:}
Let us focus on just the first three equation because the
last one is a heat equation that we can treat independently.
Applying operator $\D_{l}$ to the system $(M^{'}_{1})$, we obtain then in setting:
$$q_{l}=\Delta_{l}q,\;d_{l}=\Delta_{l}d,\;{\cal{T}}_{l}=\Delta_{l}{\cal{T}}$$
the following system:
\begin{eqnarray}
&&\p_{t}q_{l}+\Lambda d_{l}=F_{l},\label{feq1}\\
&&\p_{t}d_{l}-\nu\Delta d_{l}-\e \Lambda^{3}q_{l}-\beta \Lambda
q_{l}-\gamma\Lambda{\cal{T}}_{l}=\Lambda^{-1}{\rm div} G_{l},\label{feq2}\\
&&\p_{t}{\cal{T}}_{l}-\alpha\Delta{\cal{T}}_{l}+\delta\Lambda
d_{l}=H_{l}.\label{feq3}
\end{eqnarray}
Throughout the proof, we assume that $\delta\ne 0$: if not we have
just a heat equation on (\ref{feq3}) and we can use the proposition
\ref{fchaleur} to have the estimate on ${\cal T}$ and we have just to
deal with the  first two equations.
Denoting by $W(t)$ the semi-group associated to $(\ref{feq1}-\ref{feq3})$ we have:\\
$$\left(\begin{array}{c}
q(t)\\
u(t)\\
\theta(t)\\
\end{array}
\right) =W(t)\left(\begin{array}{c}
q_{0}\\
u_{0}\\
\theta_{0}\\
\end{array}
\right)+\int_{0}^{t}W(t-s)\left(\begin{array}{c}
F(s)\\
G(s)\\
H(s)\\
\end{array}
\right)\ ds\;.$$ We set:
$$f_{l}^{2}=\beta\|q_{l}\|_{L^{2}}^{2}+\|d_{l} \|_{L^{2}}^{2}+\frac{\gamma}{\delta}\|{\cal{T}}_{l}\|_{L^{2}}^{2}-2K\langle\Lambda q_{l},d_{l}\rangle$$
for some $K\geq0$ to be fixed hereafter and
$\langle\cdot,\cdot\rangle$ noting the $L^{2}$ inner product.
\\
To begin with, we consider the case  where $F=G=H=0$.\\
 Then we
take the inner product of (\ref{feq2}) with $d_{l}$, of (\ref{feq1})
with $\beta q_{l}$ and of (\ref{feq3}) with $\gamma{\cal T}_{l}$. We
get:
\begin{equation}
\frac{1}{2}\frac{d}{dt}\big(\|d_{l}\|_{L^{2}}^{2}+\beta\| q_{l}
\|_{L^{2}}^{2}+\frac{\gamma}{\delta}\|{\cal{T}}_{l}\|_{L^{2}}^{2}\big)+\nu\|\nabla
d_{l}\|_{L^{2}}^{2}
-\e\langle\Lambda^{3}q_{l},d_{l}\rangle+\frac{\gamma\alpha}{\delta}\|
\nabla{\cal{T}}_{l}\|_{L^{2}}^{2}=0. \label{feq4}
\end{equation}
Next, we apply the operator $\Lambda$ to (\ref{feq2}) and take the
inner product with $q_{l}$, and we take the scalar product of
(\ref{feq1}) with $\La d_{l}$ to control the term
$\frac{d}{dt}\langle\Lambda q_{l},d_{l}\rangle$. Summing the two
resulting equalities, we get:
\begin{equation}
\frac{d}{dt}\langle\Lambda q_{l},d_{l}\rangle+\| \Lambda
d_{l}\|_{L^{2}}^{2}-\nu\langle\D d_{l},\Lambda
q_{l}\rangle-\e\|\Lambda^{2} q_{l}\|_{L^{2}}^{2} -\beta\|\Lambda
q_{l}\|_{L^{2}}^{2}-\gamma\langle\Lambda {\cal{T}}_{l},\Lambda
q_{l}\rangle=0 .\label{feq5}
\end{equation}
We obtain then in summing (\ref{feq4}) and (\ref{feq5}):\\
\begin{equation}
\begin{aligned}
&\frac{1}{2}\frac{d}{dt}f_{l} ^{2}+(\nu\|\n
d_{l}\|_{L^{2}}^{2}-K\|\Lambda d_{l}\|_{L^{2}}^{2})+(K\beta\|\Lambda
q_{l}\|_{L^{2}}^{2}+K\e\|\Lambda^{2}q_{l}\|_{L^{2}}^{2})
+\frac{\gamma\alpha}{\delta}\|\n {\cal{T}}_{l}\|_{L^{2}}^{2}\\[3mm]
&+K \nu\langle\D d_{l},\Lambda q_{l}\rangle+K\gamma\langle\Lambda
{\cal{T}}_{l},\Lambda q_{l}\rangle -\e\langle \Lambda^{3}
q_{l},d_{l}\rangle=0\,.
\label{feq6}\\
\end{aligned}
\end{equation}
Like indicated, we are going to focus on low frequencies so
assume that $l\leq l_{0}$ for some $l_{0}$ to be fixed hereafter. We have then $\forall c,\,b,\,d>0$ :\\
\begin{equation}
\begin{aligned}
|\langle\D d_{l},\Lambda q_{l}\rangle|&\leq\frac{b}{2}\|\Lambda
q_{l}\|_{L^{2}}^{2}+\frac{1}{2b}\|\ \D d_{l}\|_{L^{2}}^{2}\\
&\leq\frac{b}{2}\|\Lambda q_{l}\|_{L^{2}}^{2}+\frac{C2^{2l_{0}}}{2b}
\|\Lambda d_{l}\|_{L^{2}}^{2}\,,\\
\end{aligned}
\label{feq7}
\end{equation}
$$ |\langle\Lambda^{3} q_{l},d_{l}\rangle|=|\langle\Lambda^{2} q_{l},\Lambda d_{l}\rangle|\leq\frac{C2^{2l_{0}}}{2c}\|\Lambda q_{l}\|_{L^{2}}^{2}+\frac{c}{2}\|\Lambda d_{l}\|_{L^{2}}^{2}\,.$$
\\
Moreover we have: $\|\n d_{l}\|_{L^{2}}^{2}=\|\Lambda
d_{l}\|_{L^{2}}^{2}.$ Finally we obtain:
$$
\begin{aligned}
\frac{1}{2}\frac{d}{dt}f_{l}
^{2}+\big[\nu-(K+\frac{C2^{2l_{0}}}{2b}&K\nu+\frac{c\e}{2})\big]\|\Lambda
d_{l}\|_{L^{2}}^{2}+\big[\frac{\gamma\alpha}{\delta}-\frac{K\gamma}{2d}\big]\|\Lambda
{\cal{T}}_{l}\parallel_{L^{2}}^{2}\\[2mm]
&+K\big[\beta+\e C
2^{2l_{0}}-\nu\frac{b}{2}-\e\frac{C2^{2l_{0}}}{2c}-\gamma\frac{d}{2}\big]\|\Lambda q_{l}\|_{L^{2}}^{2}\leq0.\\[2mm]
\end{aligned}
$$
Then we choose $(b,c,d)$ such that:
$$b=\frac{\beta}{2\nu}\;,\,c=\frac{\nu}{\e}\;,\,d=\frac{\beta}{2\gamma}\;,$$
which is possible if $\gamma>0$ as $\nu>0,\;\e>0\;$. In the case
where $\gamma\leq 0$, we recall that $\gamma$ and $\delta$ have
the same sign, we have then no problem because with our choice the
first and third following inequalities will be satisfied and if
$\gamma\leq 0$ in the second equation the term $\gamma\frac{d}{2}$
is positive in taking $d>0$. So we assume from now on that
$\gamma>0$ and so with this choice, we want that:
$$
\begin{aligned}
&\frac{\nu}{2}-K(1+C2^{2l_{0}}\frac{\nu^{2}}{\beta})>0,\\[2mm]
&\frac{\beta}{2}+\e
C2^{2l_{0}}-C2^{2l_{0}}\frac{\e^{2}}{2\nu}>0,\\[2mm]
&\frac{\gamma\alpha}{\delta}-K\frac{\gamma}{2}>0.\\
\end{aligned}
$$
We recall that in your case $\nu>0,\;\beta>0,\;\alpha>0$ and
$\gamma>0$, $\delta>0$.
So it suffices to choose $K$ and $l_{0}$ such that:
$$K<\min\biggl(\frac{\nu}{2(1+C2^{2l_{0}}\frac{\nu^{2}}{2\beta})},\frac{2\alpha}{\delta}\biggl)\;\;\;\mbox{and}\;\;\;2^{2l_{0}}<\min\biggl(\frac{\beta\nu}{6C\e^{2}},\frac{1}{6\e C}\biggl)\,.$$
Finally we conclude in using Proposition \ref{fderivation} part (ii) with a $c^{'}$ small enough. We get:
\begin{equation}
\frac{1}{2}\frac{d}{dt}f_{l}
^{2}+c^{'}2^{2l}f_{l}^{2}\leq0\;\;\mbox{for $l\leq l_{0}$}\;
.\label{feq8}
\end{equation}
\subsubsection*{2) Case of high frequencies:}
We are going to work with $l\geq l_{1}$ where we will determine
$l_{1}$ hereafter. We set then:\\
$$f_{l}^{2}=\e B\|\Lambda q_{l}\|_{L^{2}}^{2}+B\|d_{l}\|_{L^{2}}^{2}+\|\Lambda^{-1} {\cal{T}}_{l}\|_{L^{2}}^{2}
-2K\langle\Lambda q_{l},d_{l}\rangle,$$
and we choose $B$ and $K$ later on.\\
Then we take the  inner product of $(\ref{feq2})$ with $d_{l}$:
\begin{equation}
\frac{1}{2}\frac{d}{dt}\|d_{l}\|_{L^{2}}^{2}+\nu \|\n
d_{l}\|_{L^{2}}^{2}-\e  \langle\Lambda^{3}
q_{l},d_{l}\rangle-\beta\langle\Lambda
q_{l},d_{l}\rangle-\gamma\langle\Lambda{\cal{T}}_{l},d_{l}\rangle=0.
\label{feq9}
\end{equation}
Moreover we have in taking the scalar product of $(\ref{feq1})$ with $\Lambda^{2} q_{l}$:
\begin{equation}
\frac{1}{2}\frac{d}{dt}\|\Lambda
q_{l}\|_{L^{2}}^{2}+\langle\Lambda^{2} d_{l},\Lambda q_{l}\rangle=0.
\label{feq10}
\end{equation}
And in the same way with $(\ref{feq3})$, we have:
\begin{equation}
\frac{1}{2}\frac{d}{dt}\|\Lambda^{-1}{\cal{T}}_{l}\|_{L^{2}}^{2}+\alpha\|{\cal
T}_{l}\|_{L^{2}}^{2}+\delta\langle
 d_{l},\Lambda^{-1}{\cal{T}}_{l}\rangle=0. \label{feq11}
\end{equation}
After we sum (\ref{feq9}), (\ref{feq10}) and (\ref{feq11}) to get:
\begin{equation}
\begin{aligned}
&\frac{1}{2}\frac{d}{dt}\big(B\|d_{l}\|_{L^{2}}^{2}+\e B\|\Lambda
q_{l}\|_{L^{2}}^{2}+\|\Lambda^{-1}{\cal{T}}_{l}\|_{L^{2}}^{2}\big)
+B\nu\|\nabla d_{l}\|_{L^{2}}^{2}+\alpha\|
{\cal{T}}_{l}\|_{L^{2}}^{2}\\[2mm]
&\hspace{4cm}-B\beta\langle\Lambda
q_{l},d_{l}\rangle
-B\gamma\langle\Lambda
{\cal{T}}_{l},d_{l}\rangle+\delta\langle
d_{l},\Lambda^{-1}{\cal{T}}_{l}\rangle=0.\\
\end{aligned}
\label{feq12}
\end{equation}
Then like previously we can play with  $\langle\Lambda
q_{l},d_{l}\rangle$ to obtain a term in $\|\Lambda
q_{l}\|_{L^{2}}^{2}$.
We have then again the following equation:\\
\begin{equation}
\frac{d}{dt}\langle\Lambda q_{l},d_{l}\rangle+\|\Lambda
d_{l}\|_{L^{2}}^{2}-\nu\langle\Delta d_{l},\Lambda
q_{l}\rangle-\e\|\Lambda^{2} q_{l}\|_{L^{2}}^{2}-\beta\|\Lambda
q_{l}\|_{L^{2}}^{2}-\gamma\langle\Lambda{\cal{T}}_{l},\Lambda
q_{l}\rangle=0. \label{feq13}
\end{equation}
We sum all these expressions and get:
\begin{equation}
\begin{aligned}
&\frac{1}{2}\frac{d}{dt}f_{l}^{2}+\big[B\nu\|\nabla
d_{l}\|_{L^{2}}^{2}-K\|\Lambda
d_{l}\|_{L^{2}}^{2}\big]+\alpha\|{\cal{T}}_{l}\|_{L^{2}}^{2}+K\big[\beta
\|\Lambda q_{l}\|_{L^{2}}^{2}+\e \|\Lambda^{2}
q_{l}\|_{L^{2}}^{2}\big]\\[2mm]
&-B\beta\langle\Lambda q_{l},d_{l}\rangle-B\gamma\langle\Lambda
{\cal{T}}_{l},d_{l}\rangle+\delta\langle
d_{l},\Lambda^{-1}{\cal{T}}_{l}\rangle+K\nu\langle\Delta
d_{l},\Lambda q_{l}\rangle+\gamma K\langle\Lambda
{\cal{T}}_{l},\Lambda
q_{l}\rangle=0.\\
\end{aligned}
\label{feq141}
\end{equation}
The main term in high frequencies will be: $\|\Lambda^{2}
q_{l}\|_{L^{2}}^{2}$. The other terms may be treated by mean of
Young's inequality:
$$
\begin{aligned}
|\langle\Lambda q_{l},d_{l}\rangle|&\leq\frac{1}{2a}\|\Lambda
q_{l}\|_{L^{2}}^{2}+\frac{a}{2}
\|\Lambda d_{l}\|_{L^{2}}^{2},\\
&\leq\frac{1}{2a\,c2^{2l_{1}}}\|\Lambda^{2}q_{l}\|_{L^{2}}^{2}+\frac{a}{2}\|\Lambda d_{l}\|_{L^{2}}^{2}.\\
\end{aligned}
$$
\\
We do as before with the others terms in the second line of
(\ref{feq141}) and we obtain:
$$
\begin{aligned}
&\frac{1}{2}\frac{d}{dt}f_{l}^{2}+(B\nu-K)\|\Lambda
d_{l}\|_{L^{2}}^{2}\big]+\alpha\|{\cal{T}}_{l}\|_{L^{2}}^{2}
+K(\frac{\beta }{c2^{2l_{1}}}+\e)
\|\Lambda^{2} q_{l}\|_{L^{2}}^{2}\leq\\[2,5mm]
&B\gamma\big[\frac{1}{2a}\|
{\cal{T}}_{l}\|_{L^{2}}^{2}+\frac{a}{2}\|\Lambda
d_{l}\|_{L^{2}}^{2}\big]+K\big[\frac{\nu
b}{2}\|\Lambda^{2}q_{l}\|_{L^{2}}^{2}+\frac{2\nu}{b}\|\Lambda
d_{l}\|_{L^{2}}^{2}+
\frac{\gamma}{2c^{'}}\|{\cal{T}}_{l}\|_{L^{2}}^{2}+\frac{\gamma
c^{'}}{2}\|
\Lambda^{2} q_{l}\|_{L^{2}}^{2}\big]\\[2,5mm]
&+B\beta\big[\frac{1}{2d}\frac{1}{c2^{2l_{1}}}\|\Lambda^{2}q_{l}\|_{L^{2}}^{2}+\frac{d}{2}\frac{1}{c2^{2l_{1}}}\|
\Lambda
d_{l}\|_{L^{2}}^{2}\big]+\delta\big[\frac{1}{2e}\frac{1}{c2^{2l_{1}}}\|{\cal{T}}_{l}\|_{L^{2}}^{2}
+\frac{e}{2}\frac{1}{c2^{2l_{1}}}\|\Lambda d_{l}\|_{L^{2}}^{2}\big].\\
\end{aligned}
$$
\\
We obtain then for some $a,\;b,\;c^{'},\;d,\;e$ to be chosen:
\begin{equation}
\begin{aligned}
&\frac{1}{2}\frac{d}{dt}f_{l}^{2}+\big[B\nu-(K+B\gamma\frac{a}{2}+K\nu\frac{2}{b}+B\beta\frac{d}{2c2^{2l_{1}}}
+\delta\frac{e}{2}\frac{1}{c2^{2l_{1}}})\big]\|\Lambda
d_{l}\|_{L^{2}}^{2}\\
&+\big[\alpha-(B\gamma \frac{1}{2a}+\gamma
K\frac{1}{2c^{'}}+\delta\frac{1}{2e}\frac{1}{c2^{2l_{1}}})\big]\|{\cal{T}}_{l}\|_{L^{2}}^{2}\\
&+\big[\frac{\beta K}{c2^{2l_{1}}}+\e K-K\nu\frac{b}{2}-\gamma
K\frac{c^{'}}{2}-B\beta\frac{1}{2d}\frac{1}{2^{2l_{1}}}\big]\|\Lambda^{2}
q_{l}\|_{L^{2}}^{2}\leq0\;.\\
\end{aligned}
\label{feq15}
\end{equation}
We claim that $a$, $b$, $c^{'}$, $d$, $e$, $l_{1}$, $K$ may be
chosen so that:
\begin{eqnarray}
&&B\nu-(K+B\gamma\frac{a}{2}+K\nu\frac{2}{b}+B\beta\frac{d}{2}+\delta\frac{e}{2}\frac{1}{2^{2l_{1}}})>0,\label{f51}\\
&&\alpha-(B\gamma \frac{1}{2a}+\gamma K\frac{1}{2c^{'}}+\delta\frac{1}{2e}\frac{1}{c2^{2l_{1}}})>0,\label{f52}\\
&&\frac{\beta K}{c2^{2l_{1}}}+\e K-K\nu\frac{b}{2}-\gamma
K\frac{c^{'}}{2}-B\beta\frac{1}{2d}\frac{1}{2^{2l_{1}}}>0.\label{f53}_{}
\end{eqnarray}
We want at once that for (\ref{f51}) and (\ref{f53}):
\begin{eqnarray}
&&\nu-\gamma \frac{a}{2}-\beta\frac{d}{2}>0\;,\label{f54}\\
&&\e-\nu\frac{b}{2}-\gamma\frac{c^{'}}{2}>0\;.\label{f55}
\end{eqnarray}
So we take:
$$
\begin{aligned}
&e=1,\;a=2h\delta,\;d=2h,\;h=\frac{\nu}{2(\gamma\delta+\beta)},\;b=2\beta
h^{'},\;\\
&c^{'}=2\delta(\nu+\alpha)h^{'}\;\;\;\mbox{and}\;\;\;h^{'}=\frac{\e}{2(\gamma\delta(\nu+\alpha)+\nu\beta)}\;.\\
\end{aligned}
$$
\\
With this choice, we get (\ref{f54}) and (\ref{f55}). In what follows
it suffices to choice $B,K$ small enough and $l_{1}$ large enough.
We have then:
$$f_{l}\simeq \mbox{Max}(1,2^{l})\| q_{l}\|_{L^{2}}+\| d_{l}\|_{L^{2}}+\mbox{Min}(1,2^{l})\| {\cal{T}}_{l}\|_{L^{2}}$$
We have so obtain for $l\leq l_{0}$, $l\geq l_{1}$ and for a $c^{'}$
small enough:
$$\frac{1}{2}\frac{d}{dt}f_{l} ^{2}+c^{'}2^{2l}f_{l}^{2}\leq0.$$
\subsubsection*{3) Case of Medium frequencies:}
For $l_{0}\leq l\leq l_{1}$, there is only a
finite number of terms to treat. So it suffices to find
a $C$ such that for all these terms:
$$
\begin{aligned}
&\| q_{l}\|_{L_{T}^{r}(L^{2})}\leq C,\;\| d_{l}\|_{L_{T}^{r}(L^{2})}\leq C,\;\| {\cal{T}}_{l}\|_{L_{T}^{r}(L^{2})}\leq C\;\;\;\mbox{for all}\;\;\;T\in[0,+\infty]\\
&\hspace{8,2cm}\;\;\mbox{and}\;\;r\in[1,+\infty]\\
\end{aligned}
\leqno{(B)}
$$
with $C$ large enough independent of $T$.\\
And this is true because the system is globally stable: indeed according to proposition \ref{fanalyse}, we have:
$$\biggl\|W(t) \left(\begin{array}{c}
a\\
b\\
c\\
\end{array}
\right)\biggl\|_{L^{2}}\lesssim e^{-c_{1}(\xi)t}\biggl\|
\left(\begin{array}{c}
a\\
b\\
c\\
\end{array}
\right)\biggl\|_{L^{2}}\;\;\forall a,b,c\in L^{2}
$$
with $c_{1}(\xi)=\min_{2^{l_{0}}\leq|\xi|\leq
2^{l_{1}}}(\mbox{Re}(\lambda_{1}(\xi)),\;\mbox{Re}(\lambda_{2}(\xi)),\;\mbox{Re}(\lambda_{3}(\xi)))$
where the $\lambda_{i}(\xi)$ correspond to the eigenvalues of the
system.
We have then in using the estimate in low and high frequencies in part \ref{fbasse} and the continuity of $c_{1}(\xi)$ the fact that there exits $c_{1}$ such that:\\
$$c_{1}(\xi)\geq c_{1}>0.$$
So that we have:\\
$$
\biggl(\int^{T}_{0}\left(\begin{array}{c}
\|q_{l}(t)\|_{L^{2}}^{r}\\
\|u_{l}(t)\|_{L^{2}}^{r}\\
\|{\cal{T}}_{l}(t)\|_{L^{2}}^{r}\\
\end{array}
\right)dt \biggl)^{\frac{1}{r}}\lesssim
\biggl(\int^{T}_{0}e^{-c_{1}rs}ds\biggl)^{\frac{1}{r}}\left(\begin{array}{c}
\|(q_{0})_{l}\|_{L^{2}}\\
\|(u_{0})_{l}\|_{L^{2}}\\
\|({\cal{T}}_{0})_{l}\|_{L^{2}}\\
\end{array}
\right)\;\;\mbox{for}\;\;l_{0}\leq l\leq l_{1}.
$$
And so we have the result $(B)$.
\subsubsection*{4) Conclusion:}
In using Duhamel formula for $W$ and in taking $C$ large enough  we have for all $l$:\\
$$
\begin{aligned}
&\max (1,2^{l})\|q_{l}(t)\|_{L^{2}}+\|d_{l}(t)\|_{L^{2}}+\min
(1,2^{-l})\|{\cal{T}}_{l}(t)\|_{L^{2}}\leq Ce^{-c2^{2l}t}\big(\max
(1,2^{l})\|(q_{0})_{l}\|_{L^{2}}\\[2mm]
&+\|(d_{0})_{l}\|_{L^{2}}+\min
(1,2^{-l})\|({\cal{T}}_{0})_{l}\|_{L^{2}}\big)
+C\int^{t}_{0}e^{-c2^{2l}(t-s)}\big(\max(1,2^{2l})\|F_{l}\|_{L^{2}}+\|G_{l}\|_{L^{2}}\\[2mm]
&\hspace{10cm}+\min(1,2^{-l})\|H_{l}\|_{L^{2}})ds\,.\\
\end{aligned}
$$
Now we take the $L^{r}$ norm in time and we sum in multiplying by
$2^{l(s-1+\frac{2}{r})}$ for the low frequencies and we sum in
multiplying by $2^{l(s+\frac{2}{r})}$ for the high frequencies.\\
This yields:
$$
\begin{aligned}
&\|q\|_{\widetilde{L}^{r}_{T}(\widetilde{B}^{s-1+\frac{2}{r},s+\frac{2}{r}})}+\|{\cal{T}}\|_{\widetilde{L}^{r}_{T}(\widetilde{B}^{s-1+\frac{2}{r},s-2+\frac{2}{r}})}
+\|d\|_{\widetilde{L}^{r}_{T}(B^{s-1+\frac{2}{r}})}
\leq\|q_{0}\|_{\widetilde{B}^{s-1,s}}
+\|{\cal{T}}_{0}\|_{\widetilde{B}^{s-1,s-2}}\\[1,7mm]
&+\|d_{0}\|_{B^{s-1}}+\sum_{l\leq
0}2^{l(s-1+\frac{2}{r})}\int^{T}_{0}\biggl(\int^{T}_{0}e^{c(t-\tau)}\big(
\|F_{l}(\tau)\|_{L^{2}}+\|G_{l}(\tau)\|_{L^{2}}+\|H_{l}(\tau)\|_{L^{2}}\big)d\tau\biggl)^{r}dt\biggl)^{\frac{1}{r}}\\[1,7mm]
&+\sum_{l\geq
0}2^{l(s+\frac{2}{r})}\biggl(\int^{T}_{0}\biggl(\int^{T}_{0}e^{c(t-\tau)}\big(
\|\n F_{l}(\tau)\|_{L^{2}}+\|G_{l}(\tau)\|_{L^{2}}+\|\Lambda^{-1}
H_{l}(\tau)\|_
{L^{2}}\big)d\tau\biggl)^{r}dt\biggl)^{\frac{1}{r}}\,.\\[1,7mm]
\end{aligned}
$$
Bounding the right hand-side may be done by taking advantage of convolution inequalities.
To complete the proof of proposition \ref{flinear1}, it suffices to use that $u=-\Lambda^{-1}\n d-\Lambda^{-1}\rm div\Omega$ and to apply proposition \ref{fchaleur}. \hfill
{$\Box$}
\subsection{Global existence for temperature independent coefficients\label{fsubsection42}}
This section is devoted to the proof of theorem \ref{ftheo1} and \ref{ftheo3}. Let
us first  recall the spaces  in which we work with for the theorem
\ref{ftheo1}:
$$
\begin{aligned}
E^{s}=&[C_{b}(\R_{+},\widetilde{B}^{s-1,s})\cap
L^{1}(\R_{+}\widetilde{B}^{s+1,s+2})]\times[
C_{b}(\R_{+},B^{s-1})^{N}\cap L ^{1}(\R_{+},B^{s+1})^{N}]\\
&\times[ C_{b}(\R_{+},\widetilde{B}^{s-1,s-2})\cap
L^{1}(\R_{+},\widetilde{B}^{s+1,s})]. \\
\end{aligned}
$$
In what follows, we assume that $N\geq3$.\\
\\
{\bf Proof of  theorem \ref{ftheo1}:}\\
\\
We shall use a contracting mapping argument for the function $\psi$ defined as follows:
\begin{equation}
\psi(q,u,{\cal T})=W(t,\cdot)*\left(\begin{array}{c}
q_{0}\\
u_{0}\\
{\cal T}_{0}\\
\end{array}
\right)+\int_{0}^{t}W(t-s)\left(\begin{array}{c}
F(q,u,{\cal T})\\
G(q,u,{\cal T})\\
H(q,u,{\cal T})\\
\end{array}
\right)\ ds\; .\label{f1a1}
\end{equation}
In what follows we set:
$$\rho=\bar{\rho}(1+q)\;,\;\theta=\bar{\theta}+{\cal T}\;,\;\widetilde{T}=\Psi^{-1}(\theta).$$
The non linear terms  $F,G,H$ are defined as follows:
\begin{equation}
\begin{aligned}
F=&-{\rm div}(qu),\\[3mm]
G=&-u.\n u+\nabla(\frac{K^{'}_{\rho}}{2}|\nabla\rho|^{2}\
)+\biggl[\frac{\mu(\rho)}{\rho}-\frac{\mu(\bar{\rho})}{\bar{\rho}}\biggl]\D
u+\biggl[\frac{\zeta(\rho)}{\rho}-\frac{\zeta(\bar{\rho})}{\bar{\rho}}\biggl]\nabla
{\rm
div}\,u\\[2mm]
&+(\n((K(\rho)-K(\bar{\rho}))\D\rho)
+\biggl[\frac{P_{0}^{'}(\rho)+\w
P_{1}^{'}(\rho)}{\rho}-\frac{P_{0}^{'}(\bar{\rho})
+\bar{T}P_{1}^{'}(\bar{\rho})}{\bar{\rho}}\biggl]\n
\rho\\[2mm]
&+\biggl[\frac{P_{1}(\rho)}{\rho
\Psi^{'}(\cal{T})}-\frac{P_{1}(\bar{\rho})}{\bar{\rho}\Psi^{'}(\bar{T})}\biggl]\nabla\theta+\frac{\lambda^{'}(\rho)\n\rho{\rm
div}u}{\rho}+\frac{(du+\n u)\mu^{'}(\rho)\n\rho}{\rho}\;,
\label{f1a3}\\
\end{aligned}
\end{equation}
\\
where we note: $\zeta=\lambda+\mu$, and:
\\
\begin{equation}
H=\biggl(\frac{{\rm
div}(\chi(\rho)\nabla\theta)}{\rho}-\frac{\bar{\chi}}{\bar{\rho}}\Delta\theta
\biggl)+\biggl[\frac{\bar{T}P_{1}(\bar{\rho})}{\bar{\rho}} -\frac{\w
P_{1}(\rho )}{\rho}\biggl]{\rm div} u-u^{*}.\n\theta+\frac{D:\nabla
u}{\rho}\;. \label{f1a4}
\end{equation}
$ $
\subsubsection*{1) First step, uniform bounds:}
Let:
$$\eta=\| q_{0}\|_{\widetilde{B}^{\N-1,\N}}+\|u_{0}\|_{B^{\N-1}}+\|{\cal T}_{0}\|_{\widetilde{B}^{\N-1,\N}}\;.$$
We are going to show that $\psi$  maps the ball $B(0,R)$ into itself if $R$ is small enough.
According to proposition \ref{flinear1}, we have:
\begin{equation}
\|W(t,\cdot)*\left(\begin{array}{c}
q_{0}\\
u_{0}\\
{\cal T}_{0}\\
\end{array}
\right)\|_{E^{\N}}\leq C\eta \;.\label{f1a5}
\end{equation}
We have then according (\ref{f1a1}), proposition \ref{flinear1} and \ref{f1a5}:
\begin{equation}
\begin{aligned}
\|\psi(q,u,{\cal T})\|_{E^{\N}}\leq&\;C\eta+
\|F(q,u,{\cal T})\|_{ L^{1}(\widetilde{B}^{\N-1,\N})}+\|G(q,u,{\cal T})\|_{ L^{1}(B^{\N-1})}\\[2mm]
&\hspace{4,5cm}+\| H(q,u,{\cal T})\|_{L^{1}(\widetilde{B}^{\N-1,\N-2})}\;.\\
\end{aligned}
\label{f1a6}
\end{equation}
Moreover we suppose for the moment that:
$$\|q\|_{L^{\infty}(\R\times \R^{N})}\leq 1/2\;.\leqno{(\cal{H})}$$
We will use the different theorems on the paradifferential
calculus to obtain estimates on
$$\|F(q,u,{\cal T})\|_{ L^{1}(\widetilde{B}^{\N-1,\N})},\;\|G(q,u,{\cal T})\|_{ L^{1}(B^{\N-1})}\;\;\mbox{and}
\;\;\| H(q,u,{\cal T})\|_{L^{1}(\widetilde{B}^{\N-1,\N-2})}.$$
1) Let us first estimate $\|F(q,u,{\cal T})\|_{
L^{1}(\widetilde{B}^{\N-1,\N})}$. According to proposition
\ref{fproduit},
we have:
$$
\begin{aligned}
&\|{\rm div}(qu)\|_{
L^{1}(\widetilde{B}^{\N-1,\N})}\leq\|qu\|_{L^{1}(B^{\N})}+\|qu\|_{L^{1}(B^{\N+1})}\hspace{5cm}\\
\mbox{and:}\hspace{1cm}&\\
&\|qu\|_{L^{1}(B^{\N})}\leq\|q\|_{L^{2}(B^{\N})}\|u\|_{L^{2}(B^{\N})}\hspace{5cm}\\[2,5mm]
&\|qu\|_{L^{1}(B^{\N+1})}\leq\|q\|_{L^{\infty}(B^{\N})}\|u\|_{L^{1}(B^{\N+1})}+\|q\|_{L^{2}(B^{\N+1})}
\|u\|_{L^{2}(B^{\N})}.\hspace{5cm}\\
\end{aligned}
$$
Because $\widetilde{B}^{\N,\N+1}\h B^{\N}$ and
$\widetilde{B}^{\N,\N+1}\h B^{\N+1}$ (from proposition
\ref{fimportant}), we get:
$$\|{\rm div}(qu)\|_{L^{1}(\widetilde{B}^{\N-1,\N})}\leq\|q\|_{L^{\infty}(\widetilde{B}^{\N,\N+1})}
\|u\|_{L^{1}(B^{\N+1})}+\|q\|_{L^{2}(\widetilde{B}^{\N,\N+1})}\|u\|_{L^{2}(B^{\N})}.$$
\\
2) We have  to estimate $\|G(q,u,{\cal T})\|_{L^{1}(B^{\N-1})}$. We
see straight away that:
$$[\frac{\mu(\rho)}{\rho}-\frac{\mu(\bar{\rho})}{\bar{\rho}}]\Delta
u=K(q)\Delta u$$ for some  smooth function $K$ such that
$K(0)=0$. Hence by propositions \ref{fcomposition}, \ref{fproduit}
and \ref{fimportant} yield:
$$
\begin{aligned}
\biggl\|[\frac{\mu(\rho)}{\rho}-\frac{\mu(\bar{\rho})}{\bar{\rho}}]\Delta u\biggl\|
_{L^{1}(B^{\N-1})}\lesssim&\|K(q)\|_{L^{\infty}(B^{\N})}\|u\|_{L^{1}(B^{\N+1})},\\[2mm]
\lesssim&
\|q\|_{L^{\infty}(B^{\N})}\|u\|_{L^{1}(B^{\N+1})},\\[2mm]
\lesssim&\|q\|_{L^{\infty}(\widetilde{B}^{\N-1,\N})}\|u\|_{L^{1}(B^{\N+1})}\;.\\[2mm]
\end{aligned}
$$
In the same way we have:
$$
\begin{aligned}
&\|[\frac{\zeta(\rho)}{\rho}-\frac{\zeta(\bar{\rho})}{\bar{\rho}}]
\nabla{\rm div}u\|_{L^{1}(B^{\N-1})}\lesssim
\|q\|_{L^{\infty}(\widetilde{B}^{\N-1,\N})}\|u\|_{L^{1}(B^{\N+1})},\\[2mm]
&\|\nabla(K(\rho)-K(\bar{\rho}))\D
q)\|_{L^{1}(B^{\N-1})}\lesssim\|q\|_{L^{\infty}(B^{\N})}\|q\|_{L^{1}(B^{\N+2})},\\[2mm]
&\|[\frac{P_{0}^{'}(\rho)}{\rho}-\frac{P_{0}^{'}(\bar{\rho})}{\bar{\rho}}]\nabla
\rho\|_{L^{1}(B^{\N-1})}\lesssim\|q\|_{L^{\infty}(B^{\N-1})}\|q\|_{L^{1}(B^{\N+1})},\\[2mm]
&\;\;\;\;\;\;\;\;\;\;\;\;\;\;\;\;\;\;\;\;\;\;\;\;\;\;\;\;\;\;\;\;\;\;\;\;\;\;\;\;\;\;\;\;\;\;\;
\lesssim\|q\|_{L^{\infty}(\widetilde{B}^{\N-1,\N})}\|q\|_{L^{1}(\widetilde{B}^{\N+1,\N+2})}\;.\\
\end{aligned}
$$
After it remains two terms to treat:
$$
\begin{aligned}
&\|\biggl[\frac{\w
P_{1}^{'}(\rho)}{\rho}-\frac{\bar{T}P_{1}^{'}(\bar{\rho})}{\bar{\rho}}\biggl]\nabla
\rho\|_{L^{1}(B^{\N-1})}
\lesssim\;\|\big[(\frac{P_{1}^{'}(\rho)}{\rho}-\frac{P_{1}^{'}(\bar{\rho})}{\bar{\rho}})\,\bar{\theta}\big]\nabla
q\|_{L^{1}(B^{\N-1})}\\[2mm]
&\hspace{3,4cm}+\|\frac{P_{1}^{'}(\bar{\rho})}{\bar{\rho}}{\cal T}\n
q\|_{L^{1}(B^{\N-1})}+\|{\cal
T}\big(\frac{P_{1}^{'}(\rho)}{\rho}-\frac{P_{1}^{'}(\bar{\rho})}{\bar{\rho}}\big)\nabla
q\|_{L^{1}(B^{\N-1})},\\[2mm]
\end{aligned}
$$
$$
\begin{aligned}
\|\biggl[\frac{\w
P_{1}^{'}(\rho)}{\rho}-\frac{\bar{T}P_{1}^{'}(\bar{\rho})}{\bar{\rho}}\biggl]\nabla
\rho\|_{L^{1}(B^{\N-1})}&\lesssim\|q\|_{L^{\infty}(B^{\N-1})}\|q\|_{L^{1}(B^{\N+1})}+\|{\cal
T}\nabla q\|_{L^{1}(B^{\N-1})}\\
&\hspace{3,8cm}+\|K_{1}(q){\cal T}\nabla
q\|_{L^{1}(B^{\N-1})},\\
\end{aligned}
$$
According to proposition
\ref{fcomposition}, we have:
$$
\begin{aligned}
&\|{\cal T}\nabla
q\|_{L^{1}(B^{\N-1})}\leq\|{\cal T}\|_{L^{2}(\widetilde{B}^{\N,\N-1})}\|q\|_{L^{2}(\widetilde{B}^{\N,\N+1})}\;,\\[2mm]
&\|K_{1}(q){\cal T}\n q\|_{L^{1}(B^{\N-1})}\leq
C\|q\|_{L^{\infty}(B^{\N})}\|{\cal T}\n q\|_{L^{1}(B^{\N-1})}\;.
\end{aligned}
$$
Therefore:
$$
\begin{aligned}\|\big[\frac{\w P_{1}^{'}(\rho)}{\rho}-\frac{\bar{T}P_{1}^{'}(\bar{\rho})}{\bar{\rho}}\big]\nabla
\rho&\|_{L^{1}(B^{\N-1})}\lesssim\|q\|_{L^{\infty}(\widetilde{B}^{\N-1,\N})}\|q\|_{L^{1}(B^{\N+1})}\\[1,5mm]
&\hspace{1cm}+(1+\|q\|_{L^{\infty}(\widetilde{B}^{\N-1,\N})})
\|{\cal T}\|_{L^{2}(\widetilde{B}^{\N,\N-1})}\|q\|_{L^{2}(\widetilde{B}^{\N,\N+1})}\,.\\
\end{aligned}$$
In the same spirit:
$$\|(\frac{P_{1}(\rho)}{\rho}-\frac{P_{1}(\bar{\rho})}{\bar{\rho}})\n\theta\|_{L^{1}(B^{\N-1})}
\lesssim\|q\|_{L^{\infty}(\widetilde{B}^{\N-1,\N})}\|{\cal
T}\|_{L^{1}(\widetilde{B}^{\N+1,\N})}\;,$$
$$
\begin{aligned}
\|(\frac{P_{1}(\rho)}{A\rho}-\frac{P_{1}(\bar{\rho})}{A\bar{\rho}})\n\theta\|_{L^{1}(B^{\N-1})}\leq&
\|q\|_{L^{2}(\widetilde{B}^{\N,\N+1})}\|{\cal
T}_{BF}\|_{L^{2}(B^{\N})}\\
&\hspace{3cm}+\|q\|_{L^{\infty}(\widetilde{B}^{\N-1,\N})}\|{\cal
T}_{HF}\|_{L^{1}(B^{\N})}\;,
\end{aligned}
$$
\\
where we have:
$${\cal T}_{BF}=\sum_{l\leq0}\D_{l}{\cal T}\;\;\mbox{and}\;\;{\cal T}_{HF}=\sum_{l>0}\D_{l}{\cal T}\;.$$
Next we have the following term:
$$\|u^{*}.\n u\|_{L^{1}(B^{\N-1})}\lesssim\|u\|_{L^{\infty}(B^{\N-1})}\|u\|_{L^{1}(B^{\N+1})}.$$
And finally we have the terms coming from  ${\rm div}(D)$ which are
of the form:
$$\|\frac{\lambda^{'}(\rho)\n \rho\,
{\rm div}u}{\rho}\|_{L^{1}(B^{\N-1})}\leq\|L(q)\n\rho\,{\rm div} u
\|_{L^{1}(B^{\N-1})}+\|\frac{\lambda
^{'}(\bar{\rho})}{\bar{\rho}}\nabla\rho\, {\rm
div}u\|_{L^{1}(B^{\N-1})}$$ where we have set:
$$L(x_{1})=\frac{\lambda^{'}(\bar{\rho}(1+x_{1}))}{\bar{\rho}(1+x_{1})}-\frac{\lambda^{'}(\bar{\rho})}{\bar{\rho}}.$$
Afterwards we  can apply proposition \ref{fcomposition} to get:
$$
\begin{aligned}
&\|\n\rho\,{\rm div}u \|_{L^{1}(B^{\N-1})}\lesssim\|u \|_{L^{1}(B^{\N+1})}\|q\|_{L^{\infty}(B^{\N})}.\\[2mm]
&\|L(q)\n\rho{\rm div} u
\|_{L^{1}(B^{\N-1})}\leq\|L(q)\|_{L^{\infty}(B^{\N})}
\|\n\rho\,{\rm div} u \|_{L^{1}(B^{\N-1})}.\\
\end{aligned}
$$
As we assumed that $(\cal{H})$ is satisfied, we have in using proposition \ref{fcomposition}:
$$\|L(q)\|_{L^{\infty}(B^{\N})}\leq
C\|q\|_{L^{\infty}(B^{\N})}\;.$$ So we have:
$$\|\frac{\lambda^{'}(\rho)\n \rho\,
{\rm div}u}{\rho}\|_{L^{1}(B^{\N-1})}\lesssim\|u
\|_{L^{1}(B^{\N+1})}\|q\|_{L^{\infty}(\widetilde{B}^{\N-1,\N})}(1+\|q\|_{L^{\infty}(\widetilde{B}^{\N-1,\N})}).$$
In the same way we have in using \ref{fproduit},
\ref{fcomposition} and \ref{fimportant}:
$$\|\frac{(du+\n u)\n\rho\;
\mu^{'}(\rho)}{\rho}\|_{L^{1}(B^{\N-1})}\leq
C\|u\|_{L^{1}(B^{\N+1})}\|q\|_{L^{\infty}(B^{\N})}(1+\|q\|_{L^{\infty}(B^{\N})}).$$
$$
\begin{aligned}
\|\nabla(\frac{K^{'}_{\rho}}{2}|\n\rho|
^{2})\|_{L^{1}(B^{\N-1})}&\lesssim
\|(\frac{K^{'}_{\rho}}{2}-\frac{K^{'}_{\bar{\rho}}}{2})|\nabla\rho|
^{2})\|_{L^{1}(B^{\N})} +\|\frac{K^{'}_{\bar{\rho}}}{2}|\nabla\rho|
^{2}\|_{L^{1}(B^{\N})},\\[2mm]
&\lesssim\|L(q)\|_{L^{\infty}(B^{\N})}\||\n\rho|
^{2}\|_{L^{1}(B^{\N})}+\||\nabla\rho|
^{2}\|_{L^{1}(B^{\N})},\\[2mm]
&\lesssim\|q\|_{L^{\infty}(B^{\N})}\|\nabla\rho\|_{L^{2}(B^{\N})}^{2}+\|\nabla\rho\|_{L^{2}(B^{\N})}^{2},\\[2mm]
&\lesssim\|q\|_{L^{\infty}(B^{\N})}\|q\|_{L^{2}(B^{\N+1})}^{2}+\|q\|_{L^{2}(B^{\N+1})}^{2}.\\
\end{aligned}
$$
where $L(q)=\frac{K^{'}_{\rho}}{2}-\frac{K^{'}_{\bar{\rho}}}{2}$ .\\
\\
3) Let us finally estimate $\|H(q,u,{\cal
T})\|_{L^{1}(\widetilde{B}^{\N-1,\N-2})}$:
$$
\begin{aligned}
\|\frac{{\rm
div}(\chi(\rho)\nabla\theta)}{\rho}-\frac{\bar{\chi}}{\bar{\rho}}\D\theta\|_{L^{1}(\widetilde{B}^{\N-1,\N-2})}
\leq\,&\|K(q){\rm
div}(K_{1}(q)\n\theta)\|_{L^{1}(\widetilde{B}^{\N-1,\N-2})}\\[3mm]
\hspace{5cm}+\|{\rm
div}(K_{1}(q)\n\theta)\|&_{L^{1}(\widetilde{B}^{\N-1,\N-2})}+\|K(q)\D\theta\|_{L^{1}(\widetilde{B}^{\N-1,\N-2})},\\
\end{aligned}
$$
and we have:
$$
\begin{aligned}
\|{\rm div}(
K_{1}(q)\n\theta)\|_{L^{1}(\widetilde{B}^{\N-1,\N-2})}\leq
C\|q\|_{L^{\infty}(B^{\N})}
\|{\cal T}\|_{L^{1}(\widetilde{B}^{\N+1,\N})}.\\
\end{aligned}
$$
So finally:
$$
\begin{aligned}
\|\frac{{\rm
div}(\chi(\rho)\nabla\theta)}{\rho}-\frac{\bar{\chi}}{\bar{\rho}}
\D\theta\|_{L^{1}(\widetilde{B}^{\N-1,\N-2})} \lesssim\,\,&
\|q\|_{L^{\infty}(\widetilde{B}^{\N-1,\N})}
\|{\cal T}\|_{L^{1}(\widetilde{B}^{\N+1,\N})}\\[3mm]
&\hspace{4cm}\times(2+\|q\|_{L^{\infty}(\widetilde{B}^{\N-1,\N})})\;.\\
\end{aligned}
$$
Next we have:
$$
\begin{aligned}
&\|(\frac{\theta
P_{1}(\rho)}{A\rho}-\frac{\bar{\theta}P_{1}(\bar{\rho})}{A\bar{\rho}}){\rm
div}\,u\|_{L^{1}(\widetilde{B}^{\N-1,\N-2})}
\lesssim\;\|{\cal T}{\rm div} u\|_{L^{1}(\widetilde{B}^{\N-1,\N-2})}\\[2mm]
&\hspace{5,3cm}+ \|{\cal T}L_{1}(q){\rm div}
u\|_{L^{1}(\widetilde{B}^{\N-1,\N-2})}+\|L_{1}(q){\rm div}
u\|_{L^{1}(\widetilde{B}^{\N-1,\N})},\\
\end{aligned}
$$
where we denote:
$$L_{1}(x)=\frac{P_{1}(\bar{\rho}(1+x))}{\bar{\rho}(1+x)}-\frac{P_{1}(\bar{\rho})}{\bar{\rho}}\;.$$
On one hand,
$$
\begin{aligned}
&\|{\cal T}{\rm div}
u\|_{L^{1}(\widetilde{B}^{\N-1,\N-2})}\lesssim\|{\cal
T}\|_{L^{\infty}(\widetilde{B}^{\N-1,\N-2})}\|u\|_{L^{1}(B^{\N+1})},\\
&\|L_{1}(q){\rm div} u\|_{L^{1}(\widetilde{B}^{\N-1,\N-2})}
\lesssim\|L_{1}(q)\|_{L^{\infty}(\widetilde{B}^{\N-1,\N-2})}\|u\|_{L^{1}(B^{\N+1})},
\end{aligned}
$$
whence the desired result:
$$
\begin{aligned}
\|(\frac{\theta
P_{1}(\rho)}{A\rho}-\frac{\bar{\theta}P_{1}(\bar{\rho})}{A\bar{\rho}}){\rm
div}\,u&\|_{L^{1}(\widetilde{B}^{\N-1,\N-2})} \lesssim\,\,
\|q\|_{L^{\infty}(\widetilde{B}^{\N-1,\N-2})}\|u\|_{L^{1}(B^{\N+1})}\\[2mm]
&\hspace{0,8cm}+\|{\cal
T}\|_{L^{\infty}(\widetilde{B}^{\N-1,\N-2})}\|u\|_{L^{1}(B^{\N+1})}(1+\|q\|_{L^{\infty}(\widetilde{B}^{\N-1,\N-2})}).\\
\end{aligned}
$$
We proceed in the same way for the others terms which are similar,
and we finish with the last two following terms:
$$
\begin{aligned}
&\|u^{*}.\n\theta\|_{
L^{1}(\widetilde{B}^{\N-1,\N-2})}\lesssim\|{\cal T}\|_{
L^{1}(\widetilde{B}^{\N+1,\N})}\|u\|_{L^{\infty}(B^{\N-1})},\\[2mm]
&\|\frac{D:\nabla u}{\rho}\|_{
L^{1}(\widetilde{B}^{\N-1,\N-2})}\lesssim\|K(q)\nabla u:\nabla u\|_{
L^{1}(\widetilde{B}^{\N-1,\N-2})}+\|\nabla u:\nabla u\|_{
L^{1}(\widetilde{B}^{\N-1,\N-2})}.\\
\end{aligned}
$$
and:
$$\|K(q)\nabla u:\nabla u\|_{
L^{1}(\widetilde{B}^{\N-1,\N-2})}\lesssim\|q\|_{L^{\infty}(B^{\N})}\|u\|_{L^{2}(B^{\N})}.$$
so the result:
$$\|\frac{D:\nabla u}{\rho}\|_{
L^{1}(\widetilde{B}^{\N-1,\N-2})}\lesssim(1+\|q\|_{L^{\infty}(\widetilde{B}^{\N-1,\N})})\|u\|_{L^{2}(B^{\N})}.$$
Finally in using (\ref{f1a5}), (\ref{f1a6}) and all the previous bound,
we get:
\begin{equation}
\|\psi(q,u,{\cal T})\|_{E^{\N}}\leq C((C+1)\eta+R)^{2}. \label{f1a7}
\end{equation}
\\
Let $c$ be such that $\|\cdot\|_{B^{\N}}\leq c$ implies that:
$\|\cdot
\|_{ L^{\infty}}\leq1/3$. Then we choose R and $\eta$ such that:
$$R\leq \inf((3C)^{-1},c,1),\;\mbox{and}\;\;\eta\leq \frac{\inf(R,c)}{C+1}\;.$$
So $(\cal{H})$ is verified and we have:
$$\psi(B(0,R))\subset B(0,R)\;.$$
\subsubsection*{2) Second step: Property of contraction}
We consider $(q^{'}_{1},u^{'}_{1},{\cal T}^{'}_{1}),\,
(q^{'}_{2},u^{'}_{2},{\cal T}^{'}_{2})$ in $B(0,R)$  where we
note:
$$\theta_{i}={\cal
T}_{i}+\bar{\theta}\;,\;\widetilde{T}_{i}=\Psi^{-1}(\theta_{i})$$
and we set:
$$(\delta q=q^{'}_{2}-q^{'}_{1},\;\de u=u^{'}_{2}-u^{'}_{1},\;\de {\cal T}={\cal T}^{'}_{2}-{\cal T}^{'}_{1})\;.$$
We have according to
proposition \ref{flinear1} and (\ref{f1a1}):
\begin{equation}
\begin{aligned}
&\|\psi_{(q_{L},u_{L},{\cal T}_{L})}(q^{'}_{2},u^{'}_{2},{\cal
T}^{'}_{2})-\psi_{(q_{L},u_{L}, {\cal
T}_{L})}(q^{'}_{1},u^{'}_{1},{\cal
T}^{'}_{1})\|_{E^{\N}}\lesssim\\[2mm]
&\hspace{2cm}\|F(q_{2},u_{2},{\cal T}_{2})-F(q_{1},u_{1},{\cal
T}_{1})\|_{L^{1}(\widetilde{B}^{\N-1,\N})}\\
&\hspace{3,5cm}+\|G(q_{2},u_{2},{\cal
T}_{2})-G(q_{1},u_{1},{\cal T}_{1})\|_{L^{1}(B^{\N-1})}\\
&\hspace{5cm}+\|H(q_{2},u_{2},{\cal T}_{2})-H(q_{1},u_{1},{\cal
T}_{1})\|_{L^{1}(\widetilde{B}^{\N-1
,\N-2})}.\\
\end{aligned}
\label{f1a8}
\end{equation}
where we have:
$$
\begin{aligned}
&F(q_{2},u_{2},{\cal T}_{2})-F(q_{1},u_{1},{\cal T}_{1})=-{\rm div}(q_{2}u_{2})+{\rm div}(q_{1}u_{1})\\[5mm]
&G(q_{2},u_{2},{\cal T}_{2})-G(q_{1},u_{1},{\cal T}_{1})=
\delta u^{*}.\n u_{2}+u_{1}^{*}.\n\delta u+\n(\frac{1}{2}(K^{'}_{\rho_{2}}-K^{'}_{\rho_{1}})|\n\rho_{2}|^{2})\\[3mm]
&+\n(\frac{1}{2}K^{'}_{\rho_{1}}(|\n\rho_{2}|^{2}-|\n\rho_{1}|^{2})
-\mu(\bar{\rho})\D\delta u+\frac{\mu(\rho_{1})}{\rho_{1}}\D\de u-\zeta(\bar{\rho})\,\n{\rm div}(\de u)\\[3mm]
&+\n(K(\rho_{1})\D(\bar{\rho}\de
q))+\n((K(\rho_{2})-K(\rho_{1}))\D\rho_{2})+[\frac{P^{'}_{0}(\rho_{2})}
{\rho_{2}}-\frac{P^{'}_{0}(\rho_{1})}{\rho_{1}}]\bar{\rho}\,\n \de q\\[3mm]
\end{aligned}
$$
$$
\begin{aligned}
&+(\frac{P^{'}_{0}(\rho_{1})}{\rho_{1}}
-\frac{P^{'}_{0}(\bar{\rho})}{\bar{\rho}})\n(\de
q)+\frac{P^{'}_{0}(\bar{\rho})}{\bar{\rho}}\n(\de q)
+\bar{T}P^{'}_{1}(\bar{\rho})\n\de q+\bar{\rho}A\,\de{\cal T}P^{'}_{1}(\rho_{2})\\[3mm]
&+\bar{\rho}\,A\theta_{1}(P^{'}_{1}(\rho_{2})-P^{'}_{1}(\rho_{1}))\n
q_{2}+\frac{1}{A}[P_{1}(\rho_{2})-P_{1}(\rho_{1})]\n\theta_{2}+\frac{1}{A}[P_{1}(\rho_{1})-P_{1}(\bar{\rho})]\n\de{\cal
T}
\\[3mm]
&+\bar{\rho}\,A\theta_{1}P^{'}_{1}(\rho_{1})\n\de
q+(\frac{\lambda^{'}(\rho_{2})}{\rho_{2}}-\frac{\lambda^{'}(\rho_{1})}{\rho_{1}})\n\rho_{2}\,{\rm
div} u_{2}
+\frac{\lambda^{'}(\rho_{1})}{\rho_{1}}\,(\bar{\rho}\,\n\de\, q{\rm div} u_{2}\\[3mm]
&+\n\rho_{1}{\rm div}\,\de
u)+(\frac{\mu^{'}(\rho_{2})}{\rho_{2}}-\frac{\mu^{'}(\rho_{1})}{\rho_{1}})\,(du_{2}+\n
u_{2})\n\rho_{2}
+\frac{\mu^{'}(\rho_{1})}{\rho_{1}}(d(\de u)+\n \de u\n\rho_{1}\\[3mm]
&+\bar{\rho}\,\n(\de q)(du_{2}+\n u_{2})).\\
\end{aligned}
$$
\\
And we have for the part pertaining to $H$:
$$
\begin{aligned}
&H(q_{2},u_{2},{\cal T}_{2})-H(q_{1},u_{1},{\cal T}_{1})=
(\frac{1}{\rho_{2}}-\frac{1}{\rho_{1}})\,{\rm
div}(\chi(\rho_{2})\n\theta_{2})+\frac{1}{\rho_{1}}{\rm
div}((\chi(\rho_{2})-\chi(\rho_{1}))\n\theta_{2})\\[3mm]
&+\frac{1}{\rho_{1}}\,{\rm div}(\chi(\rho_{1})\n\de{\cal T})+\biggl[
\frac{P_{1}(\rho_{1})}{\rho_{1}}-\frac{P_{1}(\rho_{2})}{\rho_{2}}\biggl]\frac{\theta_{1}}{A}{\rm
div }u_{1}+\frac{P_{1}(\rho_{2})}{\rho_{2}}\,\frac{\de{\cal
T}}{A}{\rm div} u_{1} -\de
u^{*}.\n\theta_{2}\\[3mm]
&+\frac{P_{1}(\rho_{2})}{\rho_{2}}\,\frac{\theta_{2}}{A}{\rm div}
\de u -u_{1}^{*}\n\de{\cal T}+(\frac{1}{\rho_{2}}-
\frac{1}{\rho_{1}})D_{2}:\n u_{2}-\frac{1}{\rho_{1}}D_{1}:\n\de u-\frac{1}{\rho_{1}}(D_{2}-D_{1}):\n u_{2}\\[2mm]
\end{aligned}
$$
Let us first estimate $\|F(q_{2},u_{2},{\cal
T}_{2})-F(q_{1},u_{1},{\cal
T}_{1})\|_{L^{1}(\widetilde{B}^{\N-1,\N})}$. We have:
$$
\begin{aligned}
&\|F(q_{2},u_{2},{\cal{T}}_{2})-F(q_{1},u_{1},{\cal
T}_{1})\|_{L^{1}(\widetilde{B}^{\N-1,\N})}\leq
\|{\rm div}((q_{2}-q_{1})u_{2})\|_{L^{1}(\widetilde{B}^{\N-1,\N})}\\[2mm]
&\hspace{8,7cm}+\|{\rm div}(q_{1}(u_{2}-u_{1}))\|_{L^{1}(\widetilde{B}^{\N-1,\N})},\\[4mm]
&\lesssim\|\delta
q\|_{L^{2}(B^{\N})}\|u_{2}\|_{L^{2}(B^{\N})}+\|\delta
q\|_{L^{\infty}(B^{\N})}\|u_{2}\| _{L^{1}(B^{\N+1})}+\|\delta
q\|_{L^{2}(B^{\N+1})}\|u_{2}\|_{L^{2}(B^{\N})}\\[3mm]
&\;\;\;\;+\|q_{1}\|_{L^{2}(B^{\N})}\|\delta u\|_{L^{2}(B^{\N})}+\|
q_{1}\|_{L^{\infty}(B^{\N})}
\|\delta u\|_{L^{1}(B^{\N+1})}+\|q_{1}\|_{L^{2}(B^{\N+1})}\|\delta u\|_{L^{2}(B^{\N})}.\\
\end{aligned}
$$
\\
Next, we have to bound $\|G(q_{2},u_{2},{\cal
T}_{2})-G(q_{1},u_{1},{\cal T}_{1})\|_{L^{1}(B^{\N-1})}$. We treat
only one typical term, the others are of the same form.\\
We use essentially the proposition \ref{fcomposition} to treat the
product and the composition, so we get :
$$
\|\frac{\mu(\rho_{1})}{\rho_{1}}\D(u_{2}-u_{1})\|_{L^{1}(B^{\N-1})}\lesssim
(1+\|q_{1}\|_{L^{\infty}(B^{\N})})\|\delta u\|_{L^{1}(B^{\N+1})}.
$$
Bounding $\|H(q_{2},u_{2},{\cal T}_{2})-H(q_{1},u_{1},{\cal
T}_{1})\|_{L^{1}(\widetilde{B}^{\N-1,\N-2})}$ is left to the reader.
So we get in using the proposition \ref{flinear1} :
$$
\begin{aligned}
\|\Psi(q^{'}_{2},u^{'}_{2},{\cal T}^{'}_{2})-
\Psi(q^{'}_{1},u^{'}_{1},{\cal T}^{'}_{1})\|_{E^{\N}} \leq
&C\,\|(\de q,\de u,\de {\cal
T})\|_{E^{\N}}\,\biggl(\|(q^{'}_{1},u^{'}_{1},{\cal
T}^{'}_{1})\|_{E^{\N}}\\
&\hspace{1,5cm}+\|(q^{'}_{2},u^{'}_{2},{\cal T}^{'}_{2})\|_{E^{\N}}
+2\|(q_{L},u_{L},{\cal T}_{L})\|_{E^{\N}}\biggl).\\
\end{aligned}
$$
If one chooses $R$ small enough, we end up with in using (\ref{f1a8}) and the previous estimates:
$$\|\Psi(q^{'}_{2},u^{'}_{2},{\cal T}^{'}_{2})-\Psi(q^{'}_{1},u^{'}_{1},{\cal T}^{'}_{1})\|_{E^{\N}}
\leq\frac{3}{4}\,\|(\de q,\de u,\de {\cal T})\|_{E^{\N}}.$$ We thus
have the property of contraction and so by the fixed point theorem,
we have existence of a solution to $(NHV)$. Indeed we can see easily that
$E^{\N}$ is a Banach space.
 \subsubsection*{3)Uniqueness of the solution:}
The proof is similar to the proof of contraction, hence we will have
the same type of estimates. So consider two solutions in
$\widetilde{E}^{\N}$: $(q_{1},u_{1},{\cal T}_{1})$ and
$(q_{2},u_{2},{\cal T}_{2})$ of the system $(NHV)$ with the same
initial data. With no loss of generality, one can assume that
$(q_{1},u_{1},{\cal T}_{1})$ is the solution found in the previous
section.
\\
We thus have:
$$\|q_{1}\|_{L^{\infty}([0,T]\times\R^{N})}\leq\frac{1}{2}.
\leqno{\cal{(H)}}$$ Let $\bar{T}$ be the largest time such that
$q_{2}$ verifies $\cal{(H)}$. By continuity, we have $0<\bar{T}\leq
T$. Next we
see that:\\
$$\delta q=q_{2}-q_{1},\;\delta u=u_{2}-u_{1},\;\delta
{\cal T}={\cal T}_{2}-{\cal T}_{1}$$ verifies the system:\\
$$
\begin{cases}
\begin{aligned}
&\p_{t}\delta q+{\rm div}\delta u=F(q_{2},u_{2},{\cal T}_{2})-F(q_{1},u_{1},{\cal T}_{1}),\\[3mm]
&\p_{t}\delta u-\frac{\bar{\mu}}{\bar{\rho}}\Delta \delta
u-\frac{\bar{\zeta}}{\bar{\rho}}\nabla {\rm div}\delta
u-\bar{\rho}\bar{K}\nabla\Delta \delta
q+(P_{0}^{'}(\bar{\rho})+\bar{T}P_{1}^{'}(\bar{\rho}))\nabla \delta
q+\frac{P_{1}(\bar{\rho})}{\bar{\rho}\psi^{'}(\bar{T})}\nabla\delta{\cal
T}\\[2mm]
&\hspace{9cm}=G(q_{2},u_{2},{\cal T}_{2})-G(q_{1},u_{1},{\cal T}_{1}),\\[3mm]
&\p_{t}\delta{\cal T}-\frac{\bar{\chi}}{\bar{\rho}}\Delta\delta{\cal
T}+\frac{\bar{T}P_{1}(\bar{\rho})}{\bar{\rho}}{\rm
div}\delta u§=§H(q_{2},u_{2},{\cal T}_{2})-H(q_{1},u_{1},{\cal T}_{1}).\\
\end{aligned}
\end{cases}
$$
\\
We apply the proposition \ref{flinear1} on $[0,T_{1}]$
with $0<T_{1}\leq\bar{T}$ and we have:
$$\|(\delta q,\delta u,\delta{\cal T})\|_{\widetilde{E}^{\N}}\leq
A(T_{1})\|(\delta q,\delta u,\delta{\cal
T})\|_{\widetilde{E}^{\N}}$$
where we have for $T_{1}$ enough small $A(T_{1})\leq\frac{1}{2}$.\\
\\
And we thus have: $\de q=0$, $\de u=0$, $\de {\cal T}=0$ on
$[0,T_{1}]$ for $T_{1}$ small enough and we conclude after by
connectivity.
\hfill {$\Box$}
\\
\\
We treat now the specific case of $N=2$, where we need  more
regularity for the initial data because we cannot use the
proposition \ref{fproduit} in the case $N=2$ with the previous
initial data. Indeed we cannot treat some non-linear terms such as
$\|{\cal T}{\rm div}u\|_{L^{1}(\widetilde{B}^{0,-1})}$ or
$\|u^{*}.\n\theta\|_{L^{1}(\widetilde{B}^{0,-1})}$ because if we
want to use proposition \ref{fproduit}, we are in the case
$s_{1}+s_{2}=0$. This is the reason why more regularity is required.\\
\\
We recall the space in which we are working:
$$
\begin{aligned}
E^{'}=&[C_{b}(\R_{+},\widetilde{B}^{0,1+\e^{'}})\cap
L^{1}(\R_{+},\widetilde{B}^{2,3+\e^{'}})]\times[C_{b}(\R_{+},\widetilde{B}^{0,\e^{'}})^{N}\cap
L^{1}(\R_{+},\widetilde{B}^{2,2+\e^{'}})^{N}]\\
&\times[C_{b}(\R_{+},\widetilde{B}^{0,-1+\e^{'}})\cap
L^{1}(\R_{+},\widetilde{B}^{2,1+\e^{'}})]\\
\end{aligned}
$$
with $\e^{'}>0$, $E^{'}$ being the space in which we have a solution
. And $\widetilde{E}^{'}$ corresponds to the space where we show the
uniqueness of solution.
$$
\begin{aligned}
\widetilde{E}^{'}=&[C_{b}(\R_{+},\widetilde{B}^{0,1+\e^{'}})\cap
L^{2}(\R_{+},\widetilde{B}^{1,2+\e^{'}})]\times[C_{b}(\R_{+},\widetilde{B}^{0,\e^{'}})^{N}\cap
L^{2}(\R_{+},\widetilde{B}^{1,1+\e^{'}})^{N}]\\
&\times[C_{b}(\R_{+},\widetilde{B}^{0,-1+\e^{'}})\cap
L^{2}(\R_{+},\widetilde{B}^{1,\e^{'}})].\\
\end{aligned}
$$
\\
{\bf Proof of  theorem \ref{ftheo2}}
\\
\\
The proof is similar to the previous one except that we have changed
the functional space, in which the fixed point theorem is applied.
So we want verify that the function $\psi$ is contracting to apply
the fixed point. We denote by $(q_{L},u_{L},{\cal T}_{L})$ the
solution of the linear system $(M^{'})$ with $F=G=H=0$ and with
initial data $(q_{0},u_{0},{\cal T}_{0})$\\
Arguing as before, we get:
\begin{equation}
\begin{aligned}
&\|\psi(q,u,{\cal T})\|_{E^{'}}\leq\,C\eta+
\|F(q,u,{\cal T})\|_{ L^{1}(\widetilde{B}^{0,1+\e^{'}})}\\
&\hspace{4cm}+\|G(q,u,{\cal T})\|_{ L^{1}(\widetilde{B}^{0,\e^{'}})}
+\| H(q,u,{\cal T})\|_{L^{1}(\widetilde{B}^{0,-1+\e^{'}})}\;.\\
\end{aligned}
\label{fb6}
\end{equation}
if:
$$\|q_{0}\|_{\widetilde{B}^{0,1+\e^{'}}}+\|u_{0}\|_{\widetilde{B}^{0,\e^{'}}}+\|{\cal T}_{0}\|_{\widetilde{B}^{0,-1+\e^{'}}}\leq\eta.$$
Let us estimate $\|F(q,u,{\cal
T})\|_{ L^{1}(\widetilde{B}^{0,1+\e^{'}})}$, $\|G(q,u,{\cal T})\|_{
L^{1}(\widetilde{B}^{0,\e^{'}})}$ and $\| H(q,u,{\cal
T})\|_{L^{1}(\widetilde{B}^{0,-1+\e^{'}})}$, we just give two
examples of estimates in the space $E^{'}$,the other estimates are
left to the reader.
$$\|{\rm div}(qu)\|_{ L^{1}(\widetilde{B}^{0,1+\e^{'}})}\leq\|qu\|_{L^{1}(B^{1})}+\|qu\|_{L^{1}(B^{2+\e^{'}})},$$
and:
$$
\begin{aligned}
&\|qu\|_{L^{1}(B^{1})}\lesssim\|q\|_{L^{2}(B^{1})}\|u\|_{L^{2}(B^{1})},\\[2mm]
&\|qu\|_{L^{1}(B^{2+\e^{'}})}\lesssim\|q\|_{L^{\infty}(B^{1})}\|u\|_{L^{1}(B^{2+\e^{'}})}+
\|q\|_{L^{2}(B^{2+\e^{'}})}\|u\|_{L^{2}(B^{1})}.\\
\end{aligned}
$$
We do similarly for $\|G(q,u,{\cal
T})\|_{L^{1}(\widetilde{B}^{0,\e^{'}})}$. The new difficulty
appears on the last term $\|H(q,u,{\cal
T})\|_{L^{1}(\widetilde{B}^{0,-1+\e^{'}})}$. In fact it's only for this
term that that additional regularity is needed. Proposition \ref{fproduit} enables us to write:
$$
\begin{aligned}
&\|{\cal T}{\rm div}
u\|_{L^{1}(\widetilde{B}^{0,-1+\e^{'}})}\lesssim\|{\cal
T}\|_{L^{\infty}(\widetilde{B}^{0,-1+\e^{'}})}\|u\|_{L^{1}(B^{2})}.\\[2mm]
&\|u^{*}.\n\theta\|_{
L^{1}(\widetilde{B}^{0,-1+\e^{'}})}\lesssim\|{\cal T}\|_{
L^{1}(\widetilde{B}^{1,1+\e^{'}})}\|u\|_{L^{\infty}(B^{0})}.
\end{aligned}
$$
To conclude we follow the previous proof. Uniqueness in
$\widetilde{E}^{'}$ goes along the lines of the proof of uniqueness
in dimension $N\geq 3$.
\hfill{$\Box$}
\subsection{Existence of a solution in the general case with small initial data\label{fsubsection43}}
In this section we are interested by the general case where all the
coefficients depend on the density and the temperature except
$\kappa$. In this case to control the non-linear terms we need that
$\theta$
be bounded, that's why we need to take more regular initial data to preserve the $L^{\infty}$ bound.\\
As the initial data are more regular, we need to obtain new
estimates in Besov spaces on the linear system $(M^{'})$.
\begin{proposition}
\label{flinear2} Under conditions of proposition \ref{fconditions}
with strict inequality, let $(q,u,{\cal T})$ be a solution of the
system $(M^{'})$ on $[0,T)$ with initial conditions $
(q_{0},u_{0},{\cal T}_{0}) $ such that:
$$ q_{0}\in \widetilde{B}^{s-1,s+1},u_{0}\in \widetilde{B}^{s-1,s},{\cal T}_{0}\in\widetilde{B}^{s-1,s} .$$
Moreover  we suppose $1\leq r_{1}\leq +\infty$ and:
$$F\in\widetilde{L}_{T}^{r_{1}}(\widetilde{B}^{s-3+\frac{2}{r_{1}},s-1+\frac{2}{r_{1}}}),
\;\;G\in\widetilde{L}_{T}^{r_{1}}(\widetilde{B}^{s-3+\frac{2}{r_{1}},s-2+\frac{2}{r_{1}}}),\;\;
H\in\widetilde{L}_{T}^{r_{1}}(\widetilde{B}^{s-3+\frac{2}{r_{1}},s-2+\frac{2}{r_{1}}}).$$
We then have  the following estimate for all $r\in[r_{1},+\infty]$:
$$
\begin{aligned}
&\|q\|_{\widetilde{L}_{T}^{r}(\widetilde{B}^{s-1+\frac{2}{r},s+1+\frac{2}{r}})}+\|
u\|_{\widetilde{L}_{T}^{r}(\widetilde{B}^{s-1+\frac{2}{r},s+\frac{2}{r}})}+\|
{\cal
T}\|_{\widetilde{L}_{T}^{r}(\widetilde{B}^{s-1+\frac{2}{r},s+\frac{2}{r}})}
\lesssim\| q_{0}\|_{\widetilde{B}^{s-1,s+1}}+\|
u_{0}\|_{\widetilde{B}^{s-1,s}}\\[2mm]
&+\| {\cal T}_{0}\|_{\widetilde{B}^{s-1,s}}+\|
F\|_{L_{T}^{r_{1}}(\widetilde{B}^{s-3+\frac{2}{r_{1}},s-1+\frac{2}{r_{1}}})}+\|
G\|_{
L_{T}^{r_{1}}(\widetilde{B}^{s-3+\frac{2}{r_{1}},s-2+\frac{2}{r_{1}}})}+\|
H\|_{L_{T}^{r_{1}}(\widetilde{B}^{s-3+\frac{2}{r_{1}},s-2+\frac{2}{r_{1}}})}.\\
\end{aligned}
$$
\end{proposition}
{\bf Proof:}
\\
\\
The proof is similar to that of proposition \ref{flinear1}. Low
frequencies are treated as in proposition \ref{flinear1} because we
don't change the regularity index for the low frequencies. On the
other hand in the case of high frequencies the regularity index has
changed so that we have to see what is new. For the medium
frequencies we can proceed as in  proposition \ref{flinear1}.
\subsubsection*{Case of high frequencies:}
We are going to work
with $l\geq l_{1}$ where we will determine $l_{1}$ hereafter.
We set:
$$f_{l}^{2}=\e B\|\Lambda q_{l}\|_{L^{2}}^{2}+B\|d_{l}\|_{L^{2}}^{2}+\| {\cal{T}}_{l}\|_{L^{2}}^{2}-2K\langle\Lambda q_{l},d_{l}\rangle$$
where $B$ and $K$ will be chosen later on.\\
Then we take the scalar product of $(\ref{feq3})$ with ${\cal T}_{l}$, we get:
\begin{equation}
\frac{1}{2}\frac{d}{dt}\|{\cal T}_{l}\|_{L^{2}}^{2}+\alpha\|\n{\cal
T}_{l}\|_{L^{2}}^{2}+\delta\langle
 \Lambda d_{l},{\cal T}_{l}\rangle=0. \label{feq1111}
\end{equation}
After we sum (\ref{feq9}), (\ref{feq10}) and (\ref{feq1111}) to get:
\begin{equation}
\begin{aligned}
&\frac{1}{2}\frac{d}{dt}(B\|d_{l}\|_{L^{2}}^{2}+\e B\|\Lambda
q_{l}\|_{L^{2}}^{2}+\|{\cal T}_{l}\|_{L^{2}}^{2})+B\nu\|\nabla
d_{l}\|_{L^{2}}^{2}+\alpha\| \n{\cal
T}_{l}\|_{L^{2}}^{2}\\[2mm]
&\hspace{5cm}-B\beta\langle\Lambda
q_{l},d_{l}\rangle-B\gamma\langle\Lambda {\cal
T}_{l},d_{l}\rangle+\delta\langle\La
d_{l},{\cal T}_{l}\rangle=0.\\
\end{aligned}
\label{feq1112}
\end{equation}
We sum (\ref{feq1112}) and (\ref{feq13}) and we get:
\begin{equation}
\begin{aligned}
&\frac{1}{2}\frac{d}{dt}f_{l}^{2}+[B\nu\|\nabla
d_{l}\|_{L^{2}}^{2}-K\|\Lambda d_{l}\|_{L^{2}}^{2}]+\alpha\|\n{\cal
T}_{l}\|_{L^{2}}^{2}+[\beta K\|\Lambda q_{l}\|_{L^{2}}^{2}+\e
K\|\Lambda^{2}
q_{l}\|_{L^{2}}^{2}]\\[3mm]
&-B\beta\langle\Lambda q_{l},d_{l}\rangle-B\gamma\langle\Lambda
{\cal T}_{l},d_{l}\rangle+\delta\langle \La d_{l},{\cal
T}_{l}\rangle+K\nu\langle\Delta d_{l},\Lambda q_{l}\rangle+\gamma
K\langle\Lambda {\cal T}_{l},\Lambda
q_{l}\rangle=0.\\
\end{aligned}
\label{feq14}
\end{equation}
\\
We interest us after only to the terms of high frequencies, so
arguing as in  proposition \ref{flinear1} we get:
\begin{equation}
\begin{aligned}
&\frac{1}{2}\frac{d}{dt}f_{l}^{2}+\big[B\nu-(K+B\gamma\frac{a}{2c2^{2l_{1}}}+K\nu\frac{2}{b}+B\beta\frac{d}{2c2^{2l_{1}}}
+\delta\frac{e}{2})\big] \|\Lambda
d_{l}\|_{L^{2}}^{2}\\
&+[\alpha-(B\gamma \frac{1}{2a}+\gamma
K\frac{1}{2c^{'}}+\delta\frac{1}{2e}\frac{1}{c2^{2l_{1}}})]\|\Lambda{\cal T}_{l}\|_{L^{2}}^{2}\\
&+(\frac{\beta K}{c2^{2l_{1}}}+\e K-K\nu\frac{b}{2}-\gamma
K\frac{c^{'}}{2c2^{2l_{1}}}-B\beta\frac{1}{2d}\frac{1}{2^{2l_{1}}})\|\Lambda^{2}
q_{l}\|_{L^{2}}^{2}\leq0\;.\\
\end{aligned}
\label{feq15}
\end{equation}
Let us assume that:
$$
\begin{aligned}
&B\nu-(K+B\gamma\frac{a}{2c2^{2l_{1}}}+K\nu\frac{2}{b}+B\beta\frac{d}{2c2^{2l_{1}}}
+\delta\frac{e}{2})>0,\\[2mm]
&\alpha-(B\gamma \frac{1}{2a}+\gamma
K\frac{1}{2c^{'}}+\delta\frac{1}{2e}\frac{1}{2^{2l_{1}}})>0,\\[2mm]
&\frac{\beta K}{c2^{2l_{1}}}+\e K-K\nu\frac{b}{2}-\gamma
K\frac{c^{'}}{2c2^{2l_{1}}}-B\beta\frac{1}{2d}\frac{1}{2^{2l_{1}}}>0.\\
\end{aligned}
\leqno{(1)}
$$
We recall that $\nu>0,$ and $\alpha>0$. Next we want to have:
$$
\begin{aligned}
&\e-\nu\frac{b}{2}>0.\\
\end{aligned}
$$
So we take: $b=\frac{\nu}{\e}$ (we recall that $\e>0$). So with this
choice we get $(1)$ in taking $B,\;K$ small enough and $l_{1}$ big
enough in following the same type of estimate as in the proof of the
proposition \ref{flinear1}. We have then for $l\leq l_{0}$, $l\geq
l_{1}$ and $c^{'}$ small enough:
$$\frac{1}{2}\frac{d}{dt}f_{l} ^{2}+c^{'}2^{2l}f_{l}^{2}\leq0.$$
and:
$$f_{l}\simeq \max(1,2^{l})\| q_{l}\|_{L^{2}}+\| d_{l}\|_{L^{2}}+\| {\cal T}_{l}\|_{L^{2}}$$
Next we conclude in a similar way as in proposition \ref{flinear1}.
\hfill {$\Box$}\\
\\
In the general case the coefficients depend on the temperature and
we have to control the norm $L^{\infty}$ in order to apply the
theorems of composition. This motivates us to work in the following
spaces:
$$
\begin{aligned}
F^{\N}=&[C_{b}(\R_{+},\widetilde{B}^{\N-1,\N+1})\cap
L^{1}(\R_{+}\widetilde{B}^{\N+1,\N+3})]\times[
C_{b}(\R_{+},\widetilde{B}^{\N-1,\N})\cap\\
& L^{1}(\R_{+},\widetilde{B}^{\N+1,\N+2})]^{N}\times[
C_{b}(\R_{+},\widetilde{B}^{\N-1,\N})\cap
L^{1}(\R_{+},\widetilde{B}^{\N+1,\N+2})]. \\
\end{aligned}
$$
\\
{\bf Proof of theorem \ref{ftheo3}:}
\\
\\
The principle of the proof is similar to the previous one and we use
the same notation. We define the map $\psi$  as before with the same
$F$, $G$ and $H$ except that our coefficients depends on the density
and the temperature. We will verify only that $\psi$ maps a ball
$B(0,R)$ into itself, the end is left to the reader.
\subsubsection*{1) First step, uniform Bounds:}
We set:
$$\alpha_{0}=\| q_{0}\|_{\widetilde{B}^{\N-1,\N+1}}+\| u_{0}\|_{\widetilde{B}^{\N-1,\N}}+\|{\cal T}_{0}\|_{\widetilde{B}^{\N-1,\N}}\;.$$
\\
We denote $(q_{L},u_{L},{\cal T}_{L})$ the solution of $(M^{'})$
with initial data $(q_{0},u_{0},{\cal T}_{0})$. We have so in
accordance with proposition \ref{flinear2} the following estimates:
\begin{equation}
\|(q_{L},u_{L},{\cal T}_{L})\|_{F^{\N}}\leq C\alpha_{0}\;,\label{fa5}
\end{equation}
\vspace{0.1cm}
\begin{equation}
\begin{aligned}
&\|\psi(q,u,{\cal T})\|_{F^{\N}}\leq\,C\alpha_{0}+
\|F(q,u,{\cal T})\|_{ L_{T}^{1}(\widetilde{B}^{\N-1,\N+1})}\\[2mm]
&\hspace{3cm}+\|G(q,u,{\cal T})\|_{ L_{T}^{1}(\widetilde{B}^{\N-1,\N})}+\| H(q,u,{\cal T})\|_{L_{T}^{1}(\widetilde{B}^{\N-1,\N})}\;.\\
\end{aligned}
\label{fa6}
\end{equation}
Moreover we suppose for the moment that:
$$\|q\|_{L^{\infty}(\R\times \R^{N})}\leq 1/2\;\;\mbox{and}\;\;\|{\cal T}\|_{L^{\infty}(\R\times \R^{N})}\leq 1/2\;.\leqno{(\cal{H})}$$
We will now treat each term: $\|F(q,u,{\cal T})\|_{
L_{T}^{1}(\widetilde{B}^{\N-1,\N+1})}$, $\|G(q,u,{\cal T})\|_{
L_{T}^{1}(\widetilde{B}^{\N-1,\N})}$ and\\
 $\| H(q,u,{\cal
T})\|_{L_{T}^{1}(\widetilde{B}^{\N-1,\N})}.$
\\
\\
1) We notice that:
$$
\begin{aligned}
\|{\rm div}(qu)\|_{L^{1}(B^{\N-1})}&\leq\|q\|_{L^{2}(B^{\N})}\|u\|_{L^{2}(B^{\N})},\\[2mm]
&\leq\|q\|_{L^{2}(B^{\N+2})}\|u\|_{L^{2}(B^{\N})}+\|u\|_{L^{1}(B^{\N+2})}
\|q\|_{L^{\infty}(B^{\N})}.\\
\end{aligned}
$$
2) After we focus on $\|G(q,u,{\cal
T})\|_{L^{1}(\widetilde{B}^{\N-1,\N})}$. We have according to
proposition \ref{fcomposition}:
$$
\begin{aligned}
\biggl\|[\frac{\mu(\rho,\theta)}{\rho}-\frac{\mu(\bar{\rho},\bar{\theta})}{\bar{\rho}}]\Delta
u\biggl\|
_{L^{1}(\widetilde{B}^{\N-1,\N})}&\leq\|K(q,{\cal T})\|_{L^{\infty}(B^{\N})}\|u\|_{L^{1}(\widetilde{B}^{\N+1,\N+2})}\,\\[2mm]
&\lesssim(\|q\|_{L^{\infty}(B^{\N})}+\|{\cal T}\|_{L^{\infty}(B^{\N})})\|u\|_{L^{1}(\widetilde{B}^{\N+1,\N+2})}\,.\\[2mm]
\end{aligned}
$$
We proceed in a similar way for the term:
$$(\zeta(\rho)-\zeta(\bar{\rho})) \nabla{\rm div}u.$$
Next we have in using propositions \ref{fproduit} and
\ref{fimportant}:
$$
\begin{aligned}
&\|\rho\nabla\big(\kappa(\rho)-\kappa(\bar{\rho})\big)\D
q)\|_{L^{1}(\widetilde{B}^{\N-1,\N})}
\lesssim\biggl(\|q\|_{L^{1}(\widetilde{B}^{\N+2,\N+3})}\|q\|_{L^{\infty}(\widetilde{B}^{\N-1,\N+1})}\\
&\hspace{4,7cm}+\|q\|_{L^{\infty}(\widetilde{B}^{\N-1,\N+1})}\|q\|_{L^{1}(B^{\N+2})}\biggl)
(1+\|q\|_{L^{\infty}(\widetilde{B}^{\N-1,\N+1})})\\[3mm]
&\|\big[\,\frac{\bar{\rho}\,P_{0}^{'}(\rho)}{\rho}-P_{0}^{'}(\bar{\rho})\,\big]\nabla
q\|_{L^{1}(\widetilde{B}^{\N-1,\N})}\leq\|q\|_{L^{\infty}(\widetilde{B}^{\N-1,\N})}
\|q\|_{L^{1}(B^{\N+1})}.\\
\end{aligned}
$$
Next, we have to treat the following terms:
$$
\begin{aligned}
\|\big[\,\frac{\w\bar{\rho}\,
P_{1}^{'}(\rho)}{\rho}-\bar{T}P_{1}^{'}(\bar{\rho})\,\big]&\nabla
q\|_{L^{1}(\widetilde{B}^{\N-1,\N})}
\lesssim\,\,\|L_{1}(q)L_{2}({\cal T})\n
q\|_{L^{1}(\widetilde{B}^{\N-1,\N})}\\[2mm]
&\hspace{2cm}+\|L_{1}(q)\n
q\|_{L^{1}(\widetilde{B}^{\N-1,\N})}+\|L_{2}({\cal T})\n q\|_{L^{1}(\widetilde{B}^{\N-1,\N})},\\[3mm]
\end{aligned}
$$
where $L_{1}$ and $L_{2}$ are regular function in the sense of
proposition \ref{fcomposition}. And we have:
$$
\begin{aligned}
&\|L_{1}(q)\n
q\|_{L^{1}(\widetilde{B}^{\N-1,\N})}\lesssim\|q\|_{L^{\infty}(\widetilde{B}^{\N-1,\N})}\|q\|_{L^{1}(B^{\N+1})},\\
&\|L_{2}({\cal T})\n q\|_{L^{1}(\widetilde{B}^{\N-1,\N})}\lesssim
\|{\cal T}\|_{L^{\infty}(\widetilde{B}^{\N-1,\N})}\|q\|_{L^{1}(B^{\N+1})}.\\
\end{aligned}
$$
Finally:
$$
\begin{aligned}
\|\big[\,\frac{\w\bar{\rho}\,
P_{1}^{'}(\rho)}{\rho}-\bar{T}P_{1}^{'}(\bar{\rho})\,\big]&\nabla
q\|_{L^{1}(\widetilde{B}^{\N-1,\N})}
\lesssim\;\|q\|_{L^{\infty}(\widetilde{B}^{\N-1,\N+1})}\|{\cal
T}\|_{L^{\infty}(\widetilde{B}^{\N-1,\N})}\|q\|_{L^{1}(B^{\N+1})}\\[2mm]
&\hspace{1,3cm}+\big(\|q\|_{L^{\infty}(\widetilde{B}^{\N-1,\N+1})}+\|{\cal
T}\|_{L^{\infty}(\widetilde{B}^{\N-1,\N})}
\big)\|q\|_{L^{1}(\widetilde{B}^{\N+1,\N+3})}\\
\end{aligned}
$$
and:
$$
\begin{aligned}
\|(\frac{P_{1}(\rho)}{\rho\Psi^{'}(\w)}-\frac{P_{1}(\bar{\rho})}{\bar{\rho}\,\Psi^{'}(\bar{T})})\n\theta\|_{L^{1}(\widetilde{B}^{\N-1,\N})}
\lesssim&\,\,\|L_{1}(q)\n {\cal
T}\|_{L^{1}(\widetilde{B}^{\N-1,\N})}+\|L_{2}({\cal T})\n
\theta\|_{L^{1}(\widetilde{B}^{\N-1,\N})}\\[2mm]
&\hspace{3cm}+\|L_{1}(q)L_{2}({\cal T})\n
\theta\|_{L^{1}(\widetilde{B}^{\N-1,\N})},\\
\end{aligned}
$$
$$
\begin{aligned}
&\|L_{2}({\cal T})\n
\theta\|_{L^{1}(\widetilde{B}^{\N-1,\N})}\lesssim\|{\cal
T}\|_{L^{\infty}(\widetilde{B}^{\N-1,\N})}
\|{\cal T}\|_{L^{1}(B^{\N+1})},\\[2mm]
&\|L_{1}(q)\n
\theta\|_{L^{1}(\widetilde{B}^{\N-1,\N})}\lesssim\|q\|_{L^{\infty}(\widetilde{B}^{\N-1,\N+1})}
\|{\cal T}\|_{L^{1}(\widetilde{B}^{\N+1,\N+2})}.\\
\end{aligned}
$$
Finally:
$$
\begin{aligned}
\|(\frac{P_{1}(\rho)}{\rho\Psi^{'}(T)}-\frac{P_{1}(\bar{\rho})}{\bar{\rho}\,\Psi^{'}(\bar{T})})&\n\theta\|_{L^{1}(\widetilde{B}^{\N-1,\N})}
\lesssim\;(\|q\|_{L^{\infty}(\widetilde{B}^{\N-1,\N+1})})^{2}
\|{\cal
T}\|_{L^{1}(\widetilde{B}^{\N+1,\N+2})}\\[2mm]
&\hspace{1cm}+\big(\|{\cal
T}\|_{L^{\infty}(\widetilde{B}^{\N-1,\N})}+\|q\|_{L^{\infty}(\widetilde{B}^{\N-1,\N+1})}\big)
\|{\cal T}\|_{L^{1}(B^{\N+1,\N+2})}.\\[2mm]
\end{aligned}
$$
After we have the following terms:
$$
\|u^{*}.\n
u\|_{L^{1}(\widetilde{B}^{\N-1,\N})}\leq\|u\|_{L^{2}(\widetilde{B}^{\N,\N+1})}^{2}.$$
And  we have the terms coming from  ${\rm div}(D)$. We will treat
this one:
$$
\begin{aligned}
\|\,\frac{\lambda^{'}_{1}(\rho,\theta)\n \rho\,\, {\rm
div}u}{\rho}\,\|_{L^{1}(\widetilde{B}^{\N-1,\N})}\lesssim&\;\|L(q,{\cal
T})\n\rho\,{\rm div} u
\|_{L^{1}(\widetilde{B}^{\N-1,\N})}+\|\nabla\rho\, {\rm
div}u\|_{L^{1}(\widetilde{B}^{\N-1,\N})}\\[2mm]
\lesssim(1+\|q\|_{L^{\infty}(B^{\N-1,\N+1})}&+\|{\cal
T}\|_{L^{\infty}(B^{\N-1,\N})})\,
\|q\|_{L^{\infty}(\widetilde{B}^{\N-1,\N+1})}\|u\|_{L^{1}(\widetilde{B}^{\N+1,\N+2})}.\\
\end{aligned}
$$
Afterwards in the same way we can treat the terms of the type:
$$\frac{(du+\n u)\n\rho\;
\mu^{'}_{1}(\rho,\theta)}{\rho},\;\frac{(du+\n u)\n\theta\;
\mu^{'}_{2}(\rho,\theta)}{\rho}\;\;\mbox{and}\;\;\frac{\lambda^{'}_{2}(\rho,\theta)\n
\theta\, {\rm div}u}{\rho}.$$
Finally, we have:
$$
\begin{aligned}
\|\nabla(K^{'}_{\rho}|\n\rho|
^{2})\|_{L^{1}(\widetilde{B}^{\N-1,\N})}&\lesssim
\|(K^{'}_{\rho}-K^{'}_{\bar{\rho}})|\nabla\rho|
^{2}\|_{L^{1}(\widetilde{B}^{\N,\N+1})} +\||\nabla\rho|
^{2}\|_{L^{1}(\widetilde{B}^{\N,\N+1})}\\[2mm]
&\lesssim(1+\|q\|_{L^{\infty}(\widetilde{B}^{\N-1,\N+1})})\|q\|^{2}_{L^{2}(\widetilde{B}^{\N+1,\N+2})}.\\
\end{aligned}
$$
\\
3) Let us finally estimate $\|H(q,u,{\cal
T})\|_{L^{1}(\widetilde{B}^{\N-1,\N})}$:
$$
\begin{aligned}
&\|\,\frac{{\rm div}(\chi(\rho,\theta)\nabla\theta)}{\rho}-
\frac{\chi(\bar{\rho},\bar{\theta})}{\bar{\rho}}\D\theta\,\|_{L^{1}(\widetilde{B}^{\N-1,\N})}\leq\|K(q){\rm
div}(K_{1}(q,{\cal
T})\n\theta)\|_{L^{1}(\widetilde{B}^{\N-1,\N})}\\[2mm]
&\hspace{4,3cm}+\|{\rm
div}(K_{1}(q,{\cal T})\n\theta)\|_{L^{1}(\widetilde{B}^{\N-1,\N})}+\|K(q)\D\theta\|_{L^{1}(\widetilde{B}^{\N-1,\N})}\\
\end{aligned}
$$
\\
and in using the propositions \ref{fcomposition} and \ref{fproduit} we
get:
$$
\begin{aligned}
\|{\rm div}( K_{1}(q,{\cal
T})\n\theta)\|_{L^{1}(\widetilde{B}^{\N-1,\N})}\lesssim&\,(
\|q\|_{L^{\infty}(B^{\N})}+\|{\cal T}\|_{L^{\infty}(B^{\N})})\|{\cal
T}\|_{L^{1}(\widetilde{B}^{\N+1,\N+2})}\\[2mm]
&+(\|q\|_{L^{2}(\widetilde{B}^{\N,\N+1})}+\|{\cal
T}\|_{L^{2}(\widetilde{B}^{\N,\N+1})})\|{\cal
T}\|_{L^{2}(B^{\N+1})}.\\
\end{aligned}
$$
Next we have:
$$
\begin{aligned}
\|(\frac{TP_{1}(\rho)}{\rho}-\frac{\bar{T}P_{1}(\bar{\rho})}{\bar{\rho}}){\rm
div}&\,u\|_{L^{1}(\widetilde{B}^{\N-1,\N})}\lesssim\|L_{1}(q){\rm
div}
u\|_{L^{1}(\widetilde{B}^{\N-1,\N})}\\[2mm]
&+\|L_{1}(q)L_{2}({\cal T}){\rm div} u\|_{L^{1}(\widetilde{B}^{\N-1,\N})}+\|L_{2}({\cal T}){\rm div} u\|
_{L^{1}(\widetilde{B}^{\N-1,\N})},\\
\end{aligned}
$$
where:
$$
\begin{aligned}
&\|L_{1}(q){\rm div}
u\|_{L^{1}(\widetilde{B}^{\N-1,\N})}\lesssim\|q\|_{L^{\infty}(\widetilde{B}^{\N-1,\N})}
\|u\|_{L^{1}(B^{\N+1})}\\[2mm]
&\|L_{2}({\cal T}){\rm div}
u\|_{L^{1}(\widetilde{B}^{\N-1,\N})}\lesssim\|{\cal
T}\|_{L^{\infty}(\widetilde{B}^{\N-1,\N})}
\|u\|_{L^{1}(B^{\N+1})}\\
\end{aligned}
$$
so we get:
$$
\begin{aligned}
\|(\frac{TP_{1}(\rho)}{\rho}-\frac{\bar{T}P_{1}(\bar{\rho})}{\bar{\rho}}){\rm
div}\,u&\|_{L^{1}(\widetilde{B}^{\N-1,\N})}\lesssim\;\|{\cal
T}\|_{L^{\infty}(\widetilde{B}^{\N-1,\N})}
\|u\|_{L^{1}(B^{\N+1})}\|q\|_{L^{\infty}(\widetilde{B}^{\N-1,\N})}\\[2mm]
&\hspace{1,7cm}+\big(\|q\|_{L^{\infty}(\widetilde{B}^{\N-1,\N})}+\|{\cal
T}\|_{L^{\infty}(\widetilde{B}^{\N-1,\N})}\big)
\|u\|_{L^{1}(B^{\N+1})}\;.\\
\end{aligned}
$$
To end with, we have the last two terms:
$$
\begin{aligned}
&\|u^{*}.\n\theta\|_{ L^{1}(\widetilde{B}^{\N-1,\N})}\leq\|{\cal
T}\|_{
L^{1}(\widetilde{B}^{\N+1,\N+2})}\|u\|_{L^{\infty}(\widetilde{B}^{\N-1,\N})}\,,\\[2mm]
&\|\frac{D:\nabla u}{\rho}\|_{
L^{1}(\widetilde{B}^{\N-1,\N})}\leq\|K(q)\nabla u:\nabla u\|_{
L^{1}(\widetilde{B}^{\N-1,\N})}+\|\nabla u:\nabla u\|_{
L^{1}(\widetilde{B}^{\N-1,\N})}\\
&\hspace{3,6cm}\lesssim(1+\|q\|_{L^{\infty}(\widetilde{B}^{\N,\N+1})})
\|u\|_{ L^{2}(\widetilde{B}^{\N,\N+1})}\|u\|_{ L^{2}(B^{\N+1})}.\\
\end{aligned}
$$
\\
Finally we have in using (\ref{fa5}), (\ref{fa6}) and the previous bounds:
\begin{equation}
\|\psi(q,u^{'},{\cal T}^{'})\|_{E^{N/2}}\leq C((C+1)\eta+R)^{2}
\label{a7}
\end{equation}
Let $c$ such that $\|\cdot \|_{B^{N/2}}\leq c$ implies that:
$\|\cdot
\|_{ L^{\infty}}\leq1/3$ then we choose R and $\alpha_{0}$ such that:
$$R\leq \inf((3C)^{-1},c,1),\;\alpha_{0}\leq \inf\frac{(R,c)}{C+1}\;.$$
So $(\cal{H})$ is verified and we have then:
$$\psi(B(0,R))\subset B(0,R)\;.$$
Next one can proceed as in the proof of the theorem \ref{ftheo1}, we
have to show
the contraction of the application $\psi$ to use the theorem of the fixed point.\\
The uniqueness of the solution in the space $F^{\N}$ follows the
same lines as in theorem \ref{ftheo1}. The details are left to the
reader.
\section{Local theory for large data\label{fsection5}}
In this part we are interested in results of
existence in finite time for general initial data with density
bounded away from zero. We focus on the case where the coefficients
depend only on the density with linear specific energy, and next
we will treat the general case. As a first step, we shall study the
linear part of the system $(NHV)$ about non constant reference
density and temperature, that is:
$$
\begin{cases}
\p_{t}q+{\rm div}u=F,\\
\p_{t}u-{\rm div}(a\n u)-\n(b{\rm
div}u)-\n(c\D q)=G,\\
\p_{t}{\cal T}-{\rm div}(d\n{\cal T})=H,\\
\end{cases}
\leqno{(N)}
$$
\subsection{Study of the linearized equation}
We want to prove a priori estimates for system $(N)$
with the following hypotheses on $a,b,c,d$:
$$
\begin{aligned}
&0<c_{1}\leq a<M_{1}<\infty,\;0<c_{2}\leq
a+b<M_{2}<\infty,\;0<c_{3}\leq c<M_{3}<\infty,\\
&0<c_{4}\leq d<M_{4}<\infty.\\
\end{aligned}
$$
We remark that the last equation is just a heat equation with
variable coefficients so that one can apply the following
proposition proved in \cite{fDL}.
\begin{proposition}
\label{fchaleur3}
Let ${\cal T}$ solution
of the heat equation:
$$\p_{t}{\cal T}-{\rm div}(d\n{\cal T})=H,$$
we have so for all
index $\tau$ such that $-\N-1<\tau\leq \N-1$ the following estimate for all $\alpha\in[1,+\infty]$:
$$\|{\cal T}\|_{\widetilde{L}^{\alpha}_{T}(B^{\tau+\frac{2}{\alpha}})}\leq\|{\cal T}_{0}\|_{B^{\tau}}+
\|H\|_{\widetilde{L}^{1}_{T}(B^{\tau})} +\|\n
d\|_{\widetilde{L}^{\infty}_{T}(B^{\N-1})}\|\n {\cal
T}\|_{\widetilde{L}^{1}_{T}(B^{\tau+1})}.$$
\end{proposition}
We are now interested by the first two equations of the system
$(N)$.
$$
\begin{cases}
\p_{t}q+{\rm div}u=F\\
\p_{t}u-{\rm div}(a\n u)-\n(b{\rm
div}u)-\n(c\D q)=G\\
\end{cases}
\leqno{(N^{'})}
$$
where we keep the same hypothesis on $a$, $b$ and $c$. We have then
the following estimate of the solution in the
spaces of Chemin-Lerner:\\
\begin{proposition}
\label{flinear3}  Let $1\leq r_{1}\leq r\leq+\infty$ , $0\leq
s\leq1$, $(q_{0},u_{0})\in B^{\N+s}\times (B^{\N-1+s})^{N}$, and
$(F,G)\in\widetilde{L}^{r_{1}}_{T}(B^{\N-2+s+2/r_{1}})\times
(\widetilde{L}^{r_{1}}_{T}(B^{\N-3+s+2/r_{1}}))^{N}$.\\
\\
Suppose that $\n a$ , $\n b$ , $\n c$ belong to
$\widetilde{L}^{2}_{T}(B^{\N})$ and that $\p_{t}c\in L^{1}_{T}(L^{\infty})$.\\
\\
Let $(q,u)\in
(\widetilde{L}^{r}_{T}(B^{\N+s+2/r})\cap\widetilde{L}^{2}_{T}(B^{\N+s+1}))\times((\widetilde{L}^{r}_{T}(B^{\N+s-1+2/r}))
^{N}\cap(\widetilde{L}^{2}_{T}(B^{\N+s})^{N})$ be a solution of the
system $(N^{'})$.
\\
\\
Then there exists a constant $C$ depending only on $r,\, r_{1},\,
\bar{\lambda},\, \bar{\mu},\, \bar{\kappa},\, c_{1},\, c_{2},\,
M_{1}$ and
$M_{2}$ such that:
$$
\begin{aligned}
&\|(\n q,u)\|_{\widetilde{L}^{r}_{T}(B^{\N-1+s+2/r})}(1-C\|\n
c\|_{L^{2}_{T}(L^{\infty})})\leq\|(\n q_{0},u_{0})\|_{B^{\N}}+\|(\n
F,G)\|_{\widetilde{L}^{r_{1}}_{T}(B^{\N-3+s+2/r_{1}})}\\[3mm]
&\hspace{1,5cm}+ \|\n
q\|_{\widetilde{L}^{\infty}_{T}(B^{\N-1+s})}\|\p_{t}
c\|_{L^{1}_{T}(L^{\infty})} +\|(\n
q,u)\|_{\widetilde{L}^{2}_{T}(B^{\N+s})}(\|\n
a\|_{\widetilde{L}^{2}_{T}(B^{\N})}+\|\n
b\|_{\widetilde{L}^{2}_{T}(B^{\N})}\\[3mm]
&\hspace{12,5cm}+\|\n
c\|_{\widetilde{L}^{2}_{T}(B^{\N})})\;.\\
\end{aligned}
$$
\end{proposition}
{\bf Proof:} \\
\\
Like previously we are going to show estimates on $q_{l}$ and
$u_{l}$. So we apply to the system the operator $\D_{l}$ , and we
have then:
\begin{eqnarray}
&&\p_{t}q_{l}+{\rm div}u_{l}=F_{l}\label{fl1}\\
&&\p_{t}u_{l}-{\rm div}(a\n u_{l})-\n(b\,{\rm div}u_{l})-\n(c\D
q_{l})=G_{l}+R_{l}\label{fl2}
\end{eqnarray}
where we denote:
$$R_{l}={\rm div}([a,\D_{l}]\n
u)-\n([b,\D_{l}]{\rm div} u_{l})-\n([c,\D_{l}]\D q).$$
Performing integrations by parts and usinf (\ref{fl1}) we have:
$$
\begin{aligned}
-\int_{\R^{N}}u_{l}\n(c\D
q_{l})dx=\frac{1}{2}\frac{d}{dt}\int_{\R^{N}}c|\n
q_{l}|^{2}dx-\int_{\R^{N}}\big({\rm div}u_{l}\,(\n q_{l}.\n
c)&+\frac{|\n q|^{2}}{2}\p_{t}c\\
&+c.\n
q_{l}.\n F_{l}\big)\,dx.\\
\end{aligned}
$$
Next, we take the inner product of (\ref{fl2}) with $u_{l}$
and we use the previous equality, we have then:
$$
\begin{aligned}
\frac{1}{2}\frac{d}{dt}\big(\|u_{l}\|_{L^{2}}^{2}+\int_{\R^{N}}c|\n
q_{l}|^{2}dx\big)&+\g(a|\n u_{l}|^{2}+b|{\rm
div}u_{l}|^{2})dx=\g((G_{l}+R_{l}).u_{l}\,dx\\
&+\g\big(({\rm div}u_{l}(\n c.\n q_{l})+\frac{|\n
q_{l}^{2}|}{2}\p_{t}c+c\n q_{l}.\n F_{l})\big)\,dx\,.\\
\label{f3}
\end{aligned}
$$
In  order to recover some terms in $\D q_{l}$ we take the
inner product of the gradient of (\ref{fl1}) with $u_{l}$, the inner
product
scalar of (\ref{fl2}) with $\n q_{l}$ and we sum, we obtain then:
\begin{equation}
\begin{aligned}
\frac{d}{dt}\g\n q_{l}.u_{l}dx+\g c(\D q_{l})^{2}dx
=&\g((G_{l}+R_{l}).\n q_{l}+|{\rm div}u_{l}|^{2}+u_{l}.\n F_{l}\\
&\hspace{2cm}-a\n u_{l}:\n ^{2}q_{l}-b\D q_{l}{\rm div}u
_{l})dx.\\
\label{f4}
\end{aligned}
\end{equation}
Let $\alpha>0$ small enough. We define:
\begin{equation}
k_{l}^{2}=\|u_{l}\|_{L^{2}}^{2}+\g(\bar{\kappa}
c|\n\q|^{2}+2\alpha\n\q.\ui)dx\;.\label{f5.39}
\end{equation}
In using the previous inequality and the fact that
$a_{1}b_{1}\leq\frac{1}{2}a_{1}^{2}+\frac{1}{2}b_{1}^{2}$, we have
in summing:
\begin{equation}
\begin{aligned}
&\frac{1}{2}\frac{d}{dt}k_{l}^{2}+\frac{1}{2}\g(a|\n
u_{l}|^{2}+\alpha b|\D\q|^{2})dx
\lesssim(\|G_{l}\|_{L^{2}}+\|R_{l}\|_{L^{2}})\\[2mm]
&\times(\alpha\|\n\q\|_{L^{2}}+\|\ui\|_{L^{2}})+\|\n F_{l}\|_{L^{2}}(\alpha\|\ui\|_{L^{2}}+\|c\n\q\|_{L^{2}})
+\frac{1}{2}\|\p_{t}c\|_{L^{\infty}}\|\n\q\|_{L^{2}}^{2}\\[2mm]
&\hspace{8,5cm}+\|\n
c\|_{L^{\infty}}\|\n\q\|_{L^{2}}\|\n\ui\|_{L^{2}} \;.\\
\label{fl5}
\end{aligned}
\end{equation}
For small enough $\alpha$, we have according (\ref{f5.39}):
\begin{equation}
\frac{1}{2}k_{l}^{2}\leq\|u_{l}\|^{2}+\g\bar{\kappa}
c|\n\q|^{2}dx\leq\frac{3}{2}k_{l}^{2}\;. \label{fl6}
\end{equation}
Hence according to (\ref{fl5}) and (\ref{fl6}):
$$
\begin{aligned}
\frac{1}{2}\frac{d}{dt}k_{l}^{2}+K2^{2l}k_{l}^{2} \leq&\,\,
k_{l}\,(\|G_{l}\|_{L^{2}}+\|R_{l}\|_{L^{2}}+\|\n
F_{l}\|_{L^{2}})\|\p_{t}c\|_{L^{\infty}}\|\n\q\|_{L^{2}}+2^{2l}k_{l}^{2}\|\n
c\|_{L^{2}}\;.\\
\end{aligned}
$$
By integrating with respect to the time, we obtain:
$$
\begin{aligned}
k_{l}(t)\leq &\,e^{-K2^{2l}t}k_{l}(0)+C\int_{0}^{t}
e^{-K2^{2l}(t-\tau)}(\|\p_{t}c\|_{L^{\infty}}\|\n\q(\tau)\|_{L^{2}}+\|\n
F_{l}(\tau)\|_{L^{2}}+\|G_{l}(\tau)\|_{L^{2}}\\[2,5mm]
&\hspace{7cm}+\|R_{l}(\tau)\|_{L^{2}}+2^{l}k_{l}(\tau)\|\n
c(\tau)\|_{L^{2}})d\tau\;.\\
\end{aligned}
$$
After convolution inequalities imply that:
\begin{equation}
\begin{aligned}
\|k_{l}\|_{L^{r}([0,T])}\leq&\,(2^{-\frac{2l}{r}}k_{l}(0)+(2^{-2l(1+1/r-1/r_{1})}\|(\n
F_{l},G_{l})\|_{L^{r_{1}}_{T}(L^{2})}+2^{-\frac{2l}{r}}\|R_{l}\|_{L^{1}_{T}(L^{2})}\\[2,5mm]
&+2^{-\frac{2l}{r}}\|\n\q\|_{L^{\infty}_{T}(L^{2})}\|\p_{t}c\|_{L^{1}_{T}(L^{\infty})}+\|\n
c\|_{L^{2}_{T}(L^{\infty})}\|k_{l}\|_{L^{r}([0,T])}\;.\\
\label{f8}
\end{aligned}
\end{equation}
Moreover we have:
$$C^{-1}\,k_{l}\leq\|\n\q\|_{L^{2}}+\|\ui\|_{L^{2}}\leq C\,k_{l}.$$
Finally multiplying by $2^{(\N-1+s+\frac{2}{r})l}$ and using (\ref{fl6}), we end up with:
$$
\begin{aligned}
\|(\n q,u)&\|_{L^{r}_{T}(B^{\N-1+s+2/r})}(1-C\|\n
c\|_{L^{2}(L^{\infty})})\leq\,\|(\n
F,G)\|_{\widetilde{L}^{r_{1}}_{T}(B^{\N-3+s+2/r_{1}})}\\[3mm]
&\hspace{0,2cm}\|(\n q_{0},u_{0})\|_{B^{\N-1+s}}+\|\n
q\|_{\widetilde{L}^{\infty}_{T}(B^{\N-1+s})}\|\p_{t}c\|_{L^{1}_{T}(L^{\infty})}+\sum_{l\in\mathbb{Z}}2^{l(\N+s-1)}\|R_{l}\|_{L^{1}_{T}(L^{2})}\;.\\
\label{f27}
\end{aligned}
$$
Finally, applying lemma \ref{fcommutateur} on the appendix to bound the remainder term completes the proof
$$
\begin{aligned}
&\sum_{l\in\mathbb{Z}}2^{l(\N+s-1)}\|R_{l}\|_{L^{1}_{T}(L^{2})}\leq
C\|a\|_{\widetilde{L}^{2}_{T}(B^{\N+1})}\|u\|_{\widetilde{L}^{2}_{T}(B^{\N})}+
C\|b\|_{\widetilde{L}^{2}_{T}(B^{\N+1})}\|u\|_{\widetilde{L}^{2}_{T}(B^{\N})}\\[2mm]
&\hspace{9cm}+C\|c\|_{\widetilde{L}^{2}_{T}(B^{\N})}\|q\|_{\widetilde{L}^{2}_{T}(B^{\N+1})}\;.\\
\end{aligned}
$$
\hfill {$\Box$}
\subsection{Local existence Theorem for temperature independent coefficients}
We recall the space we will work with:
$$
\begin{aligned}
F_{T}=&[\widetilde{C}_{T}(B^{\N})\cap
L^{1}_{T}(B^{\N+2})]\times[\widetilde{C}_{T}(B^{\N-1})^{N}\cap
L^{1}_{T}(B^{\N+1})^{N}]\times[\widetilde{C}_{T}(B^{\N-2})\cap
L^{1}_{T}(B^{\N})]\\
\end{aligned}
$$
endowed with the following norm:
$$
\begin{aligned}
\|(q,u,{\cal
T})\|_{F_{T}}=&\|q\|_{L^{1}_{T}(B^{\N+2})}+\|q\|_{\widetilde{L}^{\infty}_{T}(B^{\N})}+\|u\|_{L^{1}_{T}(B^{\N+1})}
+\|u\|_{\widetilde{L}^{\infty}_{T}(B^{\N-1})}\\[2mm]
&+\|{\cal T}\|_{L^{1}_{T}(B^{\N})}+\|{\cal T}\|_{\widetilde{L}^{\infty}_{T}(B^{\N-2})}.\\
\end{aligned}
$$
We will now prove the  local existence of a solution for general
initial data with a linear specific intern energy and coefficients
independent of the temperature. The functional space we shall work with is larger
than previously, the reason why is that the low frequencies
don't play an important role as far as one is interested in {\it local} results.\\
\\
In what follows, $N\geq3$ is assumed.
\\
\\
{\bf Proof of the theorem \ref{ftheo4}:}\\
\\
Let:
$$q^{n}=q^{0}+\bar{q}^{n},\;\rho^{n}=\bar{\rho}(1+q^{n}),\;u^{n}=u^{0}+\bar{u}^{n},\;
{\cal T}^{n}={\cal T}^{0}+\bar{{\cal
T}}^{n}\;\;\mbox{and}\;\;\theta^{n}=\bar{\theta}+{\cal T}^{n}$$
where $(q^{0},u^{0},{\cal T}^{0})$ stands for the  solution of:\\
$$
\begin{cases}
\p_{t}q^{0}-\D q^{0}=0,\\
\p_{t}u^{0}-\D u^{0}=0,\\
\p_{t}{\cal T}^{0}-\D {\cal T}^{0}=0,\\
\end{cases}
$$
supplemented with initial data:
$$q^{0}(0)=q_{0}\;,\;u^{0}(0)=u_{0}\;,\;{\cal T}^{0}(0)={\cal T}_{0}.$$
Let $(\bar{q}_{n},\bar{u}_{n},\bar{{\cal T}}_{n})$ be the solution
of
the following system:
$$
\begin{cases}
\begin{aligned}
&\p_{t}\bar{q}_{n+1}+{\rm div}(\bar{u}_{n+1})=F_{n},\\[3mm]
& \p_{t}\bar{u}_{n+1}-{\rm
div}\biggl(\frac{\mu(\rho^{n})}{\rho^{n}}\n
\bar{u}^{n+1}\biggl)-\n\biggl(\frac{\zeta(\rho^{n})}{\rho^{n}} {\rm
div}(\bar{u}^{n+1})\biggl)-\n(K(\rho^{n})\D\bar{q}^{n+1})=G_{n},\\[3mm]
&\p_{t}\bar{{\cal T}}_{n+1}-{\rm
div}\biggl(\frac{\chi(\rho^{n})}{1+q^{n}}\bar{{\cal
T}}_{n+1}\biggl)=H_{n},\\[3mm]
&(\bar{q}_{n+1},\bar{u}_{n+1},\bar{{\cal T}}_{n+1})_{t=0}=(0,0,0),\\
\end{aligned}
\end{cases}
\leqno{(N_{1})}
$$
where:
$$
\begin{aligned}
F_{n}=&-{\rm div}(q^{n}u^{n})-\D q^{0}-{\rm div}(u^{0}),\\[2mm]
G_{n}=&-(u^{n})^{*}.\n
u^{n}+\n(\frac{K^{'}_{\rho^{n}}}{2}|\n\rho^{n}|^{2}
)-\n(\frac{\mu(\rho^{n})}{\rho^{n}}){\rm
div}u^{n}+\n(\frac{\zeta(\rho^{n})}{\rho^{n}}){\rm
div}u^{n}\\[2mm]
&+\frac{\lambda^{'}(\rho^{n})\n\rho^{n}{\rm
div}u^{n}}{1+q^{n}}+\frac{(du^{n}+\n
u^{n})\,\mu^{'}(\rho^{n})\n\rho^{n}}{1+q^{n}}+\biggl[\frac{P_{1}(\rho^{n})}{\rho^{n}\psi^{'}(\w^{n})}\biggl]
\nabla\theta^{n}\\[2mm]
\end{aligned}
$$
$$
\begin{aligned}
&+\frac{[P_{0}^{'}(\rho^{n})+\w^{n}P_{1}^{'}(\rho^{n})]\n
q^{n}}{1+q^{n}}-\D u^{0}+{\rm
div}\biggl(\frac{\mu(\rho^{n})}{1+q^{n}}\n
u^{0}\biggl)\\[2mm]
&+\n\biggl(\frac{\mu(\rho^{n})+\lambda(\rho^{n})}{1+q^{n}}{\rm
div}(u^{0})\biggl)+\n(K(\rho^{n})\D q^{0}),\\
\end{aligned}
$$
$$
\begin{aligned}
H_{n}=&\,\n(\frac{1}{1+q^{n}}).\n\theta^{n}\chi(\rho^{n})-\frac{\w^{n}P_{1}(\rho^{n})}{\rho^{n}}{\rm
div}u^{n} -(u^{n})^{*}.\n\theta^{n}+\frac{D_{n}:\nabla
u^{n}}{\rho^{n}}\\[2mm]
&\,-\D\theta_{0}+{\rm
div}\big(\frac{\chi(\rho^{n})}{1+q^{n}}\n\theta^{0}\big).\\
\end{aligned}
$$
\subsubsection*{1) First Step , Uniform Bound}
 Let $\e$ be a small
parameter and  choose $T$ small enough  so that
in using the estimate of the heat equation stated in proposition \ref{fchaleur} we have:
$$
\begin{aligned}
&\|{\cal T}^{0}\|_{L^{1}_{T}(B^{\N})}+\|u^{0}\|_{L^{1}_{T}(B^{\N+1})}+\|q^{0}\|_{L^{1}_{T}(B^{\N+2})}\leq\e,\\
&\|{\cal
T}^{0}\|_{\widetilde{L}^{\infty}_{T}(B^{\N-2})}+\|u^{0}\|_{\widetilde{L}^{\infty}_{T}(B^{\N-1})}+\|q^{0}\|_{\widetilde{L}^{\infty}_{T}(B^{\N})}\leq
A_{0}.\\
\end{aligned}
\leqno{({\cal{H}}_{\e})}
$$
We are going to show by induction that:
$$\|(\bar{q}^{n},\bar{u}^{n},\bar{{\cal T}}^{n})\|_{F_{T}}\leq\e.\leqno{({\cal{P}}_{n})}$$
As $(\bar{q}_{0},\bar{u}_{0},\bar{{\cal T}}_{0})=(0,0,0)$ the result
is true for $n=0$. We suppose now $({\cal{P}}_{n})$ true and we are
going to show $({\cal{P}}_{n+1})$.
\\
To begin with we are going to show that $1+q^{n}$ is  positive.
Using the fact that $B^{\N}\h L^{\infty}$ and that we take $\e$
small enough, we have for $t\in[0,T]$:
$$
\begin{aligned}
\|q^{n}-q_{0}\|_{L^{\infty}((0,T)\times\R^{N})}&\lesssim\|{\rm
div}\bar{u}^{n}\|_{L^{1}_{T}(B^{\N})}
+\|{\rm div}(q^{n-1}u^{n-1})\|_{L^{1}_{T}(B^{\N})}+\|{\rm div}u^{0}\|_{L^{1}_{T}(B^{\N})},\\[2mm]
&\lesssim2\e+\|q^{n-1}u^{n-1}\|_{L^{1}_{T}(B^{\N+1})},\\
\end{aligned}
$$
and:
$$
\|q^{n-1}u^{n-1}\|_{L^{1}_{T}(B^{\N+1})}\leq\|q^{n-1}\|_{L^{\infty}_{T}(B^{\N})}\|u^{n-1}\|_{L^{1}_{T}(B^{\N+1})}
+\|q^{n-1}\|_{L^{2}_{T}(B^{\N+1})}\|u^{n-1}\|_{L^{2}_{T}(B^{\N})}.$$
Hence:
$$\|q^{n}-q_{0}\|_{L^{\infty}((0,T)\times\R^{N})}\leq C_{1}(2\e+(A_{0}+\e)\e).$$
Finally we thus have:
$$
\begin{aligned}
\|1+q_{0}\|_{L^{\infty}((0,T)\times\R^{N})}-\|q^{n}-q_{0}\|_{L^{\infty}((0,T)\times\R^{N})}\leq
1+q^{n}\leq\|1&+q_{0}\|_{L^{\infty}((0,T)\times\R^{N})}\\
&+\|q^{n}-q_{0}\|_{L^{\infty}((0,T)\times\R^{N})},
\end{aligned}$$
whence if $\e$ is small enough:
$$\frac{c}{2\bar{\rho}}\leq1+q^{n}\leq1+\frac{\|\rho_{0}\|_{L^{\infty}}}{\bar{\rho}}.$$
In order to bound $(\bar{q}^{n},\bar{u}^{n},\bar{{\cal T}}^{n})$ in $F_{T}$, we shall use proposition \ref{flinear3}. For that
we must check that the different hypotheses of this proposition
adapted to our system $(N_{1})$ are satisfied, so we study the
following terms:
$$a^{n}=\frac{\mu(\rho^{n})}{1+q^{n}}\;,\;b^{n}=\frac{\zeta(\rho^{n})}{1+q^{n}}\;,\;c^{n}=K(\rho_{n})\;,\;
d^{n}=\frac{\chi(\rho^{n})}{1+q^{n}}.$$ In using $({\cal{P}}_{n})$
and by continuity of $\mu$ and the fact that $\mu$ is positive on
$\big[\bar{\rho}\big(1+\min(q_{0})\big)-\alpha,\bar{\rho}\big(1+\max(q_{0})\big)+\alpha\big]$,
we have:
$$0<c_{1}\leq a^{n}=\frac{\mu(\rho^{n})}{1+q^{n}}\leq M_{1}\;.$$
We proceed similarly for the others terms.
\\
Next, notice that:
$$
\begin{aligned}
\|\n
a^{n}\|_{\widetilde{L}_{T}^{2}(B^{\N})}&\leq\|\frac{\mu(\rho^{n})}{1+q^{n}}-\mu(\bar{\rho})\|_{\widetilde{L}^{2}_{T}(B^{\N+1})}
\leq C\|q^{n}\|_{\widetilde{L}^{2}_{T}(B^{\N+1})}.\\[3mm]
\|\n
b^{n}\|_{\widetilde{L}_{T}^{2}(B^{\N})}&\leq\|\frac{\zeta(\rho^{n})}{1+q^{n}}-\zeta(\bar{\rho})
\|_{\widetilde{L}^{2}_{T}(B^{\N+1})}\leq C\|q^{n}\|_{\widetilde{L}^{2}_{T}(B^{\N+1})}\\[3mm]
\|\n c^{n}\|_{\widetilde{L}_{T}^{2}(B^{\N})}&\leq
C\|q^{n}\|_{\widetilde{L}^{2}_{T}(B^{\N+1})}\;.\\[2mm]
\end{aligned}
$$
To end on our hypotheses we have to control $\p_{t}c^{n}$ in norm
$\|\cdot\|_{L^{1}_{T}(L^{\infty})}$. As $B^{\N}\h L^{\infty}$, it
actually suffices to bound $\|\p_{t}c^{n}\|_{L^{1}_{T}(B^{\N})}$. We
have:
$$\p_{t}c^{n}=K^{'}(\rho^{n})\p_{t}q^{n}=K^{'}(\rho^{n})({\rm
div}(q^{n-1}u^{n-1})-{\rm div}(u^{n}))\;.$$
And we have in using the propositions \ref{fproduit} and
\ref{fcomposition}:
$$
\begin{aligned}
&\|K^{'}(\rho^{n})({\rm div}(q^{n-1}u^{n-1})-{\rm
div}(u^{n}))\|_{L^{1}_{T}(B^{\N})}\leq\|K^{'}(\rho^{n}){\rm
div}(q^{n-1}u^{n-1})\|_{L^{1}_{T}(B^{\N})}\\[2mm]
&\hspace{9,7cm}+\|K^{'}(\rho^{n}){\rm
div}(u^{n})\|_{L^{1}_{T}(B^{\N})},\\[2mm]
&\lesssim
(1+\|q^{n}\|_{L^{\infty}_{T}(B^{\N})})(\|u^{n}\|_{L^{1}_{T}(B^{\N+1})}+\|q^{n-1}u^{n-1}\|_{L^{1}_{T}(B^{\N+1})})\\[3mm]
&\lesssim(1+\|q^{n}\|_{L^{\infty}_{T}(B^{\N})})(\|u^{n}\|_{L^{1}_{T}(B^{\N+1})}+\|q^{n-1}\|_{L^{\infty}_{T}(B^{\N})}
\|u^{n-1}\|_{L^{1}_{T}(B^{\N+1})}\\[1,8mm]
&\hspace{8,5cm}+\|q^{n-1}\|_{L^{2}_{T}(B^{\N+1})}\|u^{n-1}\|_{L^{2}_{T}(B^{\N})})\;.\\
\end{aligned}
$$
We now use proposition \ref{fchaleur3} to get the bound on
$\bar{{\cal T}}^{n}$,
so we obtain in taking $\tau=\N-2$:
\begin{equation}
\begin{aligned}
&\|{\cal T}^{n}\|_{L^{1}_{T}(B^{\N})\cap
L^{\infty}_{T}(B^{\N-2})}\lesssim\big(\|H_{n}\|_{L^{1}_{T}(B^{\N-2})}
+\|\n(\frac{\chi(\rho^{n})}{\rho^{n}})\|_{L^{\infty}_{T}(B^{\N-1})}\\
&\hspace{10cm}\times\|{\cal
T}^{n}\|_{L^{1}_{T}(B^{\N})}\big).
\end{aligned}
\label{fi2}
\end{equation}
So we need to bound $d^{n}$ in $L^{\infty}_{T}(B^{\N})$: $$\|\n
d^{n}\|_{L_{T}^{\infty}(B^{\N-1})}\leq
C\|q^{n}\|_{L_{T}^{\infty}(B^{\N})}.$$
\\
Now we show by induction $({\cal{P}}_{n+1})$. Finally, applying
the estimates of propositions \ref{flinear3} and \ref{fchaleur3}, we conclude that:
\begin{equation}
\begin{aligned}
&\|(\bar{q}^{n+1},\bar{u}^{n+1},\bar{{\cal
T}}^{n+1})\|_{F_{T}}\big(1-C(\|a^{n}\|_{L^{2}_{T}(B^{\N+1})}
+\|b^{n}\|_{L^{2}_{T}(B^{\N+1})}+\|c^{n}\|_{L^{2}_{T}(B^{\N+1})}\\[3mm]
&+\|d^{n}\|_{L^{\infty}_{T}(B^{\N})}+\|\p_{t}
c^{n}\|_{L^{1}_{T}(B^{\N})})\big)\leq\,\|(\n
F_{n},G_{n})\|_{L^{1}_{T}(B^{\N-1})}+\|H_{n}\|_{L^{1}_{T}(B^{\N-2})}\;.\\
\end{aligned}
\label{fi1}
\end{equation}
\\
Bounding the right-hand side may be
done by applying
propositions  \ref{fproduit} and \ref{fcomposition}. For instance, we have:
$$\|F_{n}\|_{L^{1}_{T}(B^{N/2})}\leq\|{\rm
div}(q^{n}u^{n})\|_{L^{1}_{T}(B^{N/2})}+\|{\rm
div}u^{0}\|_{L^{1}_{T}(B^{N/2})}+\|\D
q^{0}\|_{L^{1}_{T}(B^{N/2})}.$$
Since:
$$\|u^{n}q^{n}\|_{L^{1}_{T}(B^{N/2+1})}\lesssim\|q^{n}\|_{L^{\infty}_{T}(B^{N/2})}\|u^{n}\|_{L^{1}_{T}(B^{N/2+1})}
+\|q^{n}\|_{L^{2}_{T}(B^{N/2+1})}\|u^{n}\|_{L^{2}_{T}(L^{\infty})},$$
we can conclude that:
$$\|F_{n}\|_{L^{1}_{T}(B^{N/2})}\leq C(A_{0}+\e+\sqrt{\e})^{2}.$$
\\
Next we want to control the different terms of $G_{n}$. According to
propositions \ref{fproduit} and \ref{flinear3}, we have:
$$
\begin{aligned}
&\|(u^{n})^{*}.\n
u^{n}\|_{L^{1}_{T}(B^{\N-1})}\lesssim\|u_{n}\|_{L^{\infty}_{T}(B^{\N-1})}\|u_{n}\|_{L^{1}_{T}(B^{\N+1})}\\[2mm]
&\|\n\biggl(\frac{K^{'}_{\rho^{n}}}{2}|\n\rho_{n}|^{2}
\biggl)\|_{L^{1}_{T}(B^{\N-1})}\lesssim\|L(q^{n})|\n\rho_{n}|^{2}
\|_{L^{1}_{T}(B^{\N})}+\||\n\rho_{n}|^{2}
\|_{L^{1}_{T}(B^{\N})}\\[2mm]
&\hspace{4,7cm}\lesssim(1+\|q^{n}\|^{2}_{L^{\infty}_{T}(B^{\N})})\|q^{n}\|^{2}_{L^{2}_{T}(B^{\N+1})}.\\[2mm]
\end{aligned}
$$
After we have:
$$
\begin{aligned}
\|\n(\frac{\mu(\rho^{n})}{\rho^{n}}){\rm
div}u_{n}\|_{L^{1}_{T}(B^{\N-1})}&\lesssim\|{\rm
div}u_{n}\|_{L^{1}_{T}(B^{\N})}\|\n(\frac{\mu(\rho^{n})}{\rho^{n}}
-\frac{\mu(\bar{\rho})}{\bar{\rho}})\|_{L^{\infty}_{T}(B^{\N-1})}\\[2mm]
&\leq C\|u^{n}\|_{L^{1}_{T}(B^{\N+1})}
\|q^{n}\|_{L^{\infty}_{T}(B^{\N})}\;.\\[3mm]
\|\n(\frac{\zeta(\rho^{n})}{\rho^{n}}){\rm
div}u_{n}\|_{L^{1}_{T}(B^{\N-1})}&\leq
C\|u^{n}\|_{L^{1}_{T}(B^{\N+1})}\|q^{n}\|_{L^{\infty}_{T}(B^{\N})}\;.\\
\end{aligned}
$$
After we study the term coming from ${\rm div}(D)$:
$$
\begin{aligned}
\|\frac{\lambda^{'}(\rho^{n})\n\rho^{n}\rm
div(u^{n})}{1+q^{n}}\|_{L^{1}_{T}(B^{\N-1})}
&\lesssim\|L(q^{n})\n\rho^{n}\,{\rm
div}(u^{n})\|_{L^{1}_{T}(B^{\N-1})}+\|\n\rho^{n}\,{\rm
div}(u^{n})\|_{L^{1}_{T}(B^{\N-1})}\\[2mm]
&\lesssim(1+\|q^{n}\|_{L^{\infty}_{T}(B^{\N})})\|u^{n}\|_{L^{1}_{T}(B^{\N+1})}
\|q^{n}\|_{L^{\infty}_{T}(B^{\N})}\;.\\
\end{aligned}
$$
We proceed similarly for the following term: $$\frac{(du^{n}+\n
u^{n})\mu^{'}(\rho^{n})\n\rho^{n}}{1+q^{n}}.$$
Next we study the
last terms:
$$
\begin{aligned}
&\|\frac{[P_{0}^{'}(\rho^{n})+{\cal T}^{n}P_{1}^{'}(\rho^{n})]\n
q^{n}}{1+q^{n}}\|_{L^{1}_{T}(B^{\N-1})}\lesssim \|K(q^{n})\n
q^{n}\|_{L^{1}_{T}(B^{\N-1})}+\|K(q^{n}){\cal T}^{n}\n
q^{n}\|_{L^{1}_{T}(B^{\N-1})}\\[2,5mm]
&\hspace{11,3cm}+\|{\cal T}^{n}\n q^{n}\|_{L^{1}_{T}(B^{\N-1})},\\[2,5mm]
&\hspace{3,7cm}\lesssim
T\|q^{n}\|^{2}_{L^{\infty}_{T}(B^{\N})}+\|{\cal
T}^{n}\|_{L^{2}_{T}(B^{\N-1})}\|q^{n}\|_{L^{2}_{T}(B^{\N+1})}
(1+\|q^{n}\|_{L^{\infty}_{T}(B^{\N})})\;.\\[4mm]
&\|[\frac{P_{1}(\rho^{n})}{\rho^{n}}]\nabla\theta^{n}\|_{L^{1}_{T}(B^{\N-1})}\leq
C(\|q^{n}\|_{L^{\infty}_{T}(B^{\N})}+1)\|{\cal
T}^{n}\|_{L^{1}_{T}(B^{\N})},\\[4mm]
&\|{\rm div}(\frac{\mu(\rho^{n})}{1+q^{n}}\n
u^{0})\|_{L^{1}_{T}(B^{\N-1})}\lesssim(1+
\|q^{n}\|_{L^{\infty}_{T}(B^{\N})})\|u^{0}\|_{L^{1}_{T}(B^{\N+1})}.\\
\end{aligned}
$$
We proceed similarly with the other terms:
$$-\D u^{0},\;\;\;\n\biggl(\frac{\zeta(\rho^{n})}{1+q^{n}}{\rm
div}(u^{0})\biggl),\;\;\;\n(K(\rho^{n})\D q^{0}).$$
Let us estimate now $\|H_{n}\|_{L^{1}_{T}(B^{\N-2})}$. We obtain:
$$
\begin{aligned}
\|\n(\frac{1}{1+q^{n}}).\n\theta^{n}\chi(\rho^{n})\|_{L^{1}_{T}(B^{\N-2})}&\lesssim\|K(q^{n})\n(\frac{q^{n}}{1+q^{n}}).
\n\theta^{n}\|_{L^{1}_{T}(B^{N/2-2})}\\[2mm]
&\;\;\;\;\;\;\;\;\;\;\;\;\;\;\;\;\hspace{2cm}+C\|\n(\frac{q^{n}}{1+q^{n}}).\n\theta^{n}\|_{L^{1}_{T}(B^{\N-2})},\\[3mm]
&\lesssim(1+\|q^{n}\|_{L^{2}_{T}(B^{\N})})\|q^{n}\|_{L^{2}_{T}(B^{\N+1})}\|{\cal T}^{n}\|_{L^{2}_{T}(B^{\N-1})}\;.\\
\end{aligned}
$$
We have after these last two terms:
$$
\begin{aligned}
\|\frac{T^{n}P_{1}(\rho^{n})}{\rho^{n}}{\rm
div}u^{n}\|_{L^{1}_{T}(B^{\N-2})}\lesssim&\|K(q^{n}){\rm
div}u^{n}\|_{L^{1}_{T}(B^{\N-2})}+\|K(q^{n})K_{1}({\cal T}^{n}){\rm
div}u^{n}\|_{L_{T}^{1}(B^{\N-2})}\\[2mm]
&\hspace{4,8cm}+\|K_{1}({\cal T}^{n}){\rm
div}u^{n}\|_{L^{1}_{T}(B^{\N-2})},
\end{aligned}
$$
with $K$ and $K_{1}$ regular in sense of the proposition
\ref{fcomposition} and:
$$
\begin{aligned}
&\|K(q^{n}){\rm div}u^{n}\|_{L^{1}_{T}(B^{\N-2})}\lesssim
T\|q^{n}\|_{L^{\infty}_{T}(B^{\N})}\|u^{n}\|_{L^{\infty}_{T}(B^{\N-1})}\;,\\[2mm]
&\|K_{1}({\cal T}^{n}){\rm
div}u^{n}\|_{L^{1}_{T}(B^{N/2-2})}\lesssim\|{\cal
T}^{n}\|_{L^{\infty}_{T}(B^{\N-2})}\|u^{n}\|_{L^{1}_{T}(B^{\N+1})}\;,
\end{aligned}
$$
so finally:
$$
\begin{aligned}
\|\frac{{\cal T}^{n}P_{1}(\rho^{n})}{\rho^{n}}{\rm
div}u^{n}\|_{L^{1}_{T}(B^{\N-2})}\lesssim\;&\big(T\|q^{n}\|_{L^{\infty}_{T}(B^{\N})}+\|{\cal
T}^{n}\|_{L^{\infty}_{T}(B^{\N-2})}\big)\|u^{n}\|_{L_{T}^{\infty}(B^{\N-1})}\\
&+\|q^{n}\|_{L_{T}^{\infty}(B^{\N})}\|{\cal
T}^{n}\|_{L_{T}^{\infty}(B^{\N-2})}\|u^{n}\|_{L_{T}^{1}(B^{\N+1})}.
\end{aligned}
$$
and, since $N\geq3$:
$$
\begin{aligned}
&\|u^{n}.\n\theta^{n}\|_{L^{1}_{T}(B^{\N-2})}\leq
\|u^{n}\|_{L_{T}^{\infty}(B^{\N-1})}\|\theta^{n}\|_{L^{1}_{T}(B^{\N})},\\[2mm]
&\|\n u^{n}:\n u^{n}\|_{L^{1}_{T}(B^{\N-2})}\leq\|u^{n}\|^{2}_{L^{2}_{T}(B^{\N})}.\\
\end{aligned}
$$
We obtain in using (\ref{fi1}) and the different previous inequalities:
$$\|(\bar{q}_{n+1},\bar{u}_{n+1},\bar{{\cal T}}_{n+1})\|_{F_{T}}\big(1-C2\sqrt{\e}(A_{0}+\sqrt{\e})\big)\leq
C_{1}\big(\e(A_{0}+\sqrt{\e})^{2}+T(A_{0}+\sqrt{\e})\big).$$ In
taking $T$ and $\e$ small enough  we have $({\cal{P}}_{n+1})$, so we
have shown by induction that $(q^{n},u^{n},{\cal T}^{n})$ is bounded
in $F_{T}$.
\subsubsection*{Second Step: Convergence of the
sequence}
 We will show
that $(q^{n},u^{n},{\cal T}^{n})$ is a Cauchy sequence in the Banach
space $F_{T}$, hence converges to some
$(q,u,{\cal T})\in F_{T}$.\\
Let:
$$\delta q^{n}=q^{n+1}-q^{n},\;\delta u^{n}=u^{n+1}-u^{n},\;\delta
{\cal T}^{n}={\cal T}^{n+1}-{\cal T}^{n}.$$
The system verified by $(\de q^{n},\de u^{n},\de {\cal T}^{n})$ reads:
$$
\begin{cases}
\begin{aligned}
&\p_{t}\delta q^{n}+{\rm div}\delta u^{n}=F_{n}-F_{n-1},\\[2mm]
&\p_{t}\delta u^{n}-{\rm
div}\biggl(\frac{\mu(\rho^{n})}{\rho^{n}}\n\delta
u^{n}\biggl)-\n\biggl(\frac{\zeta(\rho^{n})}{\rho^{n}}{\rm
div}(\delta u^{n})\biggl)-\n(K(\rho^{n})\D \delta q^{n})=\\[2mm]
&\hspace{7,5cm}G_{n}-G_{n-1}+G^{'}_{n}-G^{'}_{n-1},\\[2mm]
&\p_{t}\delta{\cal T}^{n}-{\rm
div}\biggl(\frac{\chi(\rho^{n})}{1+q^{n}}\n\delta{\cal
T}^{n}\biggl)=
H_{n}-H_{n-1}+H_{n}^{'}-H_{n-1}^{'},\\[2mm]
&\delta q^{n}(0)=0\;,\;\delta u^{n}(0)=0\;,\;\delta {\cal T}^{n}(0)=0,\\
\end{aligned}
\end{cases}
$$
where we define:
$$
\begin{aligned}
G^{'}_{n}=-{\rm
div}\big((\frac{\mu(\rho^{n+1})}{\rho^{n+1}}-\frac{\mu(\rho^{n})}{\rho^{n}})
\n &u^{n+1}\big)-\n\big((K(\rho^{n+1})-K(\rho^{n}))\D
 q^{n+1}\big)\hspace{5cm}\\
&\hspace{0,6cm}-\n\big((\frac{\zeta(\rho^{n+1})}{\rho^{n+1}}
-\frac{\zeta(\rho^{n})}{\rho^{n}}){\rm
div}(u^{n+1})\big).\\
\end{aligned}
$$
In the same way we have:
$$H^{'}_{n}={\rm
div}\big((\frac{\chi(\rho^{n+1})}{1+q^{n+1}}-\frac{\chi(\rho^{n})}{1+q^{n}})\n\theta^{n+1}\big).$$
Applying propositions \ref{fchaleur3}, \ref{flinear3}, and using $({\cal{P}}_{n})$, we get:
$$
\begin{aligned}
\|(\de q^{n},\de u^{n},\de{\cal T}^{n})\|_{F_{T}}\leq\;&
C(\|F_{n}-F_{n-1}\|_{L^{1}_{T}(B^{N/2})}+\|G_{n}-G_{n-1}+G^{'}_{n}-G^{'}_{n-1}\|_{L^{1}_{T}(B^{N/2-1})}\\[2mm]
&\hspace{3,7cm}+ \|H_{n}-H_{n-1}+H^{'}_{n}-H^{'}_{n-1}\|_{L^{1}_{T}(B^{N/2-2})}),\\
\end{aligned}
$$
And by the same type of estimates as before, we get:
$$\|(\de q^{n},\de u^{n},\de {\cal T}^{n})\|_{F_{T}}\leq
C\sqrt{\e}(1+A_{0})^{3}\|(\de q^{n-1},\de u^{n-1},\de {\cal
T}^{n-1})\|_{F_{T}}.$$ So in taking $\e$ enough small we have that
$(q^{n},u^{n},{\cal T}^{n})$ is Cauchy sequence, so the limit
$(q,u,{\cal T})$ is in $F_{T}$ and we verify easily that this is a
solution of the system.
\subsubsection*{Third step: Uniqueness}
Suppose that $(q_{1},u_{1},{\cal T}_{1})$ and $(q_{2},u_{2},{\cal
T}_{2})$ are solutions with the same initial conditions, and
$(q_{1},u_{1},{\cal T}_{1})$ corresponds to the previous
solution.\\
Assume moreover that we have:
$$\|q_{1}(t)\|_{L^{\infty}}\leq\alpha,\;\;\forall t\in[0,T].$$
We set then:
$$\de q=q_{2}-q_{1},\;\de u=u_{2}-u_{1},\;\de
{\cal T}={\cal T}_{2}-{\cal T}_{1}.$$ The triplet $(\de q,\de u,\de
{\cal T})$ satisfies the following system:
$$
\begin{cases}
\begin{aligned}
&\p_{t}\delta q+{\rm div}\delta
u=F_{2}-F_{1},\\[3mm]
&\p_{t}\delta u-{\rm div}\big(\frac{\mu(\rho_{2})}{\rho^{2}}\n\delta
u\big)-\n\big(\frac{\zeta(\rho_{2})}{\rho_{2}}{\rm
div}(\delta u)\big)-\n\big(K(\rho_{2})\D \delta q\big)=G_{2}-G_{1}+G^{'}\\[2mm]
&\p_{t}\delta{\cal T}-{\rm
div}(\frac{\chi(\rho_{2})}{1+q_{2}}\n\delta{\cal T})=
H_{2}-H_{1}+H^{'},\\[2mm]
&\delta q(0)=0\;,\;\delta u(0)=0\;,\;\delta {\cal T}(0)=0\\
\end{aligned}
\end{cases}
$$
with:
$$
\begin{aligned}
G^{'}=&-{\rm
div}\biggl(\big(\frac{\mu(\rho_{2})}{\rho_{2}}-\frac{\mu(\rho_{1})}{\rho_{1}}\big)
\n u_{2}\biggl)-\n\big((K(\rho_{2})-K(\rho_{1}))\D
 q_{2}\big)\\[2mm]
&\hspace{4,2cm}-\n\biggl(\big(\frac{\zeta(\rho_{2})}{\rho_{2}}
-\frac{\zeta(\rho_{1})}{\rho_{1}}\big){\rm
div}(u_{1})\biggl),\\[2mm]
H^{'}=&{\rm
div}\biggl(\big(\frac{\chi(\rho_{2})}{1+q_{2}}-\frac{\chi(\rho_{1})}{1+q_{1}}\big)\n\theta_{2}\biggl)\;.
\end{aligned}
$$
Let $\bar{T}$ the largest time such that:
$\|q_{2}\|_{L^{\infty}((0,\bar{T})\times\R^{N}}\leq\alpha$. As
$q_{2}\in
C([0,T];B^{N/2})$, we have by continuity $0<\bar{T}\leq T$.\\
We are going to work on the interval $[0,T_{1}]$ with
$0<T_{1}\leq\bar{T}$ and we use the proposition \ref{flinear3}, so we
obtain in using the same type of estimates than in
the part on the contraction:\\
$$\|(\de q,\de u,\de{\cal T})\|_{\widetilde{F}^{\N}_{T}}\leq Z(T)\|(\de
q,\de u,\de{\cal T})\|_{\widetilde{F}^{\N}_{T}}$$ with
$Z(T)\rightarrow_{T\rightarrow0}0$.\\
We have then for $T_{1}$ small enough: $(\de q,\de u,\de{\cal
T})=(0,0,0)$ on $[0,T_{1}]$ and by connectivity we finally
conclude that:
$$q_{1}=q_{2},\;u_{1}=u_{2},\;{\cal T}_{1}={\cal
T}_{2}\;\;\mbox{on}\;\;[0,T].$$ \hfill {$\Box$}
\\
{\bf Proof of the theorem \ref{ftheo5}}\\
\\
In the special case $N=2$, we need to take more regular initial data
for the same reasons as in theorem \ref{ftheo2}. Indeed some terms
like $\Psi(\theta){\rm div}u$ or $u^{*}.\n\theta$ can't be
controlled
without more regularity.\\
The proof is similar to the previous proof of  theorem \ref{ftheo4}
except that we have changed the functional space $F_{T}(2)$, in
which the fixed point theorem is going to be applied. As we explain
above we can use the paraproduct because we have more regularity, so
we just see the term $u^{*}.\n\theta$. The other terms and the
details are
left to the reader.\\
We then have:
$$\|u^{*}.\n\theta\|_{L^{1}_{T}(\widetilde{B}^{-1,-1+\e^{'}})}\lesssim
\|{\cal
T}\|_{L^{1}_{T}(\widetilde{B}^{0,1+\e^{'}})}\|u\|_{L^{\infty}_{T}(B^{0})}.$$
\hfill {$\Box$}
\subsection{Local  existence theorem  in the general case}
Now we suppose that all the coefficients depend on the temperature
and on the density, and that conditions $(C)$ and $(D)$ are satisfied with strict inequalities.\\
One of the problem in the general case is the control of the
$L^{\infty}$ norm of the temperature $\theta$ in  order to have
control on the non linear terms where the physical coefficients
appear. Indeed in the theorem of composition
we need to control the norm $L^{\infty}$.\\
So we must impose that $\theta_{0}$ is in $B^{\N}$ to hope a
$L^{\infty}$ control. And in consequence the others initial data
have to be also more regular.\\
\\
{\bf Proof of theorem \ref{ftheo6}:}\\
\\
We proceed exactly like in  theorem \ref{ftheo4} except that we ask
more regularity for the initial data. We define then:
$$q^{n}=q^{0}+\bar{q}^{n},\;u^{n}=u^{0}+\bar{u}^{n},\;
{\cal T}^{n}={\cal T}^{0}+\bar{{\cal
T}}^{n}\;\;\mbox{and}\;\;\theta^{n}={\cal T}^{n}+\bar{\cal T}$$
where
 $(q^{0},u^{0},{\cal T}^{0})$ stands for the  solution of:
$$
\begin{cases}
\p_{t}q^{0}-\D q^{0}=0,\\
\p_{t}u^{0}-\D u^{0}=0,\\
\p_{t}{\cal T}^{0}-\D {\cal T}^{0}=0,\\
\end{cases}
$$
supplemented with initial data:
$$q^{0}(0)=q_{0}\;,\;u^{0}(0)=u_{0}\;,\;\theta^{0}(0)=\theta_{0}.$$
Let $(\bar{q}_{n},\bar{u}_{n},\bar{\theta}_{n})$ be the solution of
the following system:
$$
\begin{cases}
\begin{aligned}
&\p_{t}\bar{q}^{n+1}+{\rm div}(\bar{u}^{n+1})=F_{n},\\[3mm]
& \p_{t}\bar{u}^{n+1}-{\rm
div}\biggl(\frac{\mu(\rho^{n},\theta^{n})}{\rho^{n}}\n
\bar{u}^{n+1}\biggl)-\n\biggl(\frac{\zeta(\rho^{n},\theta^{n})}{\rho^{n}}
{\rm
div}(\bar{u}^{n+1})\biggl)-\n(\kappa(\rho^{n})\D\bar{q}^{n+1})=G_{n},\\[3mm]
&\p_{t}\bar{{\cal T}}^{n+1}-{\rm
div}\biggl(\frac{\chi(\rho^{n},\theta^{n})}{1+q^{n}}\bar{{\cal
T}}^{n+1}\biggl)
=H_{n},\\[3mm]
&(\bar{q}^{n+1},\bar{u}^{n+1},\bar{{\cal T}}^{n+1})_{t=0}=(0,0,0).\\
\end{aligned}
\end{cases}
$$
where:
$$
\begin{aligned}
F_{n}=&-{\rm div}(q^{n}u^{n})-\D q^{0}-{\rm div}(u^{0}),\\[4mm]
G_{n}=&-(u^{n})^{*}.\n
u^{n}+\n\big(\frac{K^{'}_{\rho^{n}}}{2}|\n\rho^{n}|^{2}
\big)-\n\big(\frac{\mu(\rho^{n},\theta^{n})}{\rho^{n}}\big){\rm
div}u^{n}+\n\big(\frac{\zeta(\rho^{n},\theta^{n})}{\rho^{n}}\big){\rm
div}u^{n}\\[3mm]
&+\frac{\lambda^{'}_{1}(\rho^{n},\theta^{n})\n\rho^{n}{\rm
div}u^{n}}{1+q^{n}}+\frac{(du^{n}+\n
u^{n})\,\mu_{1}^{'}(\rho^{n},\theta^{n})\n\rho^{n}}{1+q^{n}}\\[3mm]
&+\frac{\lambda_{2}^{'}(\rho^{n},\theta^{n})\n\theta^{n}{\rm
div}u^{n}}{1+q^{n}}+\frac{(du^{n}+\n
u^{n})\,\mu_{2}^{'}(\rho^{n},\theta^{n})\n\theta^{n}}{1+q^{n}}\\[3mm]
&+\frac{[P_{0}^{'}(\rho^{n})+\w^{n}P_{1}^{'}(\rho^{n})]\n
q^{n}}{1+q^{n}}+\biggl[\frac{P_{1}(\rho^{n})}{\rho^{n}\psi^{'}(T^{n})}\biggl]\nabla\theta^{n}\\[3mm]
&-\D u^{0}+{\rm
div}\biggl(\frac{\mu(\rho^{n},\theta^{n})}{1+q^{n}}\n
u^{0}\biggl)+\n\biggl(\frac{\mu(\rho^{n})+\lambda(\rho^{n})}{1+q^{n}}{\rm
div}(u^{0})\biggl)+\n(\kappa(\rho^{n})\D q^{0}),\\[4mm]
H_{n}=&\n(\frac{1}{1+q^{n}}).\n\theta^{n}\chi(\rho^{n},\theta^{n})-\frac{\w^{n}P_{1}(\rho^{n})}{\rho^{n}}{\rm
div}u^{n} -(u^{n})^{*}.\n\theta^{n}+\frac{D^{n}:\nabla
u^{n}}{\rho^{n}}\\[3mm]
&\,-\D\theta_{0}+{\rm
div}\big(\frac{\chi(\rho^{n},\theta^{n})}{1+q^{n}}\n\theta^{0}\big).\\
\end{aligned}
$$
\subsubsection*{1) First Step , Uniform Bound}
Let $\e$ be a small positive parameter and  choose $T$ small enough
so that in using the estimate of the proposition \ref{fchaleur}  we
have:
$$
\begin{aligned}
&\|{\cal
T}^{0}\|_{L^{1}_{T}(B^{\N+2})}+\|u^{0}\|_{L^{1}_{T}(B^{\N+2})}+\|q^{0}\|_{L^{1}_{T}(\widetilde{B}^{\N+2,\N+3})}
\leq\varepsilon,\\
&\|{\cal
T}^{0}\|_{\widetilde{L}^{\infty}_{T}(B^{\N})}+\|u^{0}\|_{\widetilde{L}^{\infty}_{T}(B^{\N})}+\|q^{0}\|_{\widetilde{L}^{\infty}_{T}
(\widetilde{B}^{\N,\N+1})}\leq
A_{0}.\\
\end{aligned}
\leqno{({\cal{H}}_{\e})}
$$
After we are going to show by induction that:
$$\|(\bar{q}_{n},\bar{u}_{n},\bar{{\cal T}}_{n})\|_{F_{T}}\leq\e.\leqno{({\cal{P}}_{n})}$$
As $(\bar{q}_{0},\bar{u}_{0},\bar{{\cal T}}_{0})=(0,0,0)$ the result
is true for $({\cal{P}}_{0})$. We suppose now $({\cal{P}}_{n})$ true
and we are going to show $({\cal{P}}_{n+1})$.
\\
To begin with we are going to show that $1+q^{n}$ is  positive. In
using the fact that $B^{\N}\h L^{\infty}$ and that we can take $\e$
enough small, we have for $t\in[0,T]$:
$$
\begin{aligned}
\|q^{n}-q_{0}\|_{L^{\infty}((0,T)\times\R^{N})}&\leq\;
C_{1}\big(\|{\rm div}(\bar{u}^{n})\|_{L^{1}_{T}(B^{\N})}+
\|{\rm div}(q^{n-1}u^{n-1})\|_{L^{1}_{T}(B^{\N})}\\
&\hspace{7cm}+\|{\rm div}(u^{0})\|_{L^{1}_{T}(B^{\N})}\big),\\
&\leq \;C_{1}\big(2\sqrt{T}\e+\|q^{n-1}u^{n-1}\|_{L^{1}_{T}(B^{\N+1})}\big),\\
\end{aligned}
$$
and by induction hypothesis $\big({\cal P}_{n-1}\big)$:
$$
\begin{aligned}
\|q^{n-1}u^{n-1}\|_{L^{1}_{T}(B^{\N+1})}\leq&\sqrt{T}\|q^{n-1}\|_{L^{\infty}_{T}(B^{\N})}
\|u^{n-1}\|_{L^{2}_{T}(B^{\N+1})}\\[2mm]
&\hspace{4cm}+\sqrt{T}\|q^{n-1}\|_{L^{2}_{T}(B^{\N+1})}\|u^{n-1}\|_{L^{\infty}_{T}(B^{\N})},\\
\end{aligned}
$$
thus:
$$\|q^{n}-q_{0}\|_{L^{\infty}((0,T)\times\R^{N})}\leq C_{1}\sqrt{T}(2\e+(A_{0}+\e)\e).$$
Finally we have:
$$
\begin{aligned}
&\|1+q_{0}\|_{L^{\infty}((0,T)\times\R^{N})}-\|q^{n}-q_{0}\|_{L^{\infty}((0,T)\times\R^{N})}\leq
1+q^{n}\leq\|1+q_{0}\|_{L^{\infty}((0,T)\times\R^{N})}\\
&\hspace{10,3cm}
+\|q^{n}-q_{0}\|_{L^{\infty}((0,T)\times\R^{N})},\\
\end{aligned}
$$
$$
\frac{c}{2\bar{\rho}}\leq1+q^{n}\leq1+\frac{\|\rho_{0}\|_{L^{\infty}}}{\bar{\rho}}\;.
$$
So we have shown that:
$$\|q^{n}\|_{L^{\infty}}\leq 2A_{0}\leqno{(*)}$$
and that $\rho^{n}$ is bounded away from $0$.
\\
\\
To verify the uniform bound we use the propositions \ref{fchaleur3}
and \ref{flinear3}. For that we have to  verify the different
hypotheses of these propositions, so that we study the following
terms:
$$a^{n}=\frac{\mu(\rho^{n},\theta^{n})}{1+q^{n}}\;,\;b^{n}=\frac{\zeta(\rho^{n},\theta^{n})}{1+q^{n}}\;,\;c^{n}=K(\rho_{n})\;,\;
d^{n}=\frac{\chi(\rho^{n},\theta^{n})}{1+q^{n}}\cdotp$$\\
In using $({\cal{P}}_{n})$ and by continuity of $\mu$ and the fact
that $\mu$ is positive on
$[\bar{\rho}(1+\min(q^{0}))-\alpha,\bar{\rho}(1+\max(q^{0}))+\alpha]\times
[\bar{\theta}(1+\min({\cal T}^{0}))
-\alpha,\bar{\rho}(1+\max({\cal T}^{0}))+\alpha]$, we have:
$$0<c_{1}\leq a^{n}=\frac{\mu(\rho^{n},\theta^{n})}{1+q^{n}}\leq M_{1}\;.$$
We proceed similarly to verify the bounds of the other terms.
\\
After we use the proposition \ref{fproduit} and the fact that $q^{n}$
is bounded. We get:
$$
\begin{aligned}
\|\n a^{n}\|_{\widetilde{L}_{T}^{2}(B^{\N})}&\lesssim
\|q^{n}\|_{\widetilde{L}^{2}_{T}(B^{\N+1,\N+2})}+
\|{\cal T}^{n}\|_{\widetilde{L}_{T}^{2}(B^{\N+1})}\;,\\[3mm]
\|\n b^{n}\|_{\widetilde{L}_{T}^{2}(B^{\N})}&\lesssim
\|q^{n}\|_{\widetilde{L}_{T}^{2}(B^{\N+1,\N+2})}+
\|{\cal T}^{n}\|_{\widetilde{L}_{T}^{2}(B^{\N+1})}\;.\\[3mm]
\|\n c^{n}\|_{\widetilde{L}^{2}(B^{\N})}&\lesssim
\|q^{n}\|_{\widetilde{L}_{T}^{2}(\widetilde{B}^{\N+1,\N+2})}.\\[2mm]
\end{aligned}
$$
Next we want to estimate $\p_{t}c^{n}$ in $L^{1}_{T}(B^{\N})$. For
that, we use the fact that:
$$\p_{t}c^{n}=K^{'}(\rho^{n})\p_{t}q^{n}=K^{'}(\rho^{n})({\rm
div}(q^{n-1}u^{n-1})-{\rm div}(u^{n}))$$
And we have:
$$
\begin{aligned}
&\|K^{'}(\rho^{n})({\rm div}(q^{n-1}u^{n-1})-{\rm
div}(u^{n}))\|_{L^{1}_{T}(B^{\N})}\leq\|K^{'}(\rho^{n}){\rm
div}(q^{n-1}u^{n-1})\|_{L^{1}_{T}(B^{\N})}\\
&\hspace{10cm}+\|K^{'}(\rho^{n}){\rm div}(u^{n})\|_{L^{1}_{T}(B^{\N})},\\[3mm]
&\lesssim
(1+\|q^{n}\|_{L^{\infty}_{T}(B^{\N})})(\|u^{n}\|_{L^{1}_{T}(B^{\N+1})}+\|q^{n-1}u^{n-1}\|_{L^{1}_{T}(B^{\N+1})}),\\[3mm]
&\lesssim(1+\|q^{n}\|_{L^{\infty}_{T}(B^{\N})})(\|u^{n}\|_{L^{1}_{T}(B^{\N+1})}+
\|q^{n-1}\|_{L^{\infty}_{T}(B^{\N})}\|u^{n-1}\|_{L^{1}_{T}(B^{\N+1})}\\[2mm]
&\hspace{8,5cm}+\|q^{n-1}\|_{L^{2}_{T}(B^{\N+1})}\|u^{n-1}\|_{L^{2}_{T}(B^{\N})})\;.\\
\end{aligned}
$$
Now we want to show $({\cal{P}}_{n+1})$ by induction and in this
goal we will apply the estimates of proposition \ref{fchaleur3} and
proposition \ref{flinear3}. This is possible as we have verified above
the validity of the hypotheses. We obtain:
\begin{equation}
\begin{aligned}
&\|(\bar{q}^{n+1},\bar{u}^{n+1},\bar{{\cal
T}}^{n+1})\|_{F_{T}}\big(1-C(\|a^{n}\|_{L^{2}_{T}(B^{\N+1})}
+\|b^{n}\|_{L^{2}_{T}(B^{\N+1})}+\|c^{n}\|_{L^{2}_{T}(B^{\N+1})}\\[2mm]
&\hspace{0,8cm}+\|d^{n}\|_{L^{2}_{T}(B^{\N})}+\|\p_{t}
c^{n}\|_{L^{1}_{T}(B^{\N})})\big)\lesssim\|(\n
F_{n},G_{n})\|_{L^{1}_{T}(B^{\N})}+\|H_{n}\|_{L^{1}_{T}(B^{\N})}\;.\\
\end{aligned}
\label{fB1}
\end{equation}
We want to control now the part on the right-hand side of
(\ref{fB1}), for this we do like previously in using proposition
\ref{fproduit}. We have:
$$\|F_{n}\|_{L^{1}_{T}(B^{\N+1})}\leq\|{\rm
div}(q^{n}u^{n})\|_{L^{1}_{T}(B^{\N+1})}+\|{\rm
div}(u^{0})\|_{L^{1}_{T}(B^{\N+1})}+\|\D
q^{0}\|_{L^{1}_{T}(B^{\N+1})},$$
with:
$$
\begin{aligned}
\|u^{n}q^{n}\|_{L^{1}_{T}(B^{\N+2})}&\leq\|q^{n}\|_{L^{\infty}_{T}(L^{\infty})}\|u^{n}\|_{L^{1}_{T}(B^{\N+2})}
+\|q^{n}\|_{L^{2}_{T}(B^{\N+2})}\|u^{n}\|_{L^{2}_{T}(B^{\N})}\\
&\leq\|q^{n}\|_{L^{\infty}_{T}(\widetilde{B}^{\N,\N+1})}\|u^{n}\|_{L^{1}_{T}(B^{\N+2})}
+\sqrt{T}\|q^{n}\|_{\widetilde{L}^{2}_{T}(\widetilde{B}^{\N+1,\N+2})}\|u^{n}\|_{L^{\infty}_{T}(B^{\N})}.\\
\end{aligned}
$$
One ends up with:
$$\|F_{n}\|_{L^{1}_{T}(B^{\N+1})}\leq C(A_{0}+\e+\sqrt{\e})^{2}.$$
Next we want to control the different terms of $G_{n}$. We have:
$$
\begin{aligned}
&\|(u^{n})^{*}.\n
u^{n}\|_{L^{1}_{T}(B^{\N})}\leq\sqrt{T}\|u^{n}\|_{L^{\infty}_{T}(B^{\N})}\|u^{n}\|_{L^{2}_{T}(B^{\N+1})},\\[2mm]
&\|\n\biggl(\frac{K^{'}_{\rho^{n}}}{2}|\n\rho_{n}|^{2}
\biggl)\|_{L^{1}_{T}(B^{\N})}\lesssim\sqrt{T}\|q^{n}\|_{L^{\infty}_{T}(B^{\N+1})})\|q^{n}\|_{L^{2}_{T}(B^{\N+2})}.\\[2mm]
\end{aligned}
$$
After we have:
$$
\|\n(\frac{\mu(\rho^{n},\theta^{n})}{\rho^{n}}){\rm
div}u_{n}\|_{L^{1}_{T}(B^{\N})}\leq
\|u^{n}\|_{L^{2}_{T}(B^{\N+1})}(\sqrt{T}\|q^{n}\|_{L^{\infty}_{T}(B^{\N+1})}+
\|{\cal T}^{n}\|_{L^{2}_{T}(B^{\N+1})})\;.
$$
We treat similarly the term:
$$\n\biggl(\frac{\zeta(\rho^{n},\theta^{n})}{\rho^{n}}\biggl){\rm
div}u^{n}.$$
Next we study the term:
$$
\begin{aligned}
\|\frac{\lambda^{'}(\rho^{n},\theta^{n})\n\rho^{n}\rm
div(u^{n})}{1+q^{n}}\|_{L^{1}_{T}(B^{\N})}
&\lesssim(1+\|q^{n}\|_{L^{\infty}_{T}(B^{\N})})\sqrt{T}\|u^{n}\|_{L^{2}_{T}(B^{\N+1})}
\|q^{n}\|_{L^{\infty}_{T}(B^{\N+1})}\;.\\
\end{aligned}
$$
We proceed similarly for the following term:
$$
\begin{aligned}
&\frac{(du^{n}+\n
u^{n})\mu^{'}(\rho^{n},\theta^{n})\n\rho^{n}}{1+q^{n}}\,.\\
\end{aligned}
$$
Next we have:
$$
\begin{aligned}
&\|\frac{[P_{0}^{'}(\rho^{n})+\w^{n}P_{1}^{'}(\rho^{n})]\n
q^{n}}{1+q^{n}}\|_{L^{1}_{T}(B^{\N})}\lesssim
T\|q^{n}\|_{L^{\infty}_{T}(B^{\N})}\|q^{n}\|_{L^{\infty}_{T}(B^{\N+1})}\\[2mm]
&\hspace{5,5cm}+T\|{\cal
T}^{n}\|_{L^{\infty}_{T}(B^{\N})}\|q^{n}\|_{L^{\infty}_{T}(B^{\N+1})}
\big(1+\|q^{n}\|_{L^{\infty}_{T}(B^{\N})}\big)\;,\\[2mm]
&\|[\frac{P_{1}(\rho^{n})}{\rho^{n}\Psi^{'}(\w^{n})}]\nabla\theta^{n}\|_{L^{1}_{T}(B^{\N})}\lesssim
\big(\sqrt{T}(\|q^{n}\|_{L^{\infty}_{T}(B^{\N})}+\|{\cal T}^{n}\|_{L^{\infty}_{T}(B^{\N})})\\[2mm]
&\hspace{6,3cm}+\sqrt{T} \|{\cal
T}^{n}\|_{L^{\infty}_{T}(B^{\N})}\|q^{n}\|_{L^{\infty}_{T}(B^{\N})}\big)\|{\cal
T}^{n}\|_{L^{2}_{T}(B^{\N})}\;,\\[2mm]
&\|{\rm div}(\frac{\mu(\rho^{n},\theta^{n})}{1+q^{n}}\n
u^{0})\|_{L^{1}_{T}(B^{\N})}\lesssim
\big(\|q^{n}\|_{L^{\infty}_{T}(B^{\N})}+\|{\cal T}^{n}\|_{L^{\infty}_{T}(B^{\N})}\big)\|u^{0}\|_{L^{1}_{T}(B^{\N+2})}\\[2mm]
&\hspace{5cm}+\|u^{0}\|_{L^{2}_{T}(B^{\N+1})}\big(\sqrt{T}\|q^{n}\|_{L^{\infty}_{T}(B^{\N+1})}+\|{\cal T}^{n}\|
_{L^{2}_{T}(B^{\N+1})}\big)\;.\\
\end{aligned}
$$
We proceed similarly with the other terms:
$$-\D u^{0},\;\;\;\n\biggl(\frac{\zeta(\rho^{n},\theta^{n})}{1+q^{n}}{\rm
div}(u^{0})\biggl),\;\;\;\n(\kappa(\rho^{n})\D q^{0}).$$
After we want to estimate the term $\|H^{n}\|_{L^{1}_{T}(B^{\N})}$.
So we
have:
$$
\begin{aligned}
\|\n(\frac{1}{1+q^{n}}).\n\theta^{n}\chi(\rho^{n},\theta^{n}))\|_{L^{1}_{T}(B^{\N})}\lesssim(1+
&\|q^{n}\|_{L^{\infty}_{T}(B^{\N})}+\|{\cal T}^{n}\|_{L^{\infty}_{T}(B^{\N})})\\[2mm]
&\times\sqrt{T}\|q^{n}\|_{L^{\infty}_{T}(B^{\N+1})}\|{\cal T}^{n}\|_{L^{2}_{T}(B^{\N+1})}.\\
\end{aligned}
$$
We have after these last terms:
$$
\begin{aligned}
\|\frac{\w^{n}P_{1}(\rho^{n})}{\rho^{n}}{\rm
div}u^{n}\|_{L^{1}_{T}(B^{\N})}\lesssim\;\sqrt{T}\,\big(\|q^{n}\|_{L^{\infty}_{T}(B^{\N})}+\|{\cal
T}^{n}\|_{L^{\infty}_{T}(B^{\N})}(1+&\|q^{n}\|_{L^{\infty}_{T}(B^{\N})})\big)\\[2mm]
&\times\|u^{n}\|_{L^{2}_{T}(B^{\N+1})}.\\
\end{aligned}
$$
and:
$$
\begin{aligned}
&\|u^{n}.\n\theta^{n}\|_{L^{1}_{T}(B^{\N})}\leq
\sqrt{T}\|u^{n}\|_{L^{\infty}_{T}(B^{\N})}\|{\cal T}^{n}\|_{L^{2}_{T}(B^{\N+1})}\\[2mm]
&\|\n u^{n}:\n u^{n}\|_{L^{1}_{T}(B^{\N})}\leq \|u^{n}\|^{2}_{L^{2}_{T}(B^{\N})}\\
\end{aligned}
$$
we obtain in using (\ref{fB1}), the hypothesis of recurrence to the state $n$ and the previous inequalities:\\
$$\|(\bar{q}_{n+1},\bar{u}_{n+1},\bar{{\cal T}}_{n+1})\|_{F_{T}}(1-C2\sqrt{\e}(A_{0}+\sqrt{\e})\leq
C_{1}(\e(A_{0}+\sqrt{\e})^{2}+T(A_{0}+\sqrt{\e})).$$
In taking $T$ and $\e$ small enough  we have $({\cal{P}}_{n+1})$, so
$(q^{n},u^{n},{\cal T}^{n})$ is bounded in $F_{T}$. To conclude we
proceed like in the proof of theorem \ref{ftheo4} and we show in the
same way that $(\bar{q}^{n},\bar{u}^{n},\bar{{\cal T}}^{n})$ is a
Cauchy sequence in $F_{T}$, hence converges to some $(q,u,{\cal T})$
in $F_{T}$. We verify after that $(\rho,u,\theta)$ is a solution of
the system.
\subsubsection*{Uniqueness:}
We compare the difference between two solutions with the same
initial data and we use essentially the same type of estimates than
in the part on contraction. The details are left to the reader.
\section{Appendix}
This part consists in one commutator lemma which enables us to
conclude in proposition \ref{flinear3}. Moreover we give the proof of
 proposition \ref{fcomposition} on the composition of function in
 hybrid spaces adapted from Bahouri-Chemin in \cite{fBC}.
\begin{lemme}
\label{fcommutateur} Let $0\leq s\leq1$. Suppose that
$A\in\widetilde{L}^{2}_{T}(B^{\N+1})$ and
$B\in\widetilde{L}^{2}_{T}(B^{\N-1+s})$. Then we have the following
result
:\\
$$\|\p_{k}[A,\D _{l}]B\|_{L^{1}_{T}(L^{2})}\leq
Cc_{l}2^{-l(\N-1+s)}\|A\|_{\widetilde{L}^{2}_{T}(B^{\N+1})}\|B\|_{\widetilde{L}^{2}_{T}(B^{\N-1+s})}$$
\\
with $\sum_{l\in\mathbb{Z}}c_{l}=1$.
\end{lemme}
{\bf Proof:}\\
\\
We have  the  following decomposition:
$$uv=T_{u}v+T^{'}_{v}u$$
where: $T_{u}v=\sum_{l\in\mathbb{Z}}S_{l-1}u\D_{l}v$
and: $T^{'}_{v}u=\sum_{l\in\mathbb{Z}}S_{l+2}v\D_{l}u$.\\
\\
We then have:
\begin{equation}
\p_{k}[A,\D_{l}]B=\p_{k}T^{'}_{\D_{l}B}A-\p_{k}\D_{l}T^{'}_{B}A+[T_{A},\D_{l}]\p_{k}B+T_{\p_{k}A}\D_{l}B
-\D_{l}T_{\p_{k}A}B.\label{f100}
\end{equation}
From now on, we will denote by $(c_{l})_{l\in\mathbb{Z}}$ a sequence
such that:
$$\sum_{l\in\mathbb{Z}}c_{l}\leq1.$$
\\
Now we are going to treat each term of (\ref{f100}). According to the
properties of quasi-orthogonality and the definition of
$T^{'}$ we have:\\
$$\p_{k}T^{'}_{\D_{l}B}A=\sum_{m\geq
l-2}\p_{k}(S_{m+2}\D_{l}B\D_{m}A).$$ \\
Next, in using Bernstein
inequalities, we have:
\\
$$
\begin{aligned}
\|\p_{k}T^{'}_{\D_{l}B}A\|_{L^{1}_{T}(L^{2})}&\lesssim\sum_{m\geq
l-2}2^{m}\|\D_{l}B\|_{L^{2}_{T}(L^{\infty})}\|\D_{m}A\|_{L^{2}_{T}(L^{2})}\\[2mm]
&\lesssim 2^{l\N}\|\D_{l}B\|_{L^{2}_{T}(L^{2})}\sum_{m\geq
l-2}2^{-m\N}(2^{m(\N+1)}\|\D_{m}A\|_{L^{2}_{T}(L^{2})})\\[2mm]
&\lesssim
2^{-l(N/2-1+s)}(2^{l(\N-1+s)}\|\D_{l}B\|_{L^{2}_{T}(L^{2})})\sum_{m\geq
l-2}(2^{m(\N+1)}\|\D_{m}A\|_{L^{2}_{T}(L^{2})})\\[2mm]
&\lesssim c_{l}2^{-l(N/2-1+s)}\|B\|_{\widetilde{L}^{2}_{T}(B^{\N-1+s})}\|A\|_{\widetilde{L}^{2}_{T}(B^{\N+1})}.\\
\end{aligned}
$$
\\
Next, we will use  the classic estimates on the paraproduct to bound
the second term of the right-hand side of (\ref{f100}). We obtain
then:
$$\|T^{'}_{B}A\|_{L^{1}_{T}(B^{\N+s})}\lesssim\|B\|_{L^{2}_{T}(B^{\N-1+s})}\|A\|_{L^{2}_{T}(B^{\N+1})}.$$
\\
After in using the spectral localization we have:
$$
\begin{aligned}
\|\p_{k}\D_{l}T^{'}_{B}A\|_{L^{1}_{T}(L^{2})}&\lesssim
2^{l}\|\D_{l}T^{'}_{B}A\|_{L^{1}_{T}(L^{2})}\\
&\lesssim c_{l}2^{-l(\N-1+s)}\|B\|_{\widetilde{L}^{2}_{T}(B^{\N-1+s})}\|A\|_{\widetilde{L}^{2}_{T}(B^{\N+1})}.\\
\end{aligned}
$$
\\
According to the properties of orthogonality of Littlewood-Payley decomposition we have:
$$[T_{A},\D_{l}]\p_{k}B=\sum_{|m-l|\leq
4}[S_{m-1}A,\D_{l}]\D_{m}\p_{k}B\;.$$
\\
In applying Taylor formula, we obtain for $x\in\R^{N}$:\\
$$[S_{m-1}A,\D_{l}]\D_{m}\p_{k}B(x)=2^{-l}\g\int^{1}_{0}h(y)(y.S_{m-1}\n
A(x-2^{-l}\tau y))\D_{m}\p_{k}B(x-2^{-l}y)d\tau dy\;.$$
\\
By an inequality of convolution we have:\\
$$\|[S_{m-1}A,\D_{l}]\D_{m}\p_{k}B\|_{L^{2}}\lesssim
2^{-l}\|\n A\|_{L^{\infty}}\|\D_{m}\p_{k}B\|_{L^{2}}\;.$$ So we get:
$$\|[T_{A},\D_{l}]\p_{k}B\|_{L^{1}_{T}(L^{2})}\lesssim c_{l}2^{-l(\N-1+s)}\|\n A\|_{L^{2}_{T}(L^{\infty})}
\|B\|_{\widetilde{L}^{2}_{T}(B^{\N-1+s})}\;.$$ Finally we have:
$$T_{\p_{k}A}\D_{l}B=\sum_{|l-m|\leq 4}S_{m-1}\p_{k}A\D_{l}\D_{m}B.$$
\\
Hence:
$$\|T_{\p_{k}A}\D_{l}B\|_{L^{1}_{T}(L^{2})}\leq\|\p_{k}A\|_{L^{2}_{T}(L^{\infty})}\|\D_{l}B\|_{L^{2}_{T}(L^{2})}\;.$$
\\
And the classic estimates on the paraproduct give:
$$\|T_{\p_{k}A}\D_{l}B\|_{L^{1}_{T}(B^{\N-1})}\lesssim c_{l}2^{-l(\N-1+s)}\|\p_{k}A\|_{\widetilde{L}^{2}_{T}(B^{\N})}
\|B\|_{\widetilde{L}^{2}_{T}(B^{\N-1+s})}\;.$$ The proof is
complete. \hfill {$\Box$}
\\
\\
{\bf Proof of proposition \ref{fcomposition}:}
\\
\\
To show (i) we use ``first lin\'earisation'' method
introduced by Y.Meyer in \cite{fM}, which amounts to write that:
$$
F(u_{1},u_{2},\cdots,u_{d})=\sum_{p\in\mathbb{Z}}
(F(S_{p+1}u_{1},\cdots,S_{p+1}u_{d})-F(S_{p}u_{1},\cdots,S_{p}u_{d})).
$$
According to Taylor formula, we have:
$$F(S_{p+1}u_{1},\cdots,S_{p+1}u_{d})-F(S_{p}u_{1},\cdots,S_{p}u_{d})=m^{1}_{p}u_{1}^{p}+\cdots+m^{d}_{p}u_{d}^{p}$$
with $u_{i}^{p}=\D_{p}u_{i}$ and
$$m^{i}_{p}=\int^{1}_{0}\p_{i}F(S_{p}u_{1}+su^{p}_{1},\cdots,S_{p}u_{i}+su^{p}_{i},\cdots,S_{p}u_{d}+su^{p}_{d})ds.$$
Observe that:
$$\|m^{i}_{p}\|_{L^{\infty}}\leq\|\n F\|_{L^{\infty}}.$$
We have:
$$\D_{p}F(u_{1},u_{2},\cdots,u_{d})=\D_{p}^{1}+\D_{p}^{2}$$
\\
where we have decomposed the sum into two parts:
$$\D_{p}^{(1)}=\sum_{q\geq
p}\biggl(\D_{p}(u_{1}^{q}m^{1}_{q})+\cdots+\D_{p}(u_{d}^{q}m^{1}_{q})\biggl),$$
$$\D_{p}^{(2)}=\sum_{q\leq
p-1}\biggl(\D_{p}(u_{1}^{q}m^{1}_{q})+\cdots+\D_{p}(u_{d}^{q}m^{1}_{q})\biggl).$$
Now we bound $\|\D_{p}^{(1)}\|_{L^{\rho}_{T}(L^{p})}$ in this way:
$$
\begin{aligned}
\|\D_{p}^{(1)}\|_{L^{\rho}_{T}(L^{p})}&\leq\sum_{q\geq
p}\big(\|u_{1}^{q}\|_{L^{\rho}_{T}(L^{p})}\|m^{1}_{q}\|_{L^{\infty}_{T}(L^{\infty})}+\cdots+
\|u_{d}^{q}\|_{L^{\rho}_{T}(L^{p})}\|m^{1}_{q}\|_{L^{\infty}_{T}(L^{\infty})}\big)\\[2mm]
&\leq\sum_{q\geq
p}\|\n F\|_{L^{\infty}}2^{-qs}c_{q}(\|u_{1}\|_{\widetilde{L}^{\rho}_{T}(B^{s}_{p})}+\cdots+\|u_{d}\|_{\widetilde{L}^{\rho}_{T}(B^{s}_{p})})\\
\end{aligned}
$$
with $(c_{q})\in l^{1}(\mathbb{Z})$.\\
Therefore, since $s>0$:
$$\sum_{p\in\mathbb{Z}}2^{ps}\|\D_{p}^{(1)}\|_{L^{\rho}_{T}(L^{p})}\leq
C\|\n
F\|_{L^{\infty}}(\|u_{1}\|_{\widetilde{L}^{\rho}_{T}(B^{s}_{p})}+\cdots+\|u_{d}\|_{\widetilde{L}^{\rho}_{T}(B^{s}_{p})}).$$
\\
To bound $\|\D_{p}^{(2)}\|_{L^{\rho}_{T}(L^{p})}$ we use the fact
that the support of the Fourier transform of $\D_{p}^{(2)}$ is
included in the shell $2^{p}\cal{C}$, so that according to Bernstein
inequality:\\
$$
\begin{aligned}
\|\D_{p}^{(2)}\|_{L^{\rho}_{T}(L^{p})}\leq&\sum_{q\leq
p-1}(\|\D_{p}(u_{1}^{q}m^{1}_{q})\|_{L^{\rho}_{T}(L^{p})}+\cdots+\|\D_{p}(u_{d}^{q}m^{1}_{q})\|_{L^{\rho}_{T}(L^{p})}),
\\[2mm]
\leq&C2^{-p([s]+1)}\sum_{q\leq
p-1}(\|\p^{[s]+1}(u_{1}^{q}m^{1}_{q})\|_{L^{\rho}_{T}(L^{p})}+\cdots+\|\p^{[s]+1}(u_{d}^{q}m^{1}_{q})\|_{L^{\rho}_{T}
(L^{p})}).
\end{aligned}
$$
\\
Moreover we have according to Fa\`a-di-Bruno formula:
$$\p^{k}m_{q}^{i}=\int^{1}_{0}\sum_{l_{1}+\cdots+l_{m}=k,l_{m}\ne
0}A^{k}_{l_{1}\cdots
l_{m}}F^{m+1}(S_{q}(u)+su^{q})\prod_{n=1}^{m}\p_{l_{n}}(S_{q}(u)+su^{q}))ds.$$
Hence we get for all $k\in\mathbb{N}$:
$$\|\p^{k}m^{i}_{q}\|_{L^{\infty}_{T}(L^{\infty})}\leq
C_{u_{i},k}2^{qk}$$
with: $C_{u_{i},k}=C(1+\|u_{i}\|_{L^{\infty}_{T}(L^{\infty})})$.\\
\\
We have then:
$$\|\D_{p}^{(2)}\|_{L^{\rho}_{T}(L^{p})}\leq
C2^{-p([s]+1)}\sum_{q\leq
p-1}c_{q}2^{q(-s+[s]+1)}C_{u_{1},\cdots,u_{d}}(\|u_{1}\|_{\widetilde{L}^{\rho}_{T}(B^{s}_{p})}+
\cdots+\|u_{d}\|_{\widetilde{L}^{\rho}_{T}(B^{s}_{p})}).$$
\\
Hence the result:
$$\sum_{p\in\mathbb{Z}}2^{ps}\|\D_{p}^{(2)}\|_{L^{\rho}_{T}(L^{p})}\leq
C_{u_{1},\cdots,u_{d}}(\|u_{1}\|_{\widetilde{L}^{\rho}_{T}(B^{s}_{p})}+\cdots+\|u_{d}\|_{\widetilde{L}^{\rho}_{T}(B^{s}_{p})}).$$
So  the first part of the proof is complete.
\\
\\
For proving (ii) we proceed in the same way as before. We get:\\
$$F(u)=\sum_{q\in\mathbb{Z}}m_{q}u_{q}.$$
And we have for $p>0$:
$$\D_{p}F(u)=\D_{p}^{1}+\D_{p}^{2}$$
so:
$$\|\D_{p}^{1}\|_{L_{T}^{\rho}(L^{2})}\leq\sum_{q\geq p}\|F^{'}\|
_{L^{\infty}}2^{-qs_{2}}c_{q}\|u\|_{\widetilde{L}^{\rho}_{T}(\widetilde{B}^{s_{1},s_{2}})}.$$
Hence in using convolution inequality:
$$\sum_{p>0}2^{ps_{2}}\|\D_{p}^{1}\|_{L_{T}^{\rho}(L^{2})}\leq C\|F^{'}\|_{L^{\infty}}
\|u\|_{\widetilde{L}_{T}^{\rho}(\widetilde{B}^{s_{1},s_{2}})}.$$
After we get for all $s>0$:
$$
\begin{aligned}
\|\D_{p}^{2}\|_{L_{T}^{\rho}(L^{2})}&\leq C2^{-p([s]+1)}\sum_{q\leq
p-1}\|\p^{[s]+1}(m_{q}u_{q})\|_{L_{T}^{\rho}(L^{2})},
\\[2mm]
&\leq C2^{-p([s]+1)}\sum_{q\leq p-1}2^{q([s]+1-s(q))}c_{q}\|u_{q}\|_{\widetilde{L}_{T}^{\rho}(\widetilde{B}^{s_{1},s_{2}})}.\\
\end{aligned}
$$
with $s(q)=s_{1}$ or $s_{2}$.\\
So we obtain:
\begin{equation}
\begin{aligned}
\sum_{p>0}2^{ps_{2}}\|\D_{p}^{2}\|_{L^{\rho}_{T}(L^{2})}\leq&
\sum_{p>0}2^{-p([s]+1-s_{2})}\sum_{q\leq 0}c_{q}2^{q
([s]+1-s_{1})}\|u\|_{\widetilde{L}^{\rho}_{T}(\widetilde{B}^{s_{1},s_{2}})}\\[2mm]
&+\sum_{p>0}2^{-p([s]+1-s_{2})}\sum_{0<q\leq p-1}c_{q}2^{q([s]+1-s_{2})}\|u\|_{\widetilde{L}_{T}^{\rho}(\widetilde{B}^{s_{1},s_{2}})}.\\
\end{aligned}
\label{fmoi}
\end{equation}
\\
We have to choose $s$, so for the first term of (\ref{fmoi}) we just
need that: $[s]+1-s_{2}>0$ and $[s]+1-s_{1}>0$ and for the second
term of (\ref{fmoi}) we just have a inequality of
convolution. So we can take $s=1+\max(s_{1},s_{2})$.\\
\\
We do the same for $p<0$ and we have:\\
$$
\begin{aligned}
\sum_{p\leq0}2^{ps_{1}}\|\D_{p}^{1}\|_{L^{\rho}_{T}(L^{2})}\leq&\sum_{p\leq0}2^{ps_{1}}\sum_{q\geq
p}
\|F^{'}\|_{L^{\infty}}2^{-qs_{1}}c_{q}\|u\|_{\widetilde{L}^{\rho}(\widetilde{B}^{s_{1},s_{2}})}\\
&\hspace{1cm}+\sum_{p\leq0}2^{ps_{1}}\sum_{p\leq q\leq0 }
\|F^{'}\|_{L^{\infty}}2^{-qs_{2}}c_{q}\|u\|_{\widetilde{L}^{\rho}(\widetilde{B}^{s_{1},s_{2}})}.\\
\end{aligned}
$$
We conclude by a inequality of convolution.\\
And for the term $\D_{p}^{2}$ we get:
$$
\begin{aligned}
\sum_{p\leq0}2^{ps_{1}}\|\D_{p}^{2}\|_{\widetilde{L}^{\rho}(L^{2})}\leq&\sum_{p\leq0}2^{-p([s]+1-s_{1})}
\sum_{q\leq p-1}c_{q}2^{q([s]+1-s_{1})}\|u\|_{\widetilde{L}^{\rho}(\widetilde{B}^{s_{1},s_{2}})}.\\
\end{aligned}
$$
\\
\\
For proving (iii) and (iv), one just has to use the following identity:
$$G(v)-G(u)=(v-u)\int_{0}^{1}H(u+\tau(v-u))d\tau+G^{'}(0)(v-u)$$
where $H(w)=G^{'}(w)-G^{'}(0)$, and we conclude by using (i), (ii)
and  proposition \ref{fproduit}. \hfill {$\Box$}
\\
\\
\section{Annex: Notations of differential calculus}
If $f\,:\,\R^{n}\rightarrow\R\,,$ we denote:\\
$$df=\sum_{i=1}^{n}\frac{\p f}{\p x_{i}}dx_{i}=\p_{i}fdx_{i}$$
with the summation convention on repeated indices and the simplified notation:\\
$$\p_{i}f=\frac{\p f}{\p x_{i}}.$$
The vector field associated to the differential $df$ is noted $\n
f$,\\
$$\n f=\frac{\p f}{\p x_{i}}\frac{\p }{\p x_{i}}.$$
\\
Let $f:\,\R^{n}\rightarrow\R^{n}$. Let denote $f_{i}$ the $i$ th component of $f$, and:\\
$$(df)_{i,j}=\p_{j}f_{i}.$$
By analogy with the case of the scalar, we denote:
$$\n f=(df)^{*}\,,\;\mbox{so}\;(\n f)_{i,j}=\p_{i}f_{j}.$$
The ${\rm curl}$ of f is given by:
$$({\rm curl}f)_{i,j}=\p_{i}f_{j}-\p_{j}f_{i}.$$
The divergence of the vector field f is given by:
$${\rm div}\, f=tr\,df=\p_{i}f_{i}.$$
If $A\,:\,\R^{n}\rightarrow\R^{n\times n}$, with coefficients
$a_{i,j}$,we set:
$$({\rm div}\,A)_{j}=\p_{i}a_{i,j}\,,\;{\rm div}\,A={\rm
div}(A\frac{\p}{\p x_{j}})\,dx_{j}.$$ In particular, for $f$ scalar,
we
have:\\
$${\rm div}(fI)=df.$$
And finally we set:\\
$$A:B=a_{i,j}b_{i,j}\;.$$

\end{document}